\documentclass[a4paper,12pt]{amsart}
\usepackage{yfonts,amsmath,amsfonts,amssymb,amsthm}
\usepackage{graphics}
\usepackage{epsfig}
\usepackage{istgame}
\usepackage{array}
\usepackage{tikz}
\usepackage{float}
\usepackage[font=footnotesize,labelfont=bf]{caption} 
\usepackage[hidelinks,hyperindex,breaklinks]{hyperref}

\usepackage[normalem]{ulem}

\newtheorem{theorem}{Theorem}
\newtheorem{remark}{Remark}

\usepackage{color}
\definecolor{darkred}{rgb}{0.9,0.1,0.1}
\definecolor{darkblue}{rgb}{0,0,0.7}
\definecolor{darkgreen}{rgb}{0,0.5,0}

\title[Malliavin calculus in regularity structures]{Lecture notes on Malliavin calculus in regularity structures}

\author[L.~Broux]{Lucas Broux}
\author[F.~Otto]{Felix Otto}
\author[M.~Tempelmayr]{Markus Tempelmayr}

\newcommand{\treeZero}{
\hspace{-0.8ex}
\raisebox{0.8mm}{
\begin{istgame}
\istroot(0) \endist
\end{istgame}}
}

\newcommand{\treeThree}{
\hspace{-.75ex}
\raisebox{-.5mm}{
\begin{istgame}[font=\tiny]
\xtdistance{3mm}{2mm}
\setistgrowdirection{north}
\istroot(0)[null node] \istb* \istb* \istb* \endist
\end{istgame}}
}

\newcommand{\treeThreeZero}{
\hspace{-.75ex}
\raisebox{-.5mm}{
\begin{istgame}[font=\tiny]
\xtdistance{3mm}{2mm}
\setistgrowdirection{north}
\istroot(0)[null node] \istb* \istb* \endist
\end{istgame}}
}

\newcommand{\treeThreeZeroZero}{
\hspace{-1ex}
\raisebox{-.5mm}{
\begin{istgame}[font=\tiny]
\xtdistance{3mm}{2mm}
\setistgrowdirection{north}
\istroot(0)[null node] \istb* \endist
\end{istgame}}
}

\newcommand{\treeThreeThree}{
\hspace{-.5ex}
\raisebox{-.75mm}{
\begin{istgame}[font=\tiny]
\xtdistance{2mm}{1.5mm}
\setistgrowdirection{north}
\istroot(0)[null node] \istb* \istb \istb* \endist
\istroot(a)(0-2)[null node] \istb* \istb* \istb* \endist
\end{istgame}}
}

\newcommand{\treeThreeThreeZeroA}{
\hspace{-.5ex}
\raisebox{-.75mm}{
\begin{istgame}[font=\tiny]
\xtdistance{2mm}{1.5mm}
\setistgrowdirection{north}
\istroot(0)[null node] \istb* \istb \istb* \endist
\xtdistance{2mm}{3mm}
\istroot(a)(0-2)[null node] \istb* \istb* \endist
\end{istgame}}
}

\newcommand{\treeThreeThreeZeroB}{
\hspace{-.5ex}
\raisebox{-.75mm}{
\begin{istgame}[font=\tiny]
\xtdistance{2mm}{3mm}
\setistgrowdirection{north}
\istroot(0)[null node] \istb \istb* \endist
\xtdistance{2mm}{1.5mm}
\istroot(a)(0-1)[null node] \istb* \istb* \istb* \endist
\end{istgame}}
}

\newenvironment{keyword}
{
\begin{center}
\begin{minipage}{.9\textwidth}\small\textbf{Keywords}:\noindent
}
{
\end{minipage}
\end{center}
}

\newenvironment{msc}
{
\begin{center}
\begin{minipage}{.9\textwidth}\small\textbf{MSC 2020}:\noindent
}
{
\end{minipage}
\end{center}
}

\begin{document}

\begin{abstract}
Malliavin calculus provides a characterization of the centered
model in regularity structures that is stable under removing the
small-scale cut-off. In conjunction with a spectral gap inequality,
it yields the stochastic estimates of the model.

\smallskip

This becomes transparent on the level of a notion of model that 
parameterizes the solution manifold, and thus is indexed by
multi-indices rather than trees, and which allows for
a more geometric than combinatorial perspective.
In these lecture notes, this is carried out for a PDE with heat operator,
a cubic nonlinearity, and driven by additive noise, reminiscent of the stochastic
quantization of the Euclidean $\phi^4$ model.

\smallskip

More precisely, we informally motivate our notion of the model $(\Pi,\Gamma)$ as
charts and transition maps, respectively, of the nonlinear solution manifold. 
These geometric objects are algebrized in terms of formal power series, 
and their algebra automorphisms.
We will assimilate the directional Malliavin derivative to a tangent vector of
the solution manifold. This means that it can be treated as a modelled distribution,
thereby connecting stochastic model estimates to pathwise solution theory, 
with its analytic tools of reconstruction and integration.
We unroll an inductive calculus that in an automated way applies to 
the full subcritical range.
\end{abstract}

\maketitle


\begin{msc} 
60H17, 
60L30, 
60H07, 
81T16. 
\end{msc}


\section{Motivation and setting} \label{section:motivation}

\subsection{A nonlinear partial differential equation\texorpdfstring{ for $\phi$}{} with rough right-hand side\texorpdfstring{ $\xi$}{}}\label{ss:1}

We focus on the parabolic differential operator of second order (in fact, the heat operator)
in $d$ space dimensions
\begin{align*}
L:=\partial_0-\sum_{i=1}^d\partial_i^2,
\end{align*}
where $\partial_i$ denotes the partial derivative\footnote{there would only be minor changes for other constant-coefficient elliptic or parabolic operators} w.~r.~t.~$x_i$.
Hence $x_0$ is the time-like variable and $\{x_i\}_{i=1}^d$ are the space-like variables\footnote{in fact, we treat the parabolic operator like an elliptic one}.
Given a parameter $\lambda$ and a space-time function $\xi$
we are interested in the manifold of all space-time functions $\phi$ 
that solve the partial differential equation (PDE) in the entire space-time
\begin{align}\label{ck18}
L\phi-\lambda \phi^3=\xi\quad\mbox{on}\;\mathbb{R}^{1+d},
\end{align}
which is nonlinear due to the presence of the cube\footnote{there would be few changes for another power}.
At this point, let us emphasize that our use of the word ``manifold'' throughout these notes is informal.
In particular, we will not attempt to rigorously endow the space of solutions to \eqref{ck18} with the structure of a topological manifold.

\medskip

We are interested in the situation when the right-hand side (r.~h.~s.)
$\xi$ is so rough that it is not a function but just a Schwartz distribution. 
A Schwartz distribution $\xi$ is a bounded linear form on the space of Schwartz functions.
The space of Schwartz functions in turn is the linear space of all infinitely often differentiable
functions $\zeta$ that decay so fast that the family of semi-norms
\begin{align}\label{ck52a}
\sup_{x\in\mathbb{R}^{1+d}}(|x|^k+1)|\partial^{\bf n}\zeta(x)|\quad\mbox{is finite},
\end{align}
where $k\in\mathbb{N}_0$, ${\bf n}=(n_0,\cdots,n_d)\in\mathbb{N}_0^{1+d}$, 
and $\partial^{\bf n}$ $:=\prod_{i=0}^d\partial_i^{n_i}$.
The pairing is denoted by $(\xi,\zeta)$. 

\medskip

A pertinent example for $d=0$ is the following: Take a realization of Brownian motion, 
which we think of as a function $B_{x_0}$ of our time-like variable $x_0$, and consider
\begin{align}\label{ck45}
(\xi,\zeta)=-\int_{\mathbb{R}}dx_0 B_{x_0}\partial_0\zeta.
\end{align}
Since by the law of the iterated logarithm (i.~e.~$|B_{x_0}|\lesssim|x_0|^\alpha+1$ for any 
$\alpha>\frac{1}{2}$) we have $|(\xi,\zeta)|$ 
$\lesssim\sup_{x}(|x|^2+1)|\partial_0\zeta(x)|$, $\xi$ is indeed a Schwartz distribution.
Informally, i.~e.~in a distributional sense, one writes (\ref{ck45})
as $\xi=\partial_0B$. The derivative $\xi$ of Brownian motion is called (temporal) white noise.
Note that $\phi=B$ satisfies (\ref{ck18}) with $\lambda=0$ (next to $d=0$).
Since almost surely, $B$ has infinite variation, 
$\xi$ cannot be represented as an integral against a locally integrable function,
and thus is a genuine Schwartz distribution.

\medskip

We are interested in the even worse situation when solutions\footnote{the
notation is consistent with Subsection \ref{ss:para}} $\Pi_0$ of the
corresponding linear problem, i.~e.~(\ref{ck18}) with $\lambda=0$,
\begin{align}\label{ck48}
L\Pi_0=\xi
\end{align}
are genuine Schwartz distributions.
If $\Pi_0$ and $\xi$ are distributions, (\ref{ck48}) is to be interpreted in the sense of
\begin{align}\label{ck51}
(\Pi_0,L^*\zeta)=(\xi,\zeta)\quad\mbox{for all Schwartz functions $\zeta$},
\end{align}
where $L^*:=-\partial_0-\sum_{i=1}^d\partial_i^2$ is the (informal) dual of $L$.
Once more, $d=0$ provides an easy example: if $\xi=\frac{d^2B}{dx_0^2}$
then $\Pi_0=\frac{dB}{dx_0}$ modulo an additive random constant, and thus is a genuine 
distribution.

\medskip

However, if $\Pi_0$ is a genuine distribution, 
then its cube $\Pi_0^3$ does not have a canonical sense,
which is why the equation is called  ``singular'' in this regime.
This does not bode well for the non-linear problem (\ref{ck18})
and is the challenge addressed in these lecture notes.


\subsection{Structure of these lecture notes}

In Section \ref{section:motivation} we informally motivate and rigorously introduce
our version of a centered model in the language of regularity structures. 
In doing so, we adopt a more geometric than combinatorial perspective. 
In Subsection \ref{ss:renorm}, we postulate the form of a counterterm for (\ref{ck18}),
motivated by the symmetries from Subsection \ref{ss:intro}, 
giving rise to the index set of multi-indices and the notion of ``homogeneity''.
We then introduce the concept of a parameterization of the nonlinear solution
manifold (Subsection \ref{ss:para}), informally\footnote{term-by-term in the physics jargon} 
write it as a power series
(Subsection \ref{ss:new}) 
recovering the same index set of multi-indices as in Subsection \ref{ss:renorm},
and finally ``algebrize'' it in terms of a formal power series $\Pi$ (Subsection \ref{ss:series}),
with $\Pi^{-}$ and $c$ corresponding to the r.~h.~s.~and the counterterm, respectively. 

\medskip

In the following subsections, we rigorously characterize $(\Pi,\Pi^{-},c)$:
Only some of the coefficients are allowed to be non-zero, i.~e.~``populated'' 
(Subsection \ref{ss:pop}).
Returning to the scale invariance of the solution manifold,
we 
impose a scaling invariance on the coefficients of $(\Pi,\Pi^{-},c)$ 
(Subsection \ref{ss:hom}).
Having restricted to the singular but renormalizable range (Subsection \ref{ss:intro}),
and as a consequence of a Liouville principle, $\Pi$ is unique (Subsection \ref{ss:unique}),
and $c$ and thus $\Pi^{-}$ are unique (Subsection \ref{ss:c});
the latter connects to what is called BPHZ renormalization.

\medskip

However, by imposing the scale invariance, we arbitrarily singled out an origin;
we now consider an arbitrary ``base-point'' $x$.
This gives rise to another parameterization $\Pi_x$ (Subsection \ref{ss:trans}),
and thus to a change-of-base-point transformation
$\Gamma_{x}^*$ (Subsection \ref{ss:trafo}), which is algebrized as an
endomorphism of the algebra of formal power series. The following subsections deal 
with the structure of $\Gamma_x^*$ and its pre-dual $\Gamma_x$: 
its uniqueness (Subsection \ref{ss:oneGamma}), its action on 
space-time polynomials (Subsection \ref{ss:dual}), its matrix representation
(Subsection \ref{ss:coord}), the population of its matrix entries (Subsection \ref{ss:triang}),
and its triangularity (Subsection \ref{ss:triang2}).
All this amounts to a self-contained introduction of a centered model 
$(\Pi_x,\Pi_x^{-},\Gamma_{xx'}:=\Gamma_x\Gamma_{x'}^{-1})$ in the sense of regularity structures.

\medskip

In Section \ref{section:main}, we state the stochastic 
estimates on $(\Pi_x,\Pi_x^{-},\Gamma_{x})$ and sketch
their proof, focussing here on the algebraic aspects.
The scale invariance in law that emerges in the limit 
of vanishing regularization
motivates the uniform estimates (Subsection \ref{ss:expect}).
The main result and some extensions are formulated in Subsections \ref{ss:main}
and \ref{ss:improve}.
The motivation for the usage of the (directional) Malliavin derivative $\delta\Pi$ 
is given in Subsection \ref{ss:usage}. Its control requires a further structural insight, 
arising from the (informal) parameterization of the tangent space to the solution manifold,
see Subsection \ref{ss:tangent}.
This motivates to approximate $\delta\Pi$, locally near $x$, 
by a linear combination of $\Pi_x$, with coefficients encoded in 
a linear endomorphism ${\rm d}\Gamma_x^*$. The latter can be assimilated 
with a modelled distribution in the language of regularity structures.
The next subsections are devoted to the structure of ${\rm d}\Gamma_x^*$:
its uniqueness (Subsection \ref{ss:struct}), the population of its matrix entries
(Subsection \ref{ss:order1}), and its triangularity (Subsection \ref{ss:order}),
which determines the order of induction in the proof.
While the original relation $\Pi\mapsto\Pi^{-}$ is not robust
under vanishing regularization, its counterpart
$\delta\Pi\mapsto\delta\Pi^{-}$ on the level of Malliavin derivatives is 
(Subsection \ref{ss:robust}). 

\medskip

In the following subsections we embark on the actual (stochastic) estimates.
While the construction and the estimates have to be logically intertwined, in these notes we focus on the quantitative estimates under the assumptions that the objects have been constructed.
We refer to \cite{LOTT21} for the full arguments in the context of a quasi-linear equation.
In Subsection \ref{ss:spectral_gap},
we introduce our use of the spectral gap inequality by duality, 
estimating probabilistic $L^p$-norms. The carr\'e-du-champs is
inherently linked to the space-time $L^2$ topology; this is best propagated 
by working with $L^2$-based space-time Besov norms (Subsection \ref{ss:Besov}). 
We then lay out the induction step, which is a sequence of
four arguments: a continuity property\footnote{reminiscent of modelled distributions} 
of $x\mapsto {\rm d}\Gamma_x^*$, 
namely an estimate of ${\rm d}\Gamma_{x+y}^*-{\rm d}\Gamma_x^*\Gamma_{x\,x+y}^*$, 
by an {\sc algebraic argument} (Subsection \ref{ss:algIII}), an estimate of the
``rough-path increment'' $\delta\Pi^{-}-{\rm d}\Gamma_x^*\Pi_x^{-}$ by
what in regularity structures corresponds to a {\sc reconstruction} of $\delta\Pi^{-}$
(Subsection \ref{ss:reconstr}), 
an estimate of $\delta\Pi-{\rm d}\Gamma_x^*\Pi_x$ by Schauder theory
which in regularity structures is called {\sc integration} 
(Subsection \ref{ss:intIII}),
and returning to the continuity property of $x\mapsto {\rm d}\Gamma_x^*$ by
an analytic argument we call {\sc three-point argument} (Subsection \ref{ss:3pIII}).

\medskip

This is the crucial but only the first of three rounds of these four
arguments, as explained in Subsection \ref{ss:logic}: What was done for
$(\delta\Pi$ $-{\rm d}\Gamma_x^*\Pi_x,$ 
$\delta\Pi^{-}-{\rm d}\Gamma_x^*\Pi_x^{-},$ 
${\rm d}\Gamma_{x+y}^*-{\rm d}\Gamma_x^*\Gamma_{x\,x+y}^*)$ needs to be repeated
for $(\delta\Pi,\delta\Pi^{-},{\rm d}\Gamma_x^*)$ (Subsection \ref{ss:second_round}), 
and finally for $(\Pi,\Pi^{-},\Gamma_x^*)$ itself (Subsection \ref{ss:third_round}). 
The arguments in the second and third round, which have to carried out within
the induction step in the right order, are simpler. 

\medskip

In Section \ref{s:proofs}, we provide the analytical details of the proof. 
A family of convolution operators\footnote{that provide local averaging} 
that has the semi-group property (Subsection \ref{ss:semigroup}) is convenient, 
as it allows for telescoping over dyadic space-time scales in reconstruction 
(Subsection \ref{ss:recIIIrev}).
It is also convenient 
when estimating the expectation $\mathbb{E}\Pi^{-}$ 
(Subsection \ref{ss:expectation}). Any general Schwartz convolution kernel can be recovered
(Subsection \ref{ss:change}).
We use an annealed version of Sobolev's inequality to pass from an estimate of $\delta\Pi^{-}-{\rm d}\Gamma_x^*\Pi_x^{-}$
to one of $\delta\Pi^{-}$ (Subsection \ref{ss:smearing}).
The last three subsections deal with integration, 
where the semi-group convolution now provides a decomposition of $L^{-1}$ 
into small and large scales, which is quintessential for any Schauder theory. 
Subsection~\ref{ss:abstract_int} provides an abstract representation of solutions to $Lu=f$ under appropriate growth conditions, 
which is applied to $\Pi$ and $\delta\Pi$ in Subsection~\ref{ss:intIrevisited} 
and to $\delta\Pi-{\rm d}\Gamma^*_x\Pi_x$ in Subsection~\ref{ss:intIIIrevisited}.

%
\tableofcontents
%

\subsection{Bibliographical context}

The theory of regularity structures \cite{Hai14} provides a systematic framework to treat equations in a singular regime as outlined in Subsection~\ref{ss:1}. 
Inspired by rough path theory \cite{Lyo98, Gub04}, 
it separates analytic from probabilistic arguments, 
the former being dealt with in \cite{Hai14} while the latter is addressed in \cite{CH16}.
In addition it comes with an algebraic structure \cite{BHZ19, BCCH21},
which allows to treat equations arbitrarily close to criticality\footnote{
see Subsection~\ref{ss:intro} for the notion of (sub)criticality}.
For an introduction to regularity structures 
we recommend \cite{Hai15}; 
more focus on the algebraic aspects is put in \cite{Che22b}. 
While the above-mentioned works provide a local well-posedness theory for a large class of singular stochastic PDEs, 
this setting can also be leveraged to establish global well-posedness results \cite{CMW23} 
and to prove properties of the associated invariant measure \cite{HS22}.
Some of the most recent developments include progress on the stochastic quantization 
of the two and three dimensional Yang-Mills measure \cite{CCHS22a, CCHS22b, CS23}, 
for an overview see \cite{Che22a}.

\medskip

Like \cite{OST23}, these lecture notes are intended as a reader's digest of \cite{LOTT21}.
The present notes give a more complete account of \cite{LOTT21} than the older notes \cite{OST23}, 
however adapted to a different model case,
namely the semi-linear $\phi^4$ equation rather than the previously treated
quasi-linear equation, thus confirming the flexibility of the approach, see also \cite{BL23, GT23}.
In the analytic treatment, these notes adopt a simplification from the recent \cite{HS23} (see the discussion around \eqref{ap25} below)
which in turn took inspiration from \cite{LOTT21}.

\medskip

Loosely speaking, \cite{LOTT21} constitutes an alternative to \cite{CH16}
when it comes to
establishing
the stochastic estimate of the centered model in regularity structures. 
The approach of \cite{LOTT21} is based on Malliavin calculus 
in conjunction with a spectral gap estimate, and was taken over in \cite{HS23}.
Following \cite{LOTT21}, but opposed to \cite{HS23},
these lecture notes implement this approach 
for a model that can be seen as a parameterization of the solution manifold, 
a top-down approach that leads to a more parsimonious index set than trees, namely multi-indices.
This type of model was introduced in \cite{OSSW},
and then motivated in \cite{LO22} in a non-singular setting.
Building upon \cite{OST23}, we highlight the conceptual use of Malliavin calculus,
which is brought to full fruition in \cite{Tem}.

\medskip

We note that prior to this line of research, Malliavin calculus had been used within regularity structures, but with a more classical purpose, namely the study of densities of solutions to stochastic partial differential equations, see \cite{CFG,GL20,Sch23}.
Stochastic estimates based on a spectral gap assumption were first carried out in \cite{IORT23}, 
however in a simple setting with no need of regularity structures. 
More recently, and inspired by \cite{LOTT21}, the spectral gap inequality has been adopted as a convenient tool to prove stochastic estimates in regularity structures: 
In the tree-setting but without appealing to diagrams in \cite{HS23, BH23}, 
in the tree-setting and making use of diagrammatic tools in \cite{BB23}, 
and in a rough path setting in \cite{GK23}. 

\medskip

The (pre-)Lie- and Hopf-algebraic aspects of the structure group of this multi-index based model
were first explored in \cite{LOT}, and embedded
into a post-Lie perspective in \cite{BK23,JZ}. 
In \cite{BD23} it was shown that the algebraic structure on multi-indices is also a multi-Novikov algebra, 
which is isomorphic to the free multi-Novikov algebra. 
Some algebraic aspects of renormalization in the multi-index setting are investigated in \cite{BL23}, 
and the analogue for rough paths is studied in \cite{Lin23}. 
These notes follow the hands-on and boiled-down approach of \cite{OST23} to the structure group. 

\medskip

We briefly comment on alternative solution theories to singular SPDEs.
Simultaneous to the development of regularity structures, 
another approach via paracontrolled distributions was presented in \cite{GIP15}. 
In the scope of this theory is equation \eqref{ck18} for $d=3$ and space-time white noise \cite{CC18}, 
for an introduction we refer to \cite{GP18}. 
Shortly after, yet another approach based on Wilson's renormalization group was given in 
\cite{Kup16} and applied to \eqref{ck18}, 
again for $d=3$ and space-time white noise. 
While both approaches are not capable to treat equations arbitrarily close to criticality, 
the latter one was more recently generalized to the full subcritical range \cite{Duc21}, 
based on the continuum version of the Polchinski flow equation.
An overview of the flow equation approach is given in \cite{Duc23}.
Incidentally, Malliavin calculus has been used in the paracontrolled setting to establish stochastic estimates and universality \cite{FG19b},
however not in combination with the spectral gap inequality. 

\medskip

Let us finally mention that the inductive approach presented here
has similarities with the one of Epstein-Glaser, see \cite[Section~3.1]{Sch95}.
In particular it does not suffer from the well-known difficulty of “overlapping sub-divergences” 
in Quantum Field Theory, which is also an issue in \cite{CH16}. 


\subsection{A random \texorpdfstring{$\xi$}{noise} with symmetries in law and restriction to the singular and subcritical range \texorpdfstring{$\alpha\in(-1,0)$}{}}\label{ss:intro}

In order to develop some theory for rough $\xi$, 
one approach is to randomize it;
i.~e.~to draw the space-time Schwartz distributions $\xi$ 
from a suitable ensemble/probability measure/law. 
One then seeks to capitalize on structural assumptions of the ensemble, 
namely the symmetries in law under
\begin{align}
\mbox{shift (``stationarity''):}
&\quad\xi(\cdot+x)=_{law}\xi\quad\mbox{for}\;x\in\mathbb{R}^{1+d},\label{ck19}\\
\mbox{(spatial) reflection symmetry:}
&\quad\xi(R_i\cdot)=_{law}\xi\quad\mbox{for}\;i=1,\cdots,d,\label{ck20}\\
\mbox{parity:}&\quad-\xi=_{law}\xi,\label{ck21}
\end{align}
where $R_i$ denotes the reflection at the $\{x_i=0\}$-plane.
These symmetries are valuable since they are compatible with the solution
manifold of (\ref{ck18}):
\begin{itemize}
\item Because $L$ has constant coefficients, $\phi(\cdot+x)$ solves (\ref{ck18})
with $\xi$ replaced by $\xi(\cdot+x)$; 
\item because $L$ is even in $\partial_i$ for $i=1,\cdots,d$, $\phi(R_i\cdot)$ solves (\ref{ck18})
with $\xi$ replaced by $\xi(R_i\cdot)$;
\item because the nonlinearity $\phi^3$ is odd in $\phi$, $-\phi$ solves (\ref{ck18}) with $\xi$ replaced by $-\xi$.
\end{itemize}

\medskip

One pertinent example is space-time white noise, which is a centered\footnote{i.~e.~of vanishing 
expectation} Gaussian on the space of Schwartz distributions, and as such characterized by
the covariance
\begin{align}\label{ck46}
\mathbb{E}(\xi,\zeta)(\xi,\zeta')\;=\;\int_{\mathbb{R}^{1+d}}dx\zeta\zeta'
\quad\mbox{for all Schwartz functions $\zeta,\zeta'$}.
\end{align}
Because the inner product $(\zeta,\zeta')$ $\mapsto\int dx\zeta\zeta'$ 
is invariant under shift and reflection,
white noise satisfies (\ref{ck19}) and (\ref{ck20}); 
it automatically satisfies (\ref{ck21}) as a centered Gaussian.

\medskip

Let us now address a further crucial symmetry in law, namely under scaling.
Recall that Brownian motion has the scale invariance
\begin{align*}
\{x_0\mapsto B_{r^2x_0}\}=_{law}\{x_0\mapsto rB_{x_0}\}\quad\mbox{for}\;r\in(0,\infty).
\end{align*}
By (\ref{ck45}) this translates to
\begin{align}\label{ck47}
&\mbox{in case of $d=0$:}\quad (\xi,r^{-2}\zeta(r^{-2}\cdot))=_{law}r^{-1}(\xi,\zeta)\nonumber\\
&\mbox{jointly in Schwartz functions $\zeta$},
\end{align}
which could also directly be inferred from (\ref{ck46}). The reason for expressing
this scale invariance in terms of $\zeta\mapsto r^{-2}\zeta(r^{-2}\cdot)$ is that
for $d=0$ this is the informal dual of $\xi\mapsto\xi(r^2\cdot)$, so that (\ref{ck47})
informally means
\begin{align}\label{ck46bis}
\mbox{in case of $d=0$:}\quad \xi(r^2\cdot)=_{law}r^{-1}\xi.
\end{align}
For $r\downarrow 0$, (\ref{ck46bis}) reflects that $\xi$ is not a function,
since zooming-in would increase its (typical) size.

\medskip

Not surprisingly, (\ref{ck46bis}) extends to $d>0$; however we need to package it 
in order to fit the parabolic $L$.
Thus we consider, for arbitrary $r\in(0,\infty)$, the scaling operator
\begin{align}\label{ck29}
Rx=(r^2x_0,rx_1,\cdots,rx_d),
\end{align}
which due to its anisotropy is compatible with $L$, 
or rather with $L^*$, in the sense that for any Schwartz function $\zeta$
\begin{align}\label{ck49}
L^*\{x\mapsto \zeta(R^{-1}x)\}=r^{-2}(L^*\zeta)(R^{-1}\cdot).
\end{align}
It follows from (\ref{ck46}) that (\ref{ck47}) extends to $d>0$:
jointly in the Schwartz function $\zeta$ we have
\begin{align}\label{ck50}
(\xi,{\rm det}R^{-1}\zeta(R^{-1}\cdot))&=_{law}\sqrt{{\rm det}R^{-1}}(\xi,\zeta)\nonumber\\
&=r^{-\frac{D}{2}}(\xi,\zeta)\quad\mbox{where}\;D:=2+d
\end{align}
denotes the ``effective dimension'' of our parabolic space-time.
In line with (\ref{ck46bis}), we informally rewrite this as
\begin{align*}
\xi(R\cdot)=_{law}r^{-\frac{D}{2}}\xi\quad\mbox{for all}\;r\in(0,\infty),
\end{align*}
which highlights that white noise gets rougher with increasing dimension.

\medskip

Actually,
we shall consider ensembles that have a scale invariance in law 
characterized by an exponent $s\in\mathbb{R}$ in
\begin{align}\label{ck22}
\xi(R\cdot)=_{law} r^{s-\frac{D}{2}}\xi\quad\mbox{for}\;r\in(0,\infty).
\end{align}
The exponent in (\ref{ck22}) is written such that the white noise ensemble 
satisfies (\ref{ck22}) with $s=0$.
It is very convenient to extend to $s\not=0$, as will
become clear in Subsection \ref{ss:unique}, see assumption (\ref{ck53}) below.
We will discuss an example of such an ensemble in Subsection~\ref{ss:expect} below.

\medskip

Consider now a random Schwartz distribution $\Pi_0$ that satisfies (\ref{ck48}) 
and has a scale invariance in law, which we informally write as
\begin{align}\label{lb01}
\Pi_0(R\cdot)=_{law}r^{\alpha}\Pi_0\quad\mbox{for all}\;r\in(0,\infty)
\end{align}
for some exponent $\alpha$. Working with (\ref{ck51}) and using
(\ref{ck49}), we learn that 
\eqref{ck22} translates into
\begin{align*}
\alpha=2+s-\frac{D}{2}.
\end{align*}
Hence for\footnote{the same holds in the borderline case of $d=2+2s$,
but is slightly more difficult to see} $d>2+2s$ and thus $D>4+2s$, we have 
\begin{align}\label{ao03}
\alpha<0, 
\end{align}
and $\Pi_0$ is expected to be a genuine
distribution.
We should not expect a solution $\phi$ of the nonlinear equation \eqref{ck18} to have better regularity in general, 
and hence also $\phi$ is expected to be a genuine distribution.
This is the situation we are interested in.

\medskip

We momentarily return to a discussion of the solution manifold of \eqref{ck18}.
Motivated by (\ref{ck22}), we consider the transformation $\xi\mapsto\hat\xi$:
\begin{align}\label{ck32}
\hat\xi:=r^{-(s-\frac{D}{2})}\xi(R\cdot),
\end{align}
which amounts to a ``blow-up'', or ``zoom-in'', for $r\ll 1$.
From (\ref{ck49}) we learn that (\ref{ck18}) 
is invariant under
\begin{align}\label{ap40}
\hat\phi&:=r^{-\alpha}\phi(R\cdot)\quad\mbox{where}\quad\alpha:=2+(s-\frac{D}{2}) ,
\end{align}
provided we adjust the strength of the {\it cubic} nonlinearity according to
\begin{align}
\hat\lambda&:=r^{3\times 2+(3-1)(s-\frac{D}{2})}\lambda
=r^{2(1+\alpha)}\lambda.\label{ck26}
\end{align}
The exponent $\alpha$ in (\ref{ap40}) generalizes (\ref{ck32}).
By invariance we mean an invariance of the solution manifold
in the sense that (\ref{ck18}) 
implies
$L\hat\phi=\hat\lambda\hat\phi^3 + \hat{\xi}$.
For our analysis, we have to limit ourselves to the ``(super-) renormalizable'' 
or ``subcritical'' case,
which means the effect of the nonlinearity vanishes on small scales,
as encapsulated by a positive exponent in (\ref{ck26}):
\begin{align}\label{ao04}
\alpha+1>0. 
\end{align}
Hence we restrict ourselves to the range of
\begin{align}\label{ck28}
\alpha
\in(-1,0)\quad\mbox{which by (\ref{ap40}) means}\quad s-\frac{D}{2}\in(-3,-2).
\end{align}
By (\ref{ck50}), this range for instance includes 
$d=4$ and $s\in(0,1)$, or $d=3$ and $s\in(-\frac{1}{2},\frac{1}{2})$.


\subsection{Renormalization through a counterterm \texorpdfstring{$h$}{}, multi-indices\texorpdfstring{ $\beta$}{}, and the homogeneity \texorpdfstring{$|\cdot|$}{}}\label{ss:renorm}

While white noise has the invariances (\ref{ck19}) - (\ref{ck21}),
and many more, it still does not allow to give (\ref{ck18}) a sense as such.
In fact, one needs to ``renormalize'' (\ref{ck18}), which means the following:
\begin{itemize}
\item On the one hand, one regularizes $\xi$ without affecting the invariances
(\ref{ck19}) - (\ref{ck21}).
\item On the other hand, one modifies the PDE (\ref{ck18}) by introducing a 
regularization-dependent ``counterterm'' that is postulated to be deterministic, 
i.~e.~independent of the realization of $\xi$ but dependent on the ensemble.
\end{itemize} 
For a given Schwartz function $\psi$,
we consider its parabolic rescaling $\psi_r$ to length scale $r$ 
\begin{align}\label{ck30}
\psi_r=r^{-D}\psi(R^{-1}\cdot)\quad\mbox{cf.~(\ref{ck29}), and set}\quad \Pi_r(x):=
(\Pi,\psi_r(x-\cdot))
\end{align}
for Schwartz distributions $\Pi$,
so that informally $\Pi_r$ is the convolution $\psi_r*\Pi$.

\medskip
Now fix a Schwartz function $\eta$ with $\int \eta = 1$, i.~e.~a kernel, and
consider the corresponding mollification $\{\xi_\rho = \eta_{\rho} * \xi\}_{\rho\downarrow 0}$ 
as regularization. 
Note that $\xi_\rho$ still satisfies \eqref{ck19} \& \eqref{ck21}, 
and provided $\eta$ is even in the spatial coordinates which we will henceforth assume, 
it also satisfies \eqref{ck20}.
The task is to determine the counterterm in such a way that 
\begin{itemize}
\item on the one hand, the new solution manifold converges for $\rho\downarrow0$
to a limiting manifold,
which is independent on the way of regularization (e.~g.~of $\eta$), 
\item and that on the other hand,
the new solution manifold preserves as many of the invariances (in law) of the old one as possible.
\end{itemize}

\medskip

To do so we make a general ansatz for the counterterm of the form
\begin{equation}\label{mt01}
L\phi = \lambda\phi^3+\xi_\rho
+\sum_\beta h_\beta \lambda^{\beta(3)} \prod_{{\bf n}} \Big(\tfrac{1}{{\bf n}!}\partial^{\bf n}\phi\Big)^{\beta({\bf n})},
\end{equation}
and successively reduce its degrees of freedom by suitable postulates. 
Here the sum is taken over all multi-indices $\beta$ over $\{3\}\sqcup\mathbb{N}_0^{1+d}$;
recall that a multi-index associates
to (the dummy\footnote{We choose to call it $3$ as it belongs to the cubic nonlinearity; if \eqref{ck18} had a further nonlinearity $\bar\lambda\phi^2$ we would choose the index set $\{2,3\}\sqcup\mathbb{N}_0^{1+d}$.}) $3$ and every
${\bf n}=(n_0,\dots,n_d)\in\mathbb{N}_0^{1+d}$ a non-negative integer which is non-zero
only for finitely many ${\bf n}$'s.
As usual we have set ${\bf n}!:=\prod_{i=0}^d(n_i!)$,
and if not specified otherwise, sums or products over 
${\bf n}$ extend over $\mathbb{N}_0^{1+d}$
and statements involving ${\bf n}$ are meant to hold for all such ${\bf n}$.
Furthermore, 
	\begin{align}
		\text{the $h_\beta$'s are deterministic.} \nonumber
	\end{align}
One should think of them as carrying the index $\rho$; in particular,
$h_{\beta}$ typically diverges as $\rho\downarrow 0$, but we omit this for brevity.

\medskip

To reduce the complexity of the counterterm we first note that the linear equation,
i.~e.~for $\lambda=0$, is not in need of renormalization. 
We thus postulate
\begin{equation}\label{mt03}
h_\beta=0\quad\text{unless}\quad\beta(3)>0.
\end{equation}

\medskip

We turn to scale invariance.
Note from 
\eqref{ap40}
that $\partial^{\bf n}\hat\phi = r^{|{\bf n}|-\alpha}\partial^{\bf n}\phi(R\cdot)$
provided we set
\begin{align}\label{ck30bis}
|{\bf n}|:=2n_0+n_1+\cdots+n_d.
\end{align}
Thus the scale invariance \eqref{ck32} - \eqref{ck26} carries over from \eqref{ck18} to \eqref{mt01} 
provided we set $\hat\rho=r^{-1}\rho$ and 
\begin{equation}\label{mt02}
\hat{h}_\beta = r^{2-|\beta|} h_\beta,
\end{equation}
where what we call the homogeneity $|\beta|$ of $\beta$ is defined through\footnote{note that since $|0|=\alpha < 0$, $| \beta |$ may be negative despite the notation $| \cdot |$} 
\begin{align}\label{ao07}
|\beta|-\alpha=\beta(3)2(1+\alpha)
+\sum_{{\bf n}}\beta({\bf n})(|{\bf n}|-\alpha).
\end{align}
Think now of $\hat\rho$ as being fixed, say $1$, 
so that we expect $\hat{h}_\beta$ to be finite.
Then taking $r\downarrow0$ amounts to $\rho\downarrow0$, 
i.~e.~removing the mollification from $\xi_\rho$. 
We then read off \eqref{mt02} that $h_\beta\to\infty$ for $|\beta|<2$, 
whereas $h_\beta\to0$ for $|\beta|>2$. 
Keeping only the relevant $h_\beta$'s, 
i.~e.~those diverging as $\rho\downarrow0$, we thus postulate 
\begin{equation}\label{mt05}
h_\beta=0 \quad \text{unless} \quad |\beta|<2.
\end{equation}

\medskip

We now turn to the invariances \eqref{ck19} - \eqref{ck21}.
By the shift invariance \eqref{ck19} of the law of $\xi$, 
and since the counterterm only depends on the ensemble but not on its realizations, 
\begin{equation*}
\text{the $h_\beta$'s are space-time constants.} 
\end{equation*}
Postulating that the spatial reflection invariance of \eqref{ck18} carries over to \eqref{mt01}, 
meaning that $\phi(R_i\cdot)$ solves \eqref{mt01} with $\xi_\rho(R_i\cdot)$ 
provided $\phi$ solves \eqref{mt01} with $\xi_\rho$,
we deduce from the spatial reflection invariance \eqref{ck20} 
\begin{equation}\label{mt06}
h_\beta=0 \quad\text{unless}\quad \sum_{\bf n} n_i \beta({\bf n}) \textnormal{ is even}.
\end{equation}
Similarly, by parity \eqref{ck21} we postulate 
\begin{equation}\label{mt07}
h_\beta=0 \quad\text{unless}\quad \sum_{\bf n}\beta({\bf n}) \textnormal{ is odd}.
\end{equation}

\medskip

We now put together \eqref{mt03}, \eqref{mt05}, \eqref{mt06}, and \eqref{mt07} 
in order to reduce the form of the counterterm made in the ansatz \eqref{mt01}.
Since $1+\alpha>0$ by \eqref{ao04}, 
the conditions $\beta(3)\geq1$ and $|\beta|<2$ yield by the definition \eqref{ao07} of the homogeneity $|\beta|$ that 
$-3\alpha>\sum_{\bf n}\beta({\bf n})(|{\bf n}|-\alpha)$.
On the one hand, this implies by $\alpha<0$ (cf.~\eqref{ao03}) 
that $\sum_{\bf n}\beta({\bf n})\leq2$;
Using that $\sum_{\bf n}\beta({\bf n})$ is odd we deduce $\sum_{\bf n}\beta({\bf n})=1$.
On the other hand, this also implies that $\beta({\bf n})=0$ unless $|{\bf n}|<2$.
Using furthermore that $\sum_{\bf n}n_i\beta({\bf n})$ is even we arrive at 
	\begin{align*}
		h_\beta=0 \quad\text{unless}\quad \beta = k \delta_3 + \delta_{\mathbf{0}} \text{ for some } k \geq 0,
	\end{align*}
where $\delta_{\bf 0}$ denotes the multi-index that associates the value one to ${\bf 0}\in\mathbb{N}_0^{1+d}$
and zero otherwise, and similarly for $\delta_3$.
Therefore, the counterterm in \eqref{mt01} reduces to $h^{(\rho)} \phi$,
where 
\begin{align}\label{mt04}
h^{(\rho)} = \sum_{k\geq0} c_k^{(\rho)} \lambda^k ;
\end{align}
As for $h_\beta$ we will omit from now on
the dependence of $c_k$ on $\rho$ in our notation. 
Note that $c_k$ coincides
with $h_{k\delta_3+\delta_{\bf 0}}$, 
and we thus have
\begin{align}
& \text{$c_k$ is a deterministic constant} , \text{ and} \label{mt08} \\
& c_k=0 \quad\text{unless}\quad | k\delta_3+\delta_{\bf 0}|<2 
\quad\textnormal{and}\quad 
k>0. \nonumber
\end{align}
In view of 
\eqref{ao07} the latter translates to:
\begin{align}\label{ck57}
c_k=0 \quad\text{unless}\quad
0<k<(1+\alpha)^{-1}.
\end{align}
These remaining (finitely many) degrees of freedom
are fixed in Subsection~\ref{ss:c} by the so-called BPHZ-choice of renormalization. 
Note that the number of constants increases as we approach 
the critical threshold $\alpha=-1$.
Thus, we have arrived at the renormalized equation
\begin{align}\label{ao40}
L\phi-(\lambda \phi^3+h^{(\rho)}\phi)=\xi_\rho\quad\mbox{in}\;\mathbb{R}^{1+d}.
\end{align}

\medskip

We remark that \eqref{ck18} has further symmetries that could be considered, 
e.~g.~invariance under space-like orthogonal transformations $Ox:=(x_0,\bar{O}(x_1,\dots,x_d))$ for orthogonal $\bar{O}\in\mathbb{R}^{d\times d}$:
If $\phi$ satisfies \eqref{ck18} with $\xi$, 
then $\phi(O\cdot)$ satisfies \eqref{ck18} with $\xi(O\cdot)$ by the invariance of the Laplacian under orthogonal transformations. 
Assuming the invariance $\xi(O\cdot)=_{law}\xi$, 
which is true for Gaussian ensembles,
would here not lead to further simplifications of the counterterm. 
However, this might be the case for other equations, e.~g.~the thin-film equation with thermal noise \cite{GT23}.


\subsection{Parameterization\texorpdfstring{ $(\lambda,p)\mapsto(h,\phi)$}{} of the solution manifold}\label{ss:para}

Obviously, for vanishing non-linearity, i.~e.~$\lambda=0$, the solution manifold of (\ref{ao40}) 
is an affine manifold. In view of 
\eqref{mt03}, it is an affine manifold over the 
linear space of all functions $p$ with $Lp=0$ in $\mathbb{R}^{1+d}$; 
by classical regularity theory for $L$, such solutions $p$ are 
(space-time) analytic functions, i.~e.~can be represented as convergent power series\footnote{We appeal in this heuristic to elliptic regularity theory, see \cite[Corollary~11.4.13]{Hor05}, which strictly speaking can not be applied to our parabolic $L$ as Tychonoff's example (see e.~g.~\cite[Theorem~8.6.7]{Hor03}) demonstrates. Let us emphasize however that this property will not be needed in the rigorous arguments of this article.}.
It is convenient to have the space of {\it all} analytic functions as parameter space.
We thus relax (\ref{ao40}) to hold only up to subtracting a (random) analytic function, we shall write \textit{modulo} analytic\footnote{The reader may wonder why $\xi_\rho$ is not absorbed in the additive analytic function, since by the Paley-Wiener-Schwartz theorem (see e.~g.\ \cite[Theorem~7.1.14]{Hor03}) $\xi_{\rho}$ is analytic as soon as the mollifier $\eta$ is compactly supported in Fourier space; we keep $\xi_{\rho}$ in the l.~h.~s., however, since we look for a parameterization of the solution manifold which remains robust as $\rho \to 0$.}
functions, i.~e.~
\begin{align}\label{ao40bis}
L\phi-(\lambda \phi^3+h^{(\rho)}\phi+\xi_\rho)= 0 \quad\mbox{ mod analytic functions},
\end{align}
where we now appeal to the fact that even if $Lp=\mbox{analytic}$, $p$ is analytic. 

\medskip

Let us pick a solution $\Pi_0$ of (\ref{ao40bis}) for $\lambda=0$, that is,
\begin{align}\label{ck24}
L\Pi_0-\xi_\rho=0 \quad\mbox{ mod analytic functions};
\end{align}
we will fix it in Subsection \ref{ss:unique}, 
and argue that it is canonical in Subsection \ref{ss:trans}.
This choice induces a parameterization for the solution manifold of (\ref{ao40bis}) for $\lambda=0$:
\begin{align}\label{ck35}
\mbox{for $\lambda=0$}:\quad\phi=\Pi_0+p,\quad\mbox{$p$ runs through analytic functions}.
\end{align}
It is tempting to think that 
-- and we shall do so for the purpose of this informal discussion -- 
such a parameterization persists in the presence of a non-linearity, i.~e.~for $\lambda\not=0$.
It is convenient to think of this parameterization in terms of the two components 
\begin{align}\label{ck35bis}
(\lambda,p)\mapsto(h,\phi)
\end{align}
or rather $(\lambda,p,\xi)\mapsto(h,\phi)$. 
In Subsections \ref{ss:pop} and
\ref{ss:unique} we will make natural choices which (at least informally) uniquely fix (\ref{ck35bis});
however, we will see in Subsection \ref{ss:trans} that (\ref{ck35bis}) is non-universal,
and depends on the (implicit) choice of an origin.

\begin{figure}[t]
\includegraphics[width=0.60\textwidth]{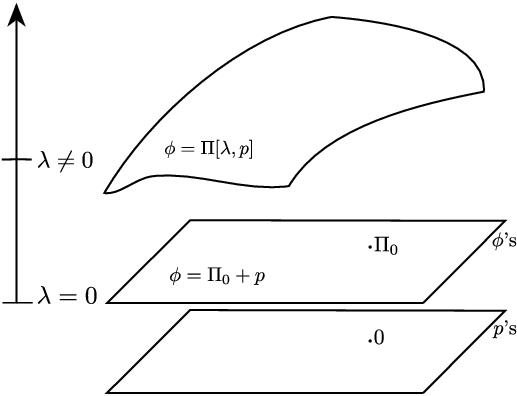}
\centering
\caption{
Heuristic visualization of the parameterization.
When $\lambda = 0$ the solution manifold is affine and parameterized by analytic functions $p$. When $\lambda \neq 0$ we expect it to be still parameterized by $p$ in a non-linear way.}
\label{fig:sol_manifold}
\end{figure}
%


\subsection{Power series representation\texorpdfstring{ $\{\Pi_\beta\}_\beta$}{} of the parameterization}\label{ss:new}

We now introduce coordinates on the parameter space of $(\lambda,p)$. 
Fixing somewhat arbitrarily a space-time origin, coordinates on the space 
of space-time analytic $p$'s are given by the coefficients of a (convergent) power
series representation, namely
\begin{align}\label{ao47}
\mathsf{z}_{\bf n}[p]:=\frac{1}{{\bf n}!}\partial^{\bf n}p(0) ,
\end{align}
where we recall that ${\bf n}=(n_0,\cdots,n_d)$ ranges over $\mathbb{N}_0^{1+d}$ and ${\bf n}!:=\prod_{i=0}^d(n_i!)$.
Since $\lambda$ multiplies a cubic term of the non-linearity, 
to be in line with the work \cite{LOTT21} on a general non-linearity,
we introduce the (here somewhat overblown) notation
\begin{align}\label{ao48}
\mathsf{z}_3[\lambda]=\lambda.
\end{align}

\medskip

We now ``algebrize'' the parameterization $(\lambda,p)\mapsto(h,\phi)$
by expressing $(h,\phi)$ as power series in the coordinates of $(\lambda,p)$: 
Recall that a multi-index $\beta$ over $\{3\}\sqcup\mathbb{N}_0^{1+d}$ 
gives rise to the monomial
\begin{align}\label{ck44}
\mathsf{z}^\beta:=\mathsf{z}_3^{\beta(3)}\prod_{{\bf n}\in\mathbb{N}_0^{1+d}}
\mathsf{z}_{\bf n}^{\beta({\bf n})}.
\end{align}
Inserting (\ref{ao47}) \& (\ref{ao48}) into 
(\ref{ck44}) defines an algebraic functional $\mathsf{z}^\beta[\lambda,p]$ on parameter space. 

\medskip

We now make 
an informal Ansatz for \eqref{ck35bis} 
by complementing the Ansatz \eqref{mt04} for $h$ with the informal 
\begin{align}\label{ao41}
\phi=\sum_{\beta}\mathsf{z}^\beta[\lambda,p]\Pi_\beta ,
\end{align}
where $\beta$ runs over all multi-indices,
and the $\Pi_\beta$'s are random space-time functions. 
Note that (\ref{ao41}) amounts to a separation of variables
into $(\lambda,p)$ on the one hand, and $x$ and the randomness on the other hand.
Even for fixed $\rho>0$,
there is no reason to believe that the 
series in (\ref{ao41}) 
is convergent. 
The main result provided in these notes is that for fixed $\beta$,
the coefficient $\Pi_\beta$ stays under control as $\rho\downarrow0$.

\medskip

For (\ref{ck35}) and (\ref{ao41}) to be consistent, it follows from
(\ref{ao47}) and Taylor's theorem that we must have
\begin{align}\label{ao46}
\Pi_{\beta}(y)=\left\{\begin{array}{cc}y^{\bf n}&\mbox{if}\;\beta=\delta_{\bf n}\\
0&\mbox{else}\end{array}\right\}\quad\mbox{provided}\;\beta(3)=0\;\mbox{but}\;\beta\not=0,
\end{align}
where 
as usual $y^{\bf n}$ $:=\prod_{i=0}^d y_i^{n_i}$,
and where $0$ denotes the multi-index that associates the value zero to all elements of $\{3\}\sqcup\mathbb{N}_0^{1+d}$.


\subsection{Characterization of\texorpdfstring{ $\Pi$}{} and\texorpdfstring{ $\Pi^{-}$}{} as formal power series}\label{ss:series}

It is convenient to compactify (\ref{ao41}) by interpreting $\Pi:=\{\Pi_\beta\}_\beta$ as
a ``formal power series''
in the (infinitely many) abstract variables $\mathsf{z}_3,\mathsf{z}_{\bf n}$
with coefficients in the space $X$ of random
space-time functions.\footnote{We note that $\Pi$ does not denote what in regularity structures is called the pre-model; 
rather it is centered at the base-point $0$ as can be seen in \eqref{ao47}, see also Subsection~\ref{ss:trans} for further details on the base-point dependence.} 
Likewise, motivated by \eqref{mt04} we interpret $c:=\{c_k\}_k$ as
a formal power series (actually just a polynomial by \eqref{ck57}) in $\mathsf{z}_3$ with deterministic scalar coefficients.
Despite its name, the notion of formal power series is rigorously defined;
and the connoisseur's notation\footnote{we should rather write $\mathbb{R}[[\mathsf{z}_3,\{\mathsf{z}_{\bf n}\}_{\bf n}]]$, but we do not for brevity} is $\Pi\in X[[\mathsf{z}_3,\mathsf{z}_{\bf n}]]$ 
and $c$ $\in\mathbb{R}[[\mathsf{z}_3]]$.
Provided the coefficient space, like here $X$ or $\mathbb{R}$, 
is an algebra, 
the formal power series space is an algebra under the multiplication rule
\begin{align}\label{ao45}
(\pi^{(1)}\pi^{(2)})_{\beta}=\sum_{\beta_1+\beta_2=\beta}\pi_{\beta_1}^{(1)}\pi_{\beta_2}^{(2)},
\end{align}
which is consistent with the usual multiplication when the power series actually converge.
Note that the unit element $\mathsf{1}$ of this algebra is characterized by 
$\mathsf{1}_\beta=0$ unless $\beta=0$, 
in which case the coefficient is given by the unit element of $X$. 

\medskip

Obviously, $\mathbb{R}[[\mathsf{z}_3]]$ can be considered as a sub-algebra
of $X[[\mathsf{z}_3,\mathsf{z}_{\bf n}]]$ so that next to
$\mathsf{z}_3\Pi^3$, 
also $c\Pi$ makes sense as an element in $X[[\mathsf{z}_3,\mathsf{z}_{\bf n}]]$,
which means that we identify
\begin{align}\label{ao45bis}
c_\beta=\left\{\begin{array}{cl}
c_k&\mbox{if}\;k=\beta(3)\;\mbox{and}\;\beta({\bf n})=0\;\mbox{for all}\;{\bf n}\\
0&\mbox{else}\end{array}\right\}.
\end{align}

\medskip

Informally, we identify $c$ with the function $h^{(\rho)}$ of $\lambda$, 
and $\Pi$ with the parameterization (\ref{ck35bis}): 
Indeed, via (\ref{ao47}) \& (\ref{ao48}) and in view of \eqref{mt04} and (\ref{ao41}), 
where we ignore the convergence issue of the latter,
$c$ associates a deterministic number to $\lambda$,
and $\Pi$ associates a random function to $(\lambda,p)$.
On this informal level, 
the PDE (\ref{ao40bis}) assumes the form\footnote{we note that $h = \sum_{\beta} h_{\beta} \mathsf{z}^{\beta}$ is related to $c$ by $c = \partial_{\mathsf{z}_{\mathbf{0}}} h$, so that the counterterm in \eqref{ao43} is in line with the corresponding exponential formula in \cite[(2.44)]{LOTT21}}
\begin{align}
L\Pi-\Pi^{-}&\;= 0 \quad\mbox{ mod analytic functions}\label{ao43bis}\\
\mbox{with}\quad\Pi^{-}&:=\mathsf{z}_3\Pi^3+c\Pi+\xi_\rho\mathsf{1},\label{ao43}
\end{align}
where \eqref{ao43bis} has to be understood component-wise.
Note that since by (\ref{ao43}), $\Pi^{-}_0=\xi_\rho$, 
(\ref{ao43bis}) is consistent with (\ref{ck24}).
According to \eqref{ck57}, (\ref{ao45}) and (\ref{ao45bis}), 
the component-wise version of identity (\ref{ao43}) reads
\begin{align}\label{ao44}
\Pi_{\beta}^{-}=\sum_{\delta_3+\beta_1+\beta_2+\beta_3=\beta}
\Pi_{\beta_1}\Pi_{\beta_2}\Pi_{\beta_3}
+\sum_{\substack{\beta_1+\beta_2=\beta\\ \beta_1(3)\not=0,\forall\;{\bf n}:\;\beta_1({\bf n})=0}}
c_{\beta_1}\Pi_{\beta_2}+\xi_\rho\delta_{\beta}^{0},
\end{align}
which reveals a strict triangularity of $\Pi\mapsto\Pi^{-}$
w.~r.~t.~the plain length $\beta(3)+\sum_{\bf n}\beta({\bf n})$ of the multi-index.
Hence (\ref{ao43bis}) \& (\ref{ao43}) suggests a hierarchical construction of $\Pi$ (at given $c$).
However, the construction will proceed by another inductive order, see Subsection \ref{ss:order};
the ingredients for this order
are the homogeneity $| \cdot |$ introduced in Subsection~\ref{ss:renorm}
and the noise homogeneity $[\cdot]$
introduced in the (following) Subsection \ref{ss:pop}.
Appealing to (\ref{ao46}) for $\Pi_{\delta_{\bf 0}}=1$, 
we find that the first examples are 
$\Pi_{0}^{-}=\xi_\rho$, $\Pi_{\delta_3}^{-}=\Pi_0^3+c_{1}\Pi_0$,
\begin{align}\label{ck68}
\Pi_{\delta_3+\delta_{\bf 0}}^{-}=3\Pi_0^2+c_{1},
\end{align}
$\Pi_{\delta_3+2\delta_{\bf 0}}^{-}=3\Pi_0$,
and $\Pi_{\delta_3+3\delta_{\bf 0}}^{-}=1$. The terms quickly become more complex:
e.~g.~$\Pi_{2\delta_3}^{-}$ $=3\Pi_0^2\Pi_{\delta_3}+c_{1}\Pi_{\delta_3}+c_{2}\Pi_0$ and
\begin{align}
\Pi_{2\delta_3+\delta_{\bf 0}}^{-}&=3\Pi_0^2\Pi_{\delta_3+\delta_{\bf 0}}+6\Pi_0\Pi_{\delta_3}
+c_{1}\Pi_{\delta_3+\delta_{\bf 0}}+c_{2}.\label{ck74bis}
\end{align}
In view of this combinatorial complexity, the task is to find an automated treatment.

\medskip

For comparison to the well-established tree-based approach \cite{Hai15}
we express these examples in the language of trees, 
where as usual the noise $\xi$ is represented by $\treeZero$,
inverting the operator $L$ is denoted by $|$,
and multiplication is denoted by attaching the trees at the root. 
Denoting the appropriately renormalized model of \cite{Hai15} by $\Pi_H$, 
we have
\begin{align*}
\Pi^-_0 &= \Pi_H (\treeZero), \ 
\Pi^-_{\delta_3} = \Pi_H(\treeThree), \ 
\Pi^-_{\delta_3+\delta_{\bf 0}} = 3\Pi_H(\treeThreeZero), \ 
\Pi^-_{\delta_3+2\delta_{\bf 0}} = 3\Pi_H(\treeThreeZeroZero), \\
\Pi^-_{2\delta_3} &= 3\Pi_H(\treeThreeThree\,), \ \
\Pi^-_{2\delta_3+\delta_{\bf 0}} = 9\Pi_H(\treeThreeThreeZeroA\,)+6\Pi_H(\treeThreeThreeZeroB\, ) .
\end{align*}
The compatibility between the (Hopf-)algebraic structures arising in recentering 
(positive renormalization) on trees and multi-indices 
was studied in \cite{LOT}, 
while the connection of the corresponding algebraic structures arising in  renormalization was investigated in \cite{BL23, Lin23}.
However, we point out that trees (and associated diagrams) do not play any role in our analysis. 


\subsection{Noise homogeneity\texorpdfstring{ $[\cdot]+1$}{} and population conditions\texorpdfstring{ on $\Pi,\Pi^{-}$}{}}\label{ss:pop}

We now motivate and make choices in the (informal) construction of the parameterization (\ref{ck35bis}).
These choices are guided by making $\Pi_\beta$ vanish  
for as many $\beta$'s as (algebraically) possible, 
thus maximizing the sparsity by minimizing the ``population'' of $\Pi$.
More precisely, we shall argue that we may postulate
\begin{align}\label{ao08}
\Pi_\beta\;
\left\{\begin{array}{ll}
=(\cdot)^{\bf n}&\mbox{if}\;\beta=\delta_{{\bf n}}\;\mbox{for some}\;{\bf n}\\
\in X&\mbox{for}\;[\beta]\ge 0\\
=0&\mbox{else}
\end{array}\right\},
\end{align}
where we introduced the notation
\begin{align}\label{ao10}
[\beta]:=2\beta(3)-\sum_{{\bf n}}\beta({\bf n})\in\mathbb{Z}.
\end{align}
This quantity is intimately related to
a simple invariance of the solution manifold of (\ref{ao40bis}), namely
\begin{align}\label{ap39}
\phi=\mu\hat\phi,\quad\xi=\mu\hat\xi\quad\mbox{and}\quad\lambda=\mu^{-2}\hat\lambda,
\end{align}
with $h^{(\rho)}$ and $\rho$ unchanged; by invariance we mean
that if $(\lambda,\phi,\xi)$ satisfies (\ref{ao40bis}), so does
$(\hat\lambda,\hat\phi,\hat\xi)$. 
In view of \eqref{ck35}, the parameterization $p$ transforms analogously to $\phi$, 
hence for $\hat p = \mu^{-1}p$ we postulate that the parameterization (\ref{ck35bis})
respects this invariance, meaning that we have 
\begin{align*}
(\hat\lambda,\hat p,\hat\xi)\mapsto(h,\hat\phi).
\end{align*}
On the level of the power series representation \eqref{mt04} and (\ref{ao41}), we read off that
this is satisfied if\footnote{this sufficient condition is not necessary since $\beta\mapsto[\beta]$ is not one-to-one}
\begin{align}\label{ap38}
\Pi_\beta=\mu^{[\beta]+1}\hat\Pi_\beta\quad\mbox{and}\quad
c_{k}=\mu^{2k}\hat c_{k}.
\end{align}
It is a good consistency check (and exercise) to verify that (\ref{ap38}) is compatible 
with (\ref{ao43bis}) \& (\ref{ao43}), leading to $\Pi_\beta^{-}$ 
$=\mu^{[\beta]+1}\hat\Pi_\beta^{-}$.
In view of the middle item in (\ref{ap39}) and the first item in (\ref{ap38}), 
$[\beta]+1$ can be interpreted as the homogeneity of $(\Pi_\beta,\Pi_\beta^{-})$ 
in the noise $\xi$. Thus (\ref{ao08}) postulates that $\Pi_\beta$
either has positive noise-homogeneity or is a polynomial.

\medskip

We shall establish (\ref{ao08}) alongside
\begin{align}\label{temp17}
\Pi_\beta^{-}\;\left\{\begin{array}{ll}
=\sigma_{\beta} (\cdot)^{\bf n}&\mbox{if}\;\beta=\delta_3+\sum_{j=1}^3\delta_{{\bf n}_j}\;
\mbox{for some}\;{\bf n}_j\\
\in X&\mbox{for}\;[\beta]\ge 0\\
=0&\mbox{else}
\end{array}\right\},
\end{align}
with ${\bf n}=\sum_{\bf m}{\bf m}\beta({\bf m})$ and the combinatorial factor
$\sigma_{\beta}:=\frac{(\sum_{\bf m} \beta ({\bf m})) !}{\prod_{\bf m}(\beta({\bf m})!)}$.
The $\beta$ appearing in (\ref{temp17}), namely
\begin{align}\label{ap65}
\beta=\delta_3+\sum_{j=1}^3\delta_{{\bf n}_j}\quad
\mbox{for some}\;{\bf n}_1,{\bf n}_2,{\bf n}_3,
\end{align}
will play a role throughout these notes.
We note that (\ref{ao08}) and (\ref{temp17}) imply
\begin{align}\label{temp24}
\begin{array}{l}\Pi_{\beta},\;\Pi^{-}_{\beta}\;\mbox{vanish
unless}\quad \beta=\delta_{\bf n},\;\beta=\delta_3+\sum_{j=1}^3\delta_{{\bf n}_j}
\;\mbox{or}\;[\beta]\ge 0,\\
\mbox{and thus}\;\Pi_{\beta},\;\Pi^{-}_{\beta}\;\mbox{vanish
unless}\quad[\beta]\ge -1,\\
\Pi_{\beta},\;\Pi^{-}_{\beta}\;\mbox{are polynomials unless}
\quad [\beta]\ge 0,
\end{array}
\end{align}
so that it is convenient to introduce the language
\begin{align}\label{ck54}
\begin{array}{l}
\beta\;\mbox{is ``purely polynomial (pp)''}\quad\mbox{iff}\quad
\beta=\delta_{\bf n}\;\mbox{for some ${\bf n}$},\\
\beta\;\mbox{is ``populated''}\quad\mbox{iff}\quad
\beta\;\mbox{is pp or of the form (\ref{ap65}) or}\;[\beta]\geq 0.
\end{array}
\end{align}

\medskip

\begin{proof}[Proof of \eqref{ao08} \textnormal{\&} \eqref{temp17}]
We now give the argument that \eqref{ao43bis} \& \eqref{ao43} are consistent with both (\ref{ao08}) and (\ref{temp17}) 
by induction in $k=\beta(3)$.
In the base case $k=0$, (\ref{ao08}) is just a reformulation of (\ref{ao46}). 
Still in the base case $k=0$,
we consider (\ref{ao44}) and note that the first and second r.~h.~s.~sums are empty,
so that $\Pi_\beta^{-}$ vanishes unless $\beta=0$, and thus $[\beta]=0$.

\medskip

We now turn to the induction step $k-1\leadsto k$ and
give ourselves a $\beta$ with $\beta(3)=k$ and $[\beta]\le-1$.
We aim at showing that $\Pi^-_\beta$ vanishes unless $\beta$ is of the form \eqref{ap65}.
We first consider $\Pi_\beta^{-}$ as given in (\ref{ao44}).
Clearly, the last term vanishes.
For the middle r.~h.~s.~term we note that the multi-indices involved satisfy
$\beta_1(3)+\beta_2(3)$ $=\beta(3)$, so that since $\beta_1(3)\ge 1$ 
and thus $\beta_2(3)<\beta(3)$, we may use the induction hypothesis on 
$\Pi_{\beta_2}$. Since the involved multi-indices also satisfy 
$[\beta_1]+[\beta_2]$ $=[\beta]$ $\le -1$, and since $[\beta_1]\ge 1$ we have $[\beta_2]\le -2$.
Hence we learn from (\ref{temp24}) that $\Pi_{\beta_2}=0$, so that the middle
r.~h.~s.~does not contribute. 
We finally turn to the first r.~h.~s.~term
in (\ref{ao44}); the involved multi-indices satisfy
$\beta_1(3)+\beta_2(3)+\beta_3(3)$ $=\beta(3)-1$, so that we may use the
induction hypothesis on $\Pi_{\beta_j}$. By definition (\ref{ao10}) they also satisfy
$[\beta_1]+[\beta_2]+[\beta_3]$ $=[\beta]-2\le-3$. Hence we learn from
(\ref{ao08}) that the $\beta_j$'s must be pp, and thus necessarily we must have 
$\beta=\delta_3+\delta_{{\bf n}_1}+\delta_{{\bf n}_2}+\delta_{{\bf n}_3}$.
The resulting contribution is $(\cdot)^{\bf n}$, and it arises $\sigma_{\beta}$ times.

\medskip

We now turn to the induction hypothesis for (\ref{ao08}).
Equipped with the one for (\ref{temp17}), in its reduced form of
(\ref{temp24}) on $\Pi^{-}_\beta$, it easily follows since we may absorb any
polynomial into the analytic function in (\ref{ao43bis}). 
\end{proof}

\medskip

Incidentally, with help of (\ref{ao10}), we may reformulate (\ref{ao07}) as
\begin{align}\label{ao12}
\lefteqn{|\beta|}\nonumber\\
&=(s-\frac{D}{2})([\beta]+1)+2\big(1+3\beta(3)-\sum_{\bf n}\beta({\bf n})\big)
+\sum_{\bf n}|{\bf n}|\beta({\bf n}).
\end{align}
The first term in (\ref{ao12}) corresponds to the effect of the noise, cf.~(\ref{ck22});
the second term corresponds to the effect of integration (i.~e.~inverting
the second-order operator $L$, whence the factor of $2$)
to which the cubic non-linearity contributes 3 units,
whereas a ``polynomial decoration'' removes one unit.
The last term captures the scaling of polynomials;
in particular we have the consistency
\begin{align}\label{temp04}
|\delta_{\bf n}|=|{\bf n}|.
\end{align}
This notion of homogeneity, 
which we motivated by scaling, 
is therefore also consistent with \cite[p.~199]{Hai15}.

\medskip

For frequent use in inductions, we retain that (\ref{ck28}) implies
\begin{align}\label{temp08}
|\cdot|-|0|\;\;\mbox{is additive},\quad\ge 0,\quad\mbox{and}\;=0\;\mbox{only if $\beta=0$},
\end{align}
and that by (\ref{ck30bis}), $|\cdot|$ is coercive, meaning
\begin{align}\label{temp15}
\{\,\beta\,|\,|\beta|<M\,\}\quad\mbox{is finite for every}\;M<\infty.
\end{align}
%


\subsection{Scale invariance of the solution manifold and homogeneity revisited}\label{ss:hom}

We would like our para\-me\-te\-ri\-zation (\ref{ck35bis})
to be consistent with 
the scaling invariance \eqref{ck32} - \eqref{ck26} and \eqref{mt02}.
In view of (\ref{ck35}) and (\ref{ap40}) we are poised to postulate the consistency in the form of
\begin{align}\label{ck34}
\mbox{for}\quad\hat p=r^{-\alpha}p(R\cdot)\quad\mbox{we have}\quad
(\hat\lambda,\hat p,\hat\xi)\mapsto(\hat h,\hat\phi).
\end{align}
We now derive the counterpart of this postulate on the level of the $\Pi_\beta$'s
and thus express (\ref{ck34}) in terms
of the coordinates (\ref{ao47}) \& (\ref{ao48}) on $(\lambda,p)$-space: 
Since (\ref{ck34}) implies that 
$\partial^{\bf n}\hat p =r^{|{\bf n}|-\alpha} (\partial^{\bf n}p)(R\cdot)$
we have
\begin{align*}
\mathsf{z}_{\bf n}[\hat p]=r^{|{\bf n}|-\alpha}\mathsf{z}_{\bf n}[p]
\quad\mbox{next to}\quad
\mathsf{z}_{3}[\hat\lambda]\stackrel{\eqref{ck26}}{=}
r^{2(1+\alpha)}\mathsf{z}_{3}[\lambda],
\end{align*}
where we recall that $|{\bf n}|$ is defined in \eqref{ck30bis}.
In terms of the monomials (\ref{ck44}), this yields
\begin{align*}
\mathsf{z}^\beta[\hat\lambda,\hat p]
=r^{\beta(3)2(1+\alpha)+\sum_{\bf n}\beta({\bf n})(|{\bf n}|-\alpha)}
\mathsf{z}^\beta[\lambda,p].
\end{align*}
Hence on the level of the power series representation \eqref{mt04} and (\ref{ao41}) the postulate (\ref{ck34})
holds if\footnote{the transformation rule (\ref{ck25}) would also be necessary if $\beta\mapsto|\beta|$ 
were one-to-one} the coefficients transform as
\begin{align}\label{ck25}
\hat c_\beta=r^{2-\beta(3)2(1+\alpha)}c_\beta\quad\mbox{and}\quad
\hat\Pi_\beta=r^{-|\beta|}\Pi_\beta(R\cdot),
\end{align}
where we recall that the homogeneity $|\beta|$ is defined in \eqref{ao07}.

\medskip

It is a good consistency check to verify that if $\hat\Pi^{-}$ is defined through
(\ref{ao43}) with $(\Pi,c,\xi)$ replaced by $(\hat\Pi,\hat c,\hat\xi)$, then
\begin{align}\label{ck33}
\hat\Pi_\beta^{-}=r^{2-|\beta|}\Pi_\beta^{-}(R\cdot).
\end{align}
\begin{proof}[Proof of \eqref{ck33}]
Indeed, for this we note that definition (\ref{ao07}) yields
\begin{align}\label{ck31}
k\delta_3+\sum_{j=1}^l\beta_j=\beta
\quad\Longrightarrow\quad
k2(1+\alpha)+\sum_{j=1}^l|\beta_j|-(l-1)\alpha=|\beta|.
\end{align}
Using (\ref{ck31}) for $k=1$ and $l=3$, we see that the first r.~h.~s.~term of (\ref{ao44})
only involves summands with $|\beta|=|\beta_1|+|\beta_2|+|\beta_3|+2$, as desired.
Using (\ref{ck31}) for $k=\beta_1(3)$ and $l=1$, we learn that the second r.~h.~s.~term
of (\ref{ao44}) only involves summands with $\beta_1(3)2(1+\alpha)+|\beta_2|$ $=|\beta|$,
which is what we want in view of the first item in (\ref{ck25}). 
For the last r.~h.~s.~term in (\ref{ao44}) it suffices to note
$|0|=\alpha=2+(s-\frac{D}{2})$, which is what we want in view of (\ref{ck32}).
\end{proof}


\subsection{Uniqueness of \texorpdfstring{$\Pi_\beta$}{Pi} given \texorpdfstring{$\Pi_\beta^{-}$}{Pi-}}\label{ss:unique}

For given $\Pi^{-}_\beta$, the solution $\Pi_\beta$ of the linear PDE (\ref{ao43bis})
is only determined up to an analytic function. In this subsection, we fix this degree of freedom
and start with the following remark:
In view of (\ref{ck32}), the second item in (\ref{ck25}), and (\ref{ck33}), 
it is natural to expect that in the limit of vanishing regularization\footnote{when $c$ diverges}, 
\begin{align}\label{ck27}
\xi(R\cdot)=_{law} r^{s-\frac{D}{2}}\xi\quad\stackrel{\rho\downarrow 0}{\Longrightarrow}\quad
\left\{\begin{array}{ccl}\Pi_\beta(R\cdot)&=_{law}&r^{|\beta|}\Pi_\beta,\\
\Pi_\beta^{-}(R\cdot)&=_{law}&r^{|\beta|-2}\Pi_\beta^{-},\end{array}\right.
\end{align}
where we note the consistency with \eqref{lb01}. 
We recall from Subsection~\ref{ss:intro} that this is an informal way of stating
\begin{align}\label{ck27bis}
&\xi(R\cdot)=_{law} r^{s-\frac{D}{2}}\xi\quad
\stackrel{\rho\downarrow 0}{\Longrightarrow}\nonumber\\
&\!\textnormal{laws of }r^{-|\beta|}\Pi_{\beta r}(0)\textnormal{ and }
r^{2-|\beta|}\Pi_{\beta r}^{-}(0)\textnormal{ do not depend on }r\!\in\!(0,\!\infty)
\end{align}
for any Schwartz function $\psi$, see (\ref{ck30}) for the notation.
This motivates the purely qualitative\footnote{in this article we understand by the term ``qualitative'' that a certain quantity is finite, while by ``quantitative'' we mean an actual estimate of this quantity independent of $\rho$} postulate 
\begin{align}\label{ck52}
\left.\begin{array}{c}
\limsup_{r\downarrow   0}r^{-|\beta|}\mathbb{E}|\Pi_{\beta r}(0)|<\infty\\
\limsup_{r\uparrow\infty}r^{-|\beta|}\mathbb{E}|\Pi_{\beta r}(0)|<\infty
\end{array}\right\}\quad\mbox{uniformly in bounded $\psi$},
\end{align}
where the boundedness refers to the semi-norms (\ref{ck52a}).

\medskip

We claim that in conjunction with the assumption
\begin{align}\label{ck53}
s\quad\mbox{is irrational},
\end{align}
the qualitative small and large scale estimate (\ref{ck52}) implies uniqueness of $\Pi_\beta$. 
Thanks to (\ref{ao08}), we may restrict ourselves to $\beta$ with $[\beta]\ge 0$.
The purpose of (\ref{ck53}) is to ensure 
\begin{align}\label{temp04bis}
[\beta]\ge0
\quad\Longrightarrow\quad|\beta|\not\in\mathbb{Z},
\end{align}
which follows from the representation (\ref{ao12}), and which can be seen as
a reverse of (\ref{temp04}).

\begin{proof}[Proof of uniqueness for $\Pi_\beta$]
Suppose there are two versions of $\Pi_\beta$; consider their difference $f$,
which by (\ref{ao43bis}) satisfies $Lf=\mbox{analytic}$, in a distributional and almost sure
sense.
Given an arbitrary bounded test random variable $F$, we consider $\bar f:=\mathbb{E}fF$
which satisfies $L\bar f=\mbox{analytic}$ in a distributional sense and thus is analytic. 
By the second part of our postulate (\ref{ck52}) 
we have $\limsup_{r\uparrow\infty}r^{-|\beta|}\int dx
\psi_r\bar f<\infty$ for any Schwartz function $\psi$. Replacing $\psi$ by $\partial^{\bf n}\psi$
and using that $(\partial^{\bf n}\psi)_r$ $=r^{|{\bf n}|}\partial^{\bf n}\psi_r$,
we see that this implies $\lim_{r\uparrow\infty}\int dx
\psi_r\partial^{\bf n}\bar f=0$ provided $|{\bf n}|>|\beta|$. By (a minor extension of)
the Liouville theorem for analytic functions this implies that $\partial^{\bf n}\bar f\equiv0$.
Hence $\bar f$ must be a polynomial of (parabolic) degree $\le|\beta|$, where
\begin{align*}
\mbox{degree of a polynomial $p$}:=\max\{\,|{\bf n}|\,|\,\partial^{\bf n}p\not\equiv0\,\}
\in\mathbb{N}_0\cup\{-\infty\}.
\end{align*} 
Hence if $|\beta|<0$, we are done; if $|\beta|\ge 0$, we turn 
to the first part of our postulate, which implies $\limsup_{r\downarrow 0}r^{-|\beta|}\int dx
\bar f\psi_r<\infty$ for any Schwartz function. Replacing $\psi$ once more by 
$\partial^{\bf n}\psi$ and fixing a $\psi$ of unit integral, we learn
that $\partial^{\bf n}\bar f(0)=\lim_{r\downarrow 0}\int dx
\partial^{\bf n}\bar f\psi_r=0$ provided $|{\bf n}|<|\beta|$.
Hence $\bar f$ has degree $|\beta|$ or vanishes. In view of (\ref{temp04bis}),
we must have the latter. Since $F$ was arbitrary, $f\equiv 0$ almost surely.
\end{proof}

By the same argument, and in line with (\ref{ck27bis}), 
we also learn that (\ref{ao43bis}) sharpens to
\begin{align}\label{ao43ter}
L\Pi_\beta-\Pi^{-}_\beta= 0 \mbox{ mod (random) polynomials of degree $\le|\beta|-2$}.
\end{align}
In fact, recalling \eqref{ao08} and \eqref{temp17}, it even sharpens\footnote{this further strengthening crucially relies on working in the whole $\mathbb{R}^{1+d}$, whereas \eqref{ao43ter} also holds in the space-time periodic setting, see e.~g.~\cite[Appendix~B]{BOS25}} to 
\begin{align*}
L\Pi_\beta-\Pi_\beta^{-}=\left\{\begin{array}{ll}
0 &\textnormal{for }[\beta]\geq0 ,\\
L(\cdot)^{\bf n} & \textnormal{for }\beta=\delta_{\bf n} ,\\
-\sigma_\beta(\cdot)^{{\bf n}_1+{\bf n}_2+{\bf n}_3} &\textnormal{for }\beta=\delta_3+\delta_{{\bf n}_1}+\delta_{{\bf n}_2}+\delta_{{\bf n}_3} ,
\end{array}\right. 
\end{align*}
with $\sigma_{\beta}=\frac{(\sum_{\bf m} \beta ({\bf m})) !}{\prod_{\bf m}(\beta({\bf m})!)}$ as in \eqref{temp17}. 
Incidentally, it follows from this and the ansatz \eqref{ao41} that (informally) \eqref{ao40bis} sharpens to 
\begin{equation*}
L\phi-(\lambda\phi^3+h^{(\rho)}\phi+\xi_\rho)
\stackrel{\eqref{ao41}}{=}\sum_\beta (L\Pi_\beta-\Pi_\beta^-)\mathsf{z}^\beta[\lambda,p]
=Lp-\lambda p^3 ,
\end{equation*}
which is consistent with \eqref{ck35}.


\subsection{Uniqueness of \texorpdfstring{$c$}{c} via BPHZ-choice of renormalization}\label{ss:c}

We now fix the last degree of freedom, 
namely the coefficients $\{c_k\}_{k\ge 1}$ of the counterterm $h^{(\rho)}$ 
by making 
one further postulate:
We impose the analogue of (\ref{ck52}) for $\Pi^{-}$:
\begin{align}\label{ck52bis}
\lim_{r\uparrow\infty}r^{2-|\beta|}\mathbb{E}|\Pi_{\beta r}^{-}(0)|<\infty
\end{align}
for any Schwartz kernel $\psi$, which again is motivated by (\ref{ck27bis}).
In fact, we shall argue that the (finitely many) degrees of freedom of \eqref{ck57}
are fixed 
by the following consequence of (\ref{ck52bis}) 
\begin{align}\label{ck52ter}
\lim_{r\uparrow\infty}\mathbb{E}\Pi_{\beta r}^{-}(0)=0\quad\mbox{provided}\;|\beta|<2.
\end{align}
Fixing the counterterm by imposing vanishing expectations is 
reminiscent of what in regularity structures is called the BPHZ-choice of renormalization.

\medskip

Before proceeding with the proof that \eqref{ck57} \textnormal{\&} \eqref{ck52ter} fix $c$, we
first identify
those populated multi-indices $\beta$, cf.~(\ref{ck54}), with $|\beta|<2$.
To this purpose we observe that by definition (\ref{ao07}) and the range (\ref{ao04}), 
we have for any multi-index $\beta$
\begin{align}\label{ap41}
\beta({\bf n})=0\;\mbox{for all ${\bf n}$}\quad\mbox{or}\quad|\beta|
\ge\sum_{\bf n}|{\bf n}|\beta({\bf n}).
\end{align}
Hence $|\beta|<2$ implies $\beta({\bf n})=0$ unless $|{\bf n}|<2$, which by
(\ref{ck30bis}) only leaves the $1+d$ cases ${\bf n}\in\{{\bf 0},{\bf e}_1,\cdots,{\bf e}_d\}$.
Moreover, we must have $\beta({\bf e}_i)\le 1$. However, (\ref{ap41}) does not put any restriction
on $\beta({\bf 0})$. For this, we note that in the range (\ref{ao04}),
\begin{align*}
\mbox{for $\beta$ with $[\beta]\ge -1$}:\quad
|\beta|+1\ge\sum_{\bf n}(|{\bf n}|+1)\beta({\bf n}),
\end{align*}
which follows from using $2\beta(3)$ $\ge -1+\sum_{\bf n}\beta({\bf n})$ on these $\beta$'s.
Hence $|\beta|<2$ implies $\sum_{\bf n}\beta({\bf n})\le 2$, and in particular 
$\beta({\bf 0})\le 2$. In conclusion,
we learn that there are only four classes of populated multi-indices with $|\beta|<2$, namely
\begin{align}\label{ap42}
\begin{array}{clccr}
{\rm (I)}&\beta=k\delta_3,
&\;|\beta|&=&\alpha+2k(1+\alpha),\\
{\rm (II)}&\beta=k\delta_3+\delta_{\bf 0},
&\;|\beta|&=&2k(1+\alpha),\\
{\rm (III)}&\beta=k\delta_3+2\delta_{\bf 0}\;\;\mbox{for}\;k\ge 1,
&\;|\beta|&=&-\alpha+2k(1+\alpha),\\
{\rm (IV)}&\beta=k\delta_3+\delta_{{\bf e}_i}\;\;\mbox{for}\;1\le i\le d,
&\;|\beta|&=&1+2k(1+\alpha).
\end{array}
\end{align}

\medskip

Incidentally, we learn from (\ref{ao44}) in conjunction with the uniqueness statement from Subsection 
\ref{ss:unique} that parity in law (\ref{ck21}) propagates
by induction in $\sum_{\bf n}\beta({\bf n})$:
\begin{align*}
(-1)^{1+\sum_{\bf n}\beta({\bf n})}\Pi_\beta=_{law}\Pi_\beta\quad\mbox{and}\quad
(-1)^{1+\sum_{\bf n}\beta({\bf n})}\Pi_\beta^{-}=_{law}\Pi_\beta^{-}.
\end{align*}
Likewise, we see that the symmetry in law (\ref{ck20}) under the reflection $R_i$,
together with the plain evenness/oddness of $\Pi_\beta$ for $\beta\not=0$ and $\beta(3)=0$ 
under $R_i$, cf.~(\ref{ao46}), propagates:
\begin{align*}
(-1)^{\sum_{\bf n}n_i\beta({\bf n})}\Pi_\beta(R_i\cdot)=_{law}\Pi_\beta\quad\mbox{and}\quad
(-1)^{\sum_{\bf n}n_i\beta({\bf n})}\Pi_\beta^{-}(R_i\cdot)=_{law}\Pi_\beta^{-}.
\end{align*}
Hence (\ref{ck52ter}) is automatically satisfied for the classes I, III, and IV 
(provided $\psi$ is spatially even).
This is the reason why we only need the counterterm $h^{(\rho)}\phi$,
and not those of the form
$h^{(\rho)}$, $h^{(\rho)}_i\partial_i\phi$, and $h^{(\rho)}\phi^2$.

\begin{proof}[Proof that \eqref{ck57} \textnormal{\&} \eqref{ck52ter} fix $c$]
We now turn to the uniqueness argument:
The following (semi-strict) triangular structure can be read off (\ref{ao44}):
\begin{align}\label{ck60}
\begin{array}{r}
\Pi_\beta^{-}\quad\mbox{depends on}\;\Pi_\gamma\;\mbox{only for}\;\gamma(3)<\beta(3),\\
\mbox{and on}\;c_{l}\;\mbox{only for}\;l\le\beta(3).
\end{array}
\end{align}
This non-strictness in the $c$-dependence is compensated by strictness for
$\beta$ of class II, cf.~(\ref{ap42}), in the sense of:
\begin{align}\label{ck59}
\Pi_{k\delta_3+\delta_{\bf 0}}^{-}-c_{k}\quad\mbox{depends on}
\;c_{l}\;\mbox{only for}\;l<k,
\end{align}
which follows from the fact that by (\ref{ao46}) \textnormal{\&} (\ref{ao45bis}), 
the middle r.~h.~s.~ term of (\ref{ao44}) can be re-written for $\beta=k\delta_3+\delta_{\bf 0}$ as
\begin{align*}
\sum_{\beta_1+\beta_2=k\delta_3+\delta_{\bf 0}}c_{\beta_1}\Pi_{\beta_2}
=c_{k}+\sum_{l=1}^{k-1}c_{l}\Pi_{(k-l)\delta_3+\delta_{\bf 0}}.
\end{align*}
We finally note that by our postulate (\ref{ck52ter}) 
and 
the fact that the space-time constant $c_{k}$ is deterministic, recall \eqref{mt08}, 
we have
\begin{align}\label{ck58}
c_{k}=-\lim_{r\uparrow\infty}\mathbb{E}(\Pi^{-}_{k\delta_3+\delta_{\bf 0}}-c_{k})_r(0)
\quad\mbox{provided}\;(1+\alpha)k<1.
\end{align}
Hence uniqueness follows by an induction in $k=\beta(3)\ge 1$ where inside each induction
step, we start with $k\delta_3+\delta_{\bf 0}$, which determines $c_{k}$
by (\ref{ck58}), since by (\ref{ck59}) and the induction hypothesis, 
$\Pi_{k\delta_3+\delta_{\bf 0}}^{-}-c_{k}$ is determined. We then deal with the
$d+2$ remaining multi-indices
$k\delta_3$, $k\delta_3+2\delta_{\bf 0}$, and $k\delta_3+\delta_{\bf e_i}$,
(in any order), where now (\ref{ck60}) is sufficient to appeal to the induction
hypothesis, because $c_{k}$ is already determined.

\medskip

We note that this argument (implicitly) relies on the following strict triangularity
\begin{align}\label{ck61}
\Pi_\beta^{-}\quad\mbox{depends on}\;\Pi_\gamma\;\mbox{only for}\;|\gamma|<|\beta|,
\end{align}
which follows from glancing at (\ref{ao44}): The multi-indices contributing
to the first r.~h.~s.~term are by (\ref{temp08}) related by
$(|\delta_3|-|0|)$ $+(|\beta_1|-|0|)$ $+(|\beta_2|-|0|)$ $+(|\beta_3|-|0|)$ 
$=|\beta|-|0|$;
again by (\ref{temp08}) the first bracket is $>0$ and all others are at least $\ge 0$,
so that necessarily $|\beta_j|<|\beta|$ for $j=1,2,3$, as desired.
The multi-indices contributing to the second r.~h.~s.~term 
are related by $(|\beta_1|-|0|)+|\beta_2|=|\beta|$, and we again obtain because
of $\beta_1\not=0$ that $|\beta_2|<|\beta|$.
\end{proof}

In conclusion, we have argued that $(\Pi,\Pi^{-},c)$ can be uniquely constructed,
so that informally, we have now fixed the parameterization $(\lambda,p)\mapsto(h,\phi)$.


\subsection{Shift-invariance of the solution manifold,
general base-points\texorpdfstring{ $x$}{}, and corresponding \texorpdfstring{$\Pi_x$}{centered model}}\label{ss:trans}

We return to the informal discussion of the solution manifold.
Already in Subsection \ref{ss:intro}, we appealed to its
invariance under shift $(\phi(\cdot+x),\xi(\cdot+x))$. We will learn over
the course of the next subsections, and ultimately at the end of Subsection \ref{ss:expect},
that the parameterization (\ref{ck35bis}),
which now is fixed thanks to the choices made in Subsections \ref{ss:pop}, \ref{ss:unique}, 
and \ref{ss:c}, does typically {\it not} respect this invariance: 
\begin{align}\label{ap30}
(\lambda,p,\xi)\mapsto\phi\quad\not\hspace{-.75ex}\Longrightarrow\quad(\lambda,p(\cdot+x),\xi(\cdot+x))
\mapsto\phi(\cdot+x),
\end{align}
a fact which we will prove over the course of Subsections~\ref{ss:trafo}, \ref{ss:coord}, and \ref{ss:improve} below, see \eqref{ap31}, \eqref{ap34}, and the discussion after \eqref{ap36}.
Hence while (\ref{ck35bis}) is unique, it is not canonical 
as it depends on the choice of an origin.

\medskip

This motivates to repeat the definition from Subsection \ref{ss:series}
with the origin $0$ replaced by a general ``base-point'' $x\in\mathbb{R}^{1+d}$
in (\ref{ao47}), which means that 
(\ref{ao41}) is replaced by 
\begin{align}\label{ck39}
\phi=\sum_{\beta}\lambda^{\beta(3)}\prod_{\bf n}
\big(\tfrac{1}{{\bf n}!}\partial^{\bf n}p(x)\big)^{\beta({\bf n})}\Pi_{x\beta}.
\end{align}
At least a priori (and in fact as we shall see), this defines a parameterization 
$(\lambda,p)\mapsto(h,\phi)$ of the solution manifold that is different from (\ref{ao41}): 
The same $p$'s (and same $\lambda$) will give different $\phi$'s. 
Now the analogue of (\ref{ao46}) reads
\begin{align}\label{ck38}
\Pi_{x\beta}(y)=\left\{\begin{array}{cc}(y-x)^{\bf n}&\mbox{if}\;\beta=\delta_{\bf n}\\
0&\mbox{else}\end{array}\right\}\quad\mbox{for}\;\beta(3)=0\;\mbox{and}\;\beta\not=0.
\end{align}
We assemble these coefficients into $\Pi_x\in X[[\mathsf{z}_{3},\mathsf{z}_{\bf n}]]$,
which in view of (\ref{ck38}) does depend on $x$.
Since $c$ provides the $p$-independent power series representation of $h$,
it stays unaffected by the change of base-point; in line with (\ref{ao43}) we set
\begin{align}\label{ck84}
\Pi_x^{-}:=\mathsf{z}_3\Pi_x^3+c\Pi_x+\xi_\rho\mathsf{1},
\end{align}
and have the analogue of (\ref{ao43ter}), namely
\begin{align}\label{ck93}
L\Pi_{x\beta}=\Pi_{x\beta}^{-} \quad\mbox{ mod polynomials of degree $\le|\beta|-2$}.
\end{align}
In view of the uniqueness of the construction, we obtain that
the family $\xi\mapsto\{\Pi_x\}_x$ is covariant under shift in the sense of
\begin{align}\label{ap32}
\Pi[\xi(\cdot+x)](y)=\Pi_x[\xi](y+x)\quad\mbox{for all}\;x,y\in\mathbb{R}^{1+d}.
\end{align}
From the covariance (\ref{ap32}) in conjunction with the stationarity 
assumption (\ref{ck19}) we obtain an analogue of (\ref{ck27})
\begin{align}\label{temp30}
\Pi_x(\cdot+x)=_{law}\Pi\quad\mbox{and}\quad\Pi_x^{-}(\cdot+x)=_{law}\Pi^{-}.
\end{align}

\medskip

Clearly, by (\ref{ck84}) we have $\Pi_{x\,\beta=0}^{-}=\xi_\rho=\Pi_{\beta=0}^{-}$;
we now claim that this translates into 
\begin{align}\label{ck62}
\Pi_{x\,\beta=0}=\Pi_{\beta=0}.
\end{align}
This shows that the definition of $\Pi_0$, and thus at least the anchoring (\ref{ck35}) 
of our parameterization (\ref{ck35bis}) at $\lambda=0$, was canonical.
The argument for (\ref{ck62}) is similar to the one given in Subsection \ref{ss:unique}.

\begin{proof}[Proof of \eqref{ck62}]
We first note that by (\ref{temp30}), the second item in (\ref{ck52})
translates into $\limsup_{r\uparrow\infty}$ $r^{-|\beta|}$ $\mathbb{E}|\Pi_{x\beta r}(x)|$
$<\infty$, uniformly for bounded $\psi$. Writing $\psi_r(-y)=\psi^{(r, x)}_r(x-y)$ for some Schwartz function $\psi^{(r,x)}$ such that $\{\psi^{(r,x)}\}_{r\uparrow\infty}$ is bounded in terms of (\ref{ck52a}), the
above implies $\limsup_{r\uparrow\infty}r^{-|\beta|}\mathbb{E}|\Pi_{x\beta r}(0)|<\infty$.
Together with the second item in (\ref{ck52}) in its original version
and with (\ref{ck28}) we obtain for $f=:\Pi_{x\,\beta=0}-\Pi_{\beta=0}$ that
$\lim_{r\uparrow\infty}\mathbb{E}|f_{r}(0)|=0$. On the other hand,
we have by (\ref{ao43ter}) and $\Pi_{x\,\beta=0}^{-}=\Pi_{\beta=0}^{-}$ that $Lf=0$.
We thus may argue as in Subsection \ref{ss:unique} that $f=0$.
\end{proof}


\subsection{The change-of-base-point transformation\texorpdfstring{ $\Gamma_x^*$}{}}\label{ss:trafo}

We continue with the informal discussion of the solution manifold.
By construction and at fixed $\lambda$ and realization $\xi$, 
the r.~h.~s.~of (\ref{ck39}) captures all 
solutions of (\ref{ao40bis}) when $p$ runs through all analytic functions.
Replacing $p$ by ${q(\cdot-x)}$, and letting $q$ run through all analytic functions,
we obviously again obtain a parameterization $(\lambda,q)\mapsto(h,\phi)$ of the solution manifold;
since this means replacing $\frac{1}{{\bf n}!}\partial^{\bf n}p(x)$ by $\mathsf{z}_{\bf n}[q]$ 
it takes the compact form of
\begin{align*}
\phi=\sum_{\beta}\mathsf{z}^\beta[\lambda,q]\Pi_{x\beta}.
\end{align*}
This and (\ref{ao41}) provide two different parameterizations; hence there exists a (random)
parameter transformation 
\begin{align}\label{ck66}
(\lambda,p)\mapsto(\lambda,q=p_{\lambda\;x})
\end{align}
such that 
\begin{align}\label{ck63}
\sum_{\beta}\mathsf{z}^\beta[\lambda,p]\Pi_{\beta}
=\sum_{\beta}\mathsf{z}^\beta[\lambda,p_{\lambda\;x}]\Pi_{x\beta}.
\end{align}
We remark that (\ref{ap30}) translates into
\begin{align}\label{ap31}
\mbox{it is {\it not} true that}\quad p_{\lambda\;x}=p(\cdot+x)\quad\mbox{for all}\;(\lambda,p) .
\end{align}
Indeed, by \eqref{ao41}, (\ref{ap30}) means that
$\sum_{\beta}\mathsf{z}^\beta[\lambda,p(\cdot+x)]\Pi_{\beta}[\xi(\cdot+x)](y)$
does not agree with
$\sum_{\beta}\mathsf{z}^\beta[\lambda,p]\Pi_{\beta}[\xi](y+x)$.
According to (\ref{ap32}), the former coincides with
$\sum_{\beta}\mathsf{z}^\beta[\lambda,p(\cdot+x)]\Pi_{x\beta}[\xi](y+x)$,
while by (\ref{ck63}), the latter can be written as
$\sum_{\beta}\mathsf{z}^\beta[\lambda,p_{\lambda\;x}]\Pi_{x\beta}[\xi](y+x)$.

\medskip
 
Note that in view of (\ref{ao47}) \& (\ref{ao48}), elements $\pi$ of 
$\mathbb{R}[[\mathsf{z}_3,\mathsf{z}_{\bf n}]]$ can informally be considered
as function(al)s on the $(\lambda,p)$-space.
Hence the nonlinear transformation (\ref{ck66}) induces by pull back a linear endomorphism\footnote{implicitly, $\Gamma^*_x$ is also indexed by the base-point $0$, similarly to $\Pi$; we drop this dependence in the notation for convenience, see Subsection~\ref{ss:triang2} for further details on the base-point dependence} $\Gamma_x^*$ 
of $\mathbb{R}[[\mathsf{z}_3,\mathsf{z}_{\bf n}]]$ via
\begin{align}\label{ck70}
(\Gamma_x^*\pi)[\lambda,p]=\pi[\lambda,p_{\lambda\;x}];
\end{align}
the $*$-notation will be motivated in the (next) Subsection \ref{ss:dual}.
Since the product (\ref{ao45}) on $\mathbb{R}[[\mathsf{z}_3,\mathsf{z}_{\bf n}]]$ extends the
product of function(al)s on $(\lambda,p)$-space, $\Gamma_x^*$ is
an algebra endomorphism, which means
\begin{align}\label{temp05}
\Gamma_x^*\pi\pi'=(\Gamma_x^*\pi)\Gamma_x^*\pi'\quad\mbox{and}\quad\Gamma_x^*\mathsf{1}=\mathsf{1}.
\end{align}
By (\ref{ao48}) and (\ref{ck70}), 
the triviality of the first component of (\ref{ck66}) translates into
\begin{align}\label{temp06}
\Gamma_x^*\mathsf{z}_3=\mathsf{z}_3,
\end{align}
and once more by (\ref{ck70}), (\ref{ck63}) translates into 
\begin{align}\label{temp02}
\Pi=\Gamma_x^*\Pi_x.
\end{align}
Since $c\in\mathbb{R}[[\mathsf{z}_3]]$, 
we immediately obtain from the rules (\ref{temp05}) \& (\ref{temp06}) that 
\begin{align}\label{ck94}
\Gamma_x^*c=c ,
\end{align}
which is not surprising since $c\in\mathbb{R}[[\mathsf{z}_3]]$ encodes that $c$ is independent of $p$, and by \eqref{ck66} the transformation acts only on the $p$ variable. 
Using \eqref{ck94} we learn from applying $\Gamma_x^*$ to (\ref{ao43}) and comparing with (\ref{ck84})
that (\ref{temp02}) transmits to $\Pi^{-}$:
\begin{align}
\Pi^{-}=\Gamma_x^*\Pi_x^{-}\label{temp01}.
\end{align}
\begin{figure}[t]
\includegraphics[width=0.9\textwidth]{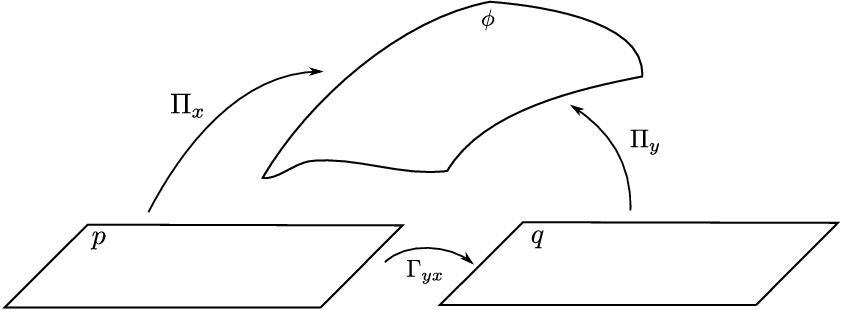}
\centering
\caption{
Heuristic visualization of the change-of-base-point transformation $\Gamma$.
Informally, $\Pi_x$ and $\Pi_y$ act as (inverse) ``charts'' on the solution manifold of the $\phi$'s, while $\Gamma_{y x}$ acts as a ``transition function'' between these two ``charts''.}
\label{fig:transition}
\end{figure}


\subsection{\texorpdfstring{$\Gamma_x$}{Gamma} acts by shift on space-time polynomials}\label{ss:dual}

In fact, $\Gamma_x^*$ is the algebraic transpose of a linear endomorphism $\Gamma_x$ of
$\mathbb{R}[\mathsf{z}_3,\mathsf{z}_{\bf n}]$, where the latter denotes
the space of polynomials in the variables $\mathsf{z}_3$ and 
$\{\mathsf{z}_{\bf n}\}_{\bf n}$, of which $\mathbb{R}[[\mathsf{z}_3,\mathsf{z}_{\bf n}]]$ 
is the canonical (algebraic) dual.
The linear space of space-time polynomials $p\in\mathbb{R}[x_0,\cdots,x_d]$ is
canonically embedded into $\mathbb{R}[\mathsf{z}_3,\mathsf{z}_{\bf n}]$ 
by specifying how the dual basis acts on them
\begin{align}\label{ck67}
\mathsf{z}^\beta.p=\mathsf{z}^\beta[\lambda=0,p],
\end{align}
where we note that the r.~h.~s.~makes (more than informal) sense since $p$
is a polynomial and not just an analytic function.
We now informally argue that $\Gamma_x$ acts on this subspace by translation:
\begin{align}\label{ck40}
\Gamma_x p=p(\cdot+x).
\end{align}

\medskip

In order to informally establish (\ref{ck40}), we return to
the parameter transformation (\ref{ck66}), and consider the case of $\lambda=0$:
By (\ref{ck35}) and (\ref{ao41}) we have $\sum_{\beta}\mathsf{z}^\beta[0,p]\Pi_\beta$ $=\Pi_0+p$.
For the general base-point $x$ this takes the form of 
$\sum_{\beta}\mathsf{z}^\beta[0,p(\cdot+x)]\Pi_{x\beta}$ $=\Pi_{x0}+p$.
Using (\ref{ck63}) to rewrite the former as 
$\sum_{\beta}\mathsf{z}^\beta[0,p_{\lambda=0\;x}]\Pi_{x\beta}$ $=\Pi_0+p$
and (\ref{ck62}) to formulate the latter as 
$\sum_{\beta}\mathsf{z}^\beta[0,p(\cdot+x)]\Pi_{x\beta}$ $=\Pi_{0}+p$ we informally deduce
\begin{align}\label{ck64}
p_{\lambda=0\;x}=p(\cdot+x),
\end{align}
so that (\ref{ap31}) holds at least partially.
Inserting (\ref{ck64}) into (\ref{ck70}) yields $(\Gamma_x^*\pi)[0,p]$
$=\pi[0,p(\cdot+x)]$, which in view of (\ref{ck67}) is to be interpreted as (\ref{ck40}).


\subsection{Matrix representation\texorpdfstring{ $\{(\Gamma_x^*)_\beta^\gamma\}$ 
of $\Gamma_x^*$}{}}\label{ss:coord}

Since the polynomial space $\mathbb{R}[\mathsf{z}_3,\mathsf{z}_{\bf n}]$ has a natural
(algebraic) basis\footnote{as opposed to its dual $\mathbb{R}[[\mathsf{z}_3,\mathsf{z}_{\bf n}]]$
of which the monomials are not a basis} indexed by multi-indices $\beta$,
$\Gamma_x$ admits a matrix representation $\{(\Gamma_x)_\gamma^\beta\}$.
Its transpose $(\Gamma_x^*)_\beta^\gamma$ $=(\Gamma_x)_\gamma^\beta$ allows us to 
express the action of $\Gamma_x^*$ coordinate-wise:
\begin{align}\label{ck74}
(\Gamma_x^*\pi)_\beta=\sum_{\gamma}(\Gamma_x^*)_\beta^\gamma\pi_\gamma, 
\quad\mbox{in particular}\quad(\Gamma_x^*)^\gamma_\beta=(\Gamma_x^*\mathsf{z}^\gamma)_\beta.
\end{align}
We note that the sum is effectively finite since by the nature of a matrix representation
\begin{align}\label{ck88}
\{\,\gamma\,|\,(\Gamma^*_x)_\beta^\gamma\not=0\,\}
=\{\,\gamma\,|\,(\Gamma_x)_\gamma^\beta\not=0\,\}\quad\mbox{is finite for every $\beta$}.
\end{align}

\medskip

We now capture (\ref{ck40}) on the level of this matrix representation.
We have by Leibniz' formula
\begin{align}\label{ap33}
\mathsf{z}_{\bf n}[p(\cdot+x)]
\stackrel{\eqref{ao47}}{=}
\sum_{\bf m} \tbinom{\bf m}{\bf n} x^{{\bf m}-{\bf n}}\mathsf{z}_{\bf m}[p],
\end{align}
where as usual $\binom{\bf m}{\bf n}$ $:=\prod_{i=0}^d \binom{m_i}{n_i}$,
with the understanding that $\binom{m_i}{n_i}=0$ unless $m_i\ge n_i$.
Hence in view of (\ref{ck67}), (\ref{ck40}) takes the form of
\begin{align}\label{temp26}
(\Gamma_x^*)_{\delta_{\bf m}}^\gamma=\left\{\begin{array}{cc}
\tbinom{\bf m}{\bf n} x^{{\bf m}-{\bf n}}&\mbox{provided}\;\gamma=\delta_{\bf n}\;\mbox{for some
{\bf n}}\\
0&\mbox{else}\end{array}\right\}.
\end{align}
Note that this is consistent with (\ref{ao46}), (\ref{ck38}), and (\ref{temp02}).

\medskip

We remark that (\ref{ap31}) translates to
\begin{align}\label{ap34}
\lefteqn{\mbox{it is {\it not} true that}\quad
(\Gamma_x^*)_{\beta}^{\delta_{\bf n}}}\nonumber\\
&=\left\{\begin{array}{cc}
\tbinom{\bf m}{\bf n} x^{{\bf m}-{\bf n}}&\mbox{provided}\;\beta=\delta_{\bf m}\;\mbox{for some
{\bf m}}\\
0&\mbox{else}\end{array}\right\}.
\end{align}
Indeed, (informally) testing the hypothetical identity in (\ref{ap31}) with $\pi$ 
$\in\mathbb{R}[[\mathsf{z}_3,\mathsf{z}_{\bf n}]]$ and appealing to (\ref{ck70}),
we would obtain from this identity that $(\Gamma_x^*\pi)[\lambda,p]$ $=\pi[\lambda,p(\cdot+x)]$.
Restricting to $\pi=\mathsf{z}_{\bf n}$ and appealing to the second item in
(\ref{ck74}) and to (\ref{ap33}) 
would give the identity in (\ref{ap34}). We will argue at the end of Subsection \ref{ss:improve}
that (\ref{ap34}) generically holds.


\subsection{Uniqueness of \texorpdfstring{$\Gamma_x$}{Gamma} given \texorpdfstring{$\Pi$}{Pi} and \texorpdfstring{$\Pi^{-}$}{Pi-}}\label{ss:oneGamma}

We claim that the random endomorphism $\Gamma_x$ of 
$\mathbb{R}[\mathsf{z}_3,\mathsf{z}_{\bf n}]$ is uniquely determined by 
$\Pi,\Pi_x$ $\in X[[\mathsf{z}_3,\mathsf{z}_{\bf n}]]$:
\begin{align}\label{ap04}
\Gamma_x\quad\mbox{is determined by}\;\Pi\;\mbox{and}\;\Pi_x\;\mbox{via (\ref{temp02})}.
\end{align}
Statement (\ref{ap04}) only relies on the algebraic rules (\ref{temp05}) \& (\ref{temp06}).
For later purpose, we note that by the uniqueness (\ref{ap04}), 
the identity (\ref{ck62}) yields for $\beta=0$
\begin{align}\label{ap68}
(\Gamma_x^*)_0^\gamma=\delta_0^\gamma.
\end{align}

\begin{proof}[Proof of \eqref{ap04}]
Since for $\rho>0$, the components of
both $\Pi$ and $\Pi_x$ are smooth space-time functions, we will use (\ref{temp02})
in form of
\begin{align}\label{ap05}
\partial^{\bf n}\Pi(x)=\Gamma_x^*\partial^{\bf n}\Pi_x(x)\quad\mbox{for all}\;{\bf n}.
\end{align}
Hence the argument for (\ref{ap04}) relies on the fact that the jet
$\{\partial^{\bf n}\Pi_x(x)\}_{\bf n}$ is rich enough. In fact, we shall establish
\begin{align}\label{ap01}
\tfrac{1}{{\bf n}!} \partial^{\bf n}\Pi_x(x)-\mathsf{z}_{\bf n}
\in\mathbb{R}[[\mathsf{z}_3,\{\mathsf{z}_{\bf m}\}_{|{\bf m}|<|{\bf n}|}]],
\end{align}
where the space denotes the (sub-)algebra of formal power series in the finitely
many variables $\mathsf{z}_3$ and $\{\mathsf{z}_{\bf m}\}_{|{\bf m}|<|{\bf n}|}$.
In view of our postulate (\ref{ck52}) (in conjunction with (\ref{temp30}) to pass
to the general base-point $x$) and the smoothness of $\Pi_{x\beta}$, we have
\begin{align}\label{t10}
\partial^{\bf n}\Pi_{x\beta}(x)=0\quad\mbox{for}\;|\beta|>|{\bf n}|.
\end{align}
According to (\ref{temp24}) and (\ref{temp04bis}), the case $|\beta|=|{\bf n}|$
reduces to $\beta=\delta_{\bf m}$ for some ${\bf m}$, so that (\ref{ck38}) implies the sharpening
\begin{align*}
\tfrac{1}{{\bf n}!} \partial^{\bf n}\Pi_{x\beta}(x)=\delta_{\beta}^{\delta_{\bf n}}\quad\mbox{that is}\quad
(\tfrac{1}{{\bf n}!}\partial^{\bf n}\Pi_x(x)-\mathsf{z}_{\bf n})_{\beta}=0\quad\mbox{for}\;|\beta|\ge|{\bf n}|.
\end{align*}
To obtain (\ref{ap01}) it remains to realize that by (\ref{ap41}),
$|\beta|<|{\bf n}|$ implies that $\beta({\bf m})=0$ unless $|{\bf m}|<|{\bf n}|$.

\medskip

Equipped with (\ref{ap05}) and (\ref{ap01}),
the statement (\ref{ap04}) is established by induction in $|{\bf n}|$, starting with the
base case of ${\bf n}={\bf 0}$.
From (\ref{ap05}) and (\ref{ap01}) for ${\bf n}={\bf 0}$ together with
(\ref{temp06}) we learn that $\Gamma_x^*\mathsf{z}_{\bf 0}$ is determined.
Hence by (\ref{temp05}) and (\ref{temp06}), $\Gamma_x^*$ is determined
on any monomial in the variables $\mathsf{z}_3,\mathsf{z}_{\bf 0}$.
Thanks to the finiteness properties (\ref{ck88}), this determines $\Gamma_x^*$ 
on $\mathbb{R}[[\mathsf{z}_3,\mathsf{z}_{\bf 0}]]$. 
We now turn to the induction step, giving ourselves an ${\bf n}$ with $|{\bf n}|\ge 1$.
Once more from (\ref{ap05}) and (\ref{ap01}), this time in conjunction with
the induction hypothesis, we see that $\Gamma_x^*\mathsf{z}_{\bf n}$ is determined.
Hence together with the induction step, it is determined on the coordinates
in $\mathsf{z}_3$ and $\{\mathsf{z}_{\bf m}\}_{|{\bf m}|\le|{\bf n}|}$.
The outcome of the induction is that $\Gamma_x^*$ is determined on
all the coordinates $\mathsf{z}_3$ and $\{\mathsf{z}_{\bf m}\}$.
Again by multiplicativity and finiteness, this determines $\Gamma_x^*$ and thus $\Gamma_x$.
\end{proof}
 

\subsection{Population of \texorpdfstring{$\Gamma_x^*$}{Gamma}}\label{ss:triang}

Recalling the language from (\ref{ck54}), we claim that $\Gamma_x^*$ is sparse in the sense that 
\begin{align}\label{cw47}
\mbox{for populated $\gamma$}:\quad
(\Gamma^*_x)_{\beta}^{\gamma}&=0
\quad\mbox{unless $\beta$ is populated}.
\end{align}
We shall split this into the two sharper statements
that distinguish between purely polynomial $\gamma$ 
and those of the form (\ref{ap65}):
\begin{align}
\mbox{for $\gamma$ pp or $[\gamma]\ge 0$}:
\quad(\Gamma^*_x)_{\beta}^{\gamma}=0\quad
\mbox{unless $\beta$ pp or $[\beta]\ge 0$},\label{cw46}\\
\mbox{for $\gamma\in$ (\ref{ap65})}:
\quad(\Gamma^*_x)_{\beta}^{\gamma}=0\quad
\mbox{unless $\beta\in$ (\ref{ap65}) or $[\beta]\ge 0$}.\label{cw46ter}
\end{align}
We shall establish (\ref{cw46ter}) alongside the following extension of (\ref{temp26})
\begin{align}
\mbox{for $\gamma\in$ (\ref{ap65})}:
\quad(\Gamma^*_x)_{\beta}^{\gamma}=\sigma_\beta^\gamma x^{{\bf m}-{\bf n}}
\quad\mbox{unless $[\beta]\ge 0$}\label{cw46bis}
\end{align}
where, in line with (\ref{temp17}), ${\bf m}:=\sum_{{\bf m}'}{\bf m}'\beta({\bf m}')$,
${\bf n}:=\sum_{{\bf n}'}{\bf n}'\gamma({\bf n}')$,
and $\sigma_\beta^\gamma$ is some (deterministic) combinatorial factor
that vanishes unless ${\bf m}\ge{\bf n}$ (which means $m_i\ge n_i$ for $i=0,\cdots,d$).

\medskip

\begin{proof}[Proof of \eqref{cw46} \textnormal{\&} \eqref{cw46ter} \textnormal{\&} \eqref{cw46bis}]
Appealing to (\ref{temp15}), we argue by induction in $|\beta|$.
According to (\ref{temp08}), $\beta=0$ is the (only) base case. 
Since $\beta=0$ satisfies $[\beta]=0$, 
(\ref{cw46}) \& (\ref{cw46ter}) \& (\ref{cw46bis}) are automatically satisfied. 
We now turn to the induction step and note that for $\gamma=0$,
in view of the last item in (\ref{temp05}),
(\ref{cw46}) is trivially satisfied while (\ref{cw46ter}) \& (\ref{cw46bis}) are empty.
Hence we consider $\gamma\not=0$ and write it as
\begin{align}\label{cw71bis}
\gamma=k\delta_3+\sum_{j=1}^l\delta_{{\bf n}_j}\quad\mbox{for}\;(k,l)\not=(0,0)
\end{align}
for some $\{{\bf n}_j\}_{j=1,\dots,l}$.
We distinguish the three cases $[\gamma]\ge 0$, $\gamma\in$ (\ref{ap65}), and $\gamma$ pp:
\begin{align}
& 2k-l\ge 0\quad\mbox{and thus}\;k\not=0,\label{cw71}\\
& k=1\quad\mbox{and}\quad l=3,\label{cw71ter}\\
& k=0\quad\mbox{and}\quad l=1.\label{ap07}
\end{align}

\medskip

Since by (\ref{ck44}), (\ref{cw71bis}) translates into
$\mathsf{z}^\gamma$ $=\mathsf{z}_3^{k}\mathsf{z}_{{\bf n}_1} \cdots
\mathsf{z}_{{\bf n}_l}$, we obtain
by (\ref{ao45}), (\ref{temp05}) \& (\ref{temp06}), and (\ref{ck74})
\begin{align}\label{ap09}
(\Gamma_x^*)^\gamma_\beta=\sum_{k\delta_3+\beta_1+\cdots+\beta_l=\beta}
(\Gamma_x^*)^{\delta_{{\bf n}_1}}_{\beta_1}\cdots(\Gamma_x^*)^{\delta_{{\bf n}_l}}_{\beta_l},
\end{align}
with the understanding that the empty sum equals $0$ and the empty product equals $1$.
We now consider a summand in (\ref{ap09}); by (\ref{temp08}) we have
\begin{align*}
(|k\delta_3|-|0|)+\sum_{j=1}^l(|\beta_j|-|0|)=|\beta|-|0|\quad\mbox{and}\quad
|\beta_j|-|0|\ge 0.
\end{align*}
In the cases (\ref{cw71}) \& (\ref{cw71ter}) we have $k\not=0$ thus $|k\delta_3|-|0|>0$. 
Therefore $|\beta_j|<|\beta|$ for all $j=1,\cdots,l$. 
Hence we may appeal to the induction hypothesis (\ref{cw46})
for the factors $(\Gamma_x^*)^{\delta_{{\bf n}_j}}_{\beta_j}$: 
they vanish unless $[\beta_j]\ge 0$ or $\beta_j$ is pp, 
which in view of (\ref{ao10}) implies that the summand in (\ref{ap09}) vanishes unless
\begin{align}\label{ap67}
(2k-l)+\sum_{j=1}^l([\beta_j]+1)=[\beta]\quad\mbox{and}\quad
[\beta_j]+1\ge0.
\end{align}
In the case of (\ref{cw71}) we thus have $[\beta]\ge 0$, as desired. 
In the case of (\ref{cw71ter}) we have either $[\beta]\ge 0$,
in which case we are done, or $[\beta_j]=-1$ for $j=1,2,3$, 
in which case we have by induction hypothesis (\ref{cw46}) that
$\beta_j=\delta_{{\bf m}_j}$ for some ${\bf m}_j$ for $j=1,2,3$. 
Hence $\beta$ is of the form (\ref{ap65}), in line with (\ref{cw46ter}).
Moreover, in this case by (\ref{temp26}) we have 
$(\Gamma_x^*)_{\delta_{{\bf m}_j}}^{\delta_{{\bf n}_j}}$
$=\tbinom{{\bf m}_j}{{\bf n}_j} x^{{\bf m}_j-{\bf n}_j}$. In view of ${\bf m}=\sum_{j=1}^3{\bf m}_j$ and ${\bf n}=\sum_{j=1}^3{\bf n}_j$ this implies (\ref{cw46bis}).

\medskip

We now turn to the $\gamma$'s of the form (\ref{ap07}),
and rewrite (\ref{temp02}) component-wise,
cf.~(\ref{ck74}), as $\Pi_\beta$ $=\sum_\gamma(\Gamma_x^*)^\gamma_\beta\Pi_{x\gamma}$.
We split the sum according to whether $\gamma$ is purely polynomial, on which we use
(\ref{ck38}), or not:
\begin{align}\label{ap08}
p:=\sum_{\bf n}(\Gamma_x^*)_\beta^{\delta_{\bf n}}(\cdot-x)^{\bf n}
=\Pi_\beta-\sum_{\gamma\;\textnormal{not pp}}(\Gamma_x^*)_\beta^\gamma\Pi_{x\gamma}.
\end{align}
According to (\ref{ao08}), $\Pi_\beta$ vanishes unless $\beta$ is pp
or $[\beta]\ge 0$. By analogy, the factor
$\Pi_{x\gamma}$ vanishes unless $\gamma$ satisfies $[\gamma]\ge 0$ .
According to what we just showed in case \eqref{cw71}, the factor $(\Gamma_x^*)_\beta^\gamma$
vanishes unless $\beta$ is pp or $[\beta]\ge 0$. Hence the r.~h.~s.~of (\ref{ap08}) 
vanish unless $\beta$ is of this type. 
Hence also the (random) polynomial $p$ on the l.~h.~s., and thus all its coefficients, 
vanish unless $\beta$ is of this type. 
This establishes the induction step of (\ref{cw46}) for $\gamma$ of the 
remaining case (\ref{ap07}).
\end{proof}


\subsection{Strict triangularity of \texorpdfstring{$\Gamma_x^*$}{Gamma}}\label{ss:triang2}

Equipped with the results of Subsection \ref{ss:triang},
we shall establish that $\Gamma_x^*$ is strictly triangular w.~r.~t.~$|\cdot|$:
\begin{align}
\mbox{for all $\gamma$}:\quad
(\Gamma^*_x-{\rm id})_\beta^\gamma&=0\quad\mbox{unless $|\gamma|<|\beta|$}.\label{cw44}
\end{align}
Incidentally, the triangularity (\ref{cw44}) w.~r.~t.~$|\cdot|$ and the coercivity
(\ref{temp15}) of the latter imply the finiteness property (\ref{ck88}).

\medskip

As a consequence of (\ref{cw44}) and (\ref{temp15}), $\Gamma_x$ is invertible,
and thus a linear automorphism of $\mathbb{R}[\mathsf{z}_k,\mathsf{z}_{\bf n}]$.
As a consequence, $\Gamma_x^*$ is an algebra automorphism of
$\mathbb{R}[[\mathsf{z}_k,\mathsf{z}_{\bf n}]]$. This prompts us to introduce
\begin{align}\label{ap45}
\Gamma_{xx'}:=\Gamma_{x'}\Gamma_x^{-1}\quad\mbox{and thus}\quad
\Gamma_{xx'}^*=\Gamma_x^{-*}\Gamma_{x'}^*,\;\;
\Gamma_{xx'}^*\Gamma_{x'x''}^*=\Gamma_{xx''}^*.
\end{align}
Then (\ref{temp02}) \& (\ref{temp01}) extend to
\begin{align}\label{ap34bis}
\Pi_{x}=\Gamma_{xx'}^*\Pi_{x'}\quad\mbox{and}\quad\Pi_{x}^{-}=\Gamma_{xx'}^*\Pi_{x'}^{-}.
\end{align}
Compare now the first item of (\ref{ap34bis}) in form of
$\Pi_{x}=\Gamma_{x\;x+y}^*\Pi_{x+y}$ with (\ref{temp02}) in form of
$\Pi[\xi(\cdot+x)]=\Gamma_{y}^*[\xi(\cdot+x)]\Pi_{y}[\xi(\cdot+x)]$.
By the uniqueness statement of Subsection \ref{ss:oneGamma},
the shift covariance (\ref{ap32}) of $\Pi$ implies the shift covariance
of $\Gamma_{xx'}$, that is, $\Gamma_{y}^*[\xi(\cdot+x)]$ $=\Gamma_{x x+y}^*[\xi]$.
Together with the stationarity assumption (\ref{ck19}) this yields
\begin{align}\label{ap47}
\Gamma_{x x+y}^*=_{law}\Gamma_y^*.
\end{align}
It is easy to show that all the $\Gamma_{x x'}$
lie in a group that is characterized by (\ref{temp05}), (\ref{temp06}),
(\ref{temp26}) for some shift vector $x$, (\ref{cw47}), and (\ref{cw44}).
This is the \textit{structure group} \cite[Definition~2.1]{Hai15} of regularity structures; we will not make
any explicit use of it.
The data of $\{\Pi_x, \Pi_{x}^-, \Gamma_{x x'}\}_{x,x'}$ is called a \textit{model} in Hairer's language, see \cite[Definition~2.5]{Hai15}.

\medskip

\begin{proof}[Proof of \eqref{cw44}]
We establish (\ref{cw44}) by induction in $|\beta|$.
We start with the base case $\beta=0$:
By (\ref{temp08}), $|\gamma|\le|0|$ implies $\gamma=0$,
so that (\ref{cw44}) amounts to (\ref{ap68}).
We now turn to the induction step and distinguish the cases
\begin{align}\label{ap69}
\gamma\not=\mbox{pp}\quad\mbox{and}\quad
\gamma=\mbox{pp}.
\end{align}
We first tackle $\gamma\neq\mbox{pp}$ by an algebraic argument,
and then treat (II) by an analytic argument.
This structure foreshadows the proceeding in the proof
of Theorem \ref{th:main}, see Subsections \ref{ss:algIII} and \ref{ss:3pIII}.

\medskip

We start with $\gamma\neq\mbox{pp}$;
the special cases $\gamma=0$ and $\gamma=\delta_3$ are trivial because of
the second item in (\ref{temp05}) and of (\ref{temp06}),
respectively. 
It remains to treat $\gamma$ is with $\gamma(3)+\sum_{\bf n}\gamma({\bf n})\ge 2$.
These we may split
$\gamma=\gamma_1+\gamma_2$ with $\gamma_j\not=0$ for $j=1,2$.
We obtain as for (\ref{ap09})
\begin{align}\label{cw57}
(\Gamma_x^*)^\gamma_\beta=\sum_{\beta_1+\beta_2=\beta}
(\Gamma_x^*)^{\gamma_1}_{\beta_1}(\Gamma_x^*)^{\gamma_2}_{\beta_2}.
\end{align}
A summand in (\ref{cw57}) vanishes unless $\beta_j\not=0$ for $j=1,2$
since otherwise by the base case we would have $\gamma_j=0$,
which is ruled out by the above splitting.
Hence in view of (\ref{temp08}) we have
\begin{align*}
(|\beta_1|-|0|)+(|\beta_2|-|0|)=|\beta|-|0|\;\;\mbox{and}\;\;
(|\beta_1|-|0|),\,(|\beta_2|-|0|)>0,
\end{align*}
so that in particular $|\beta_1|,|\beta_2|$ $<|\beta|$. This allows us to appeal
to the induction hypothesis for the two factors in (\ref{cw57}): They vanish unless
$|\gamma_j|\le|\beta_j|$ for $j=1,2$. Since by (\ref{temp08}) we also have
\begin{align*}
(|\gamma_1|-|0|)+(|\gamma_2|-|0|)=|\gamma|-|0|;
\end{align*}
this yields that the product vanishes unless $|\gamma|\le|\beta|$.
Moreover, in case of $|\gamma|=|\beta|$ we must have $|\gamma_j|=|\beta_j|$ for $j=1,2$.
By induction hypothesis this implies that necessarily $\gamma_j=\beta_j$ for $j=1,2$
and that both factors in (\ref{cw57}) are $=1$, and thus the product $=1$.
Hence the sum in (\ref{cw57}) consists of a single summand $=1$ and thus
is $=1$; finally, we must have $\gamma$ $=\gamma_1+\gamma_2$ $=\beta_1+\beta_2$ $=\beta$.

\medskip

We now turn to $\gamma=\mbox{pp}$ in (\ref{ap69}) and reconsider (\ref{ap08}).
It follows from (\ref{cw44}), which we just established for
the current $\beta$ and all $\gamma$ not pp, that the r.~h.~s.~of (\ref{ap08})
effectively only involves $\gamma$'s with $|\gamma|\le|\beta|$.
We apply the mollification $(\cdot)_r$ to (\ref{ap08}), evaluate at $x$, and
use the triangle and Cauchy-Schwarz inequalities in probability for
\begin{align*}
\mathbb{E}|p_r(x)|
&\le\mathbb{E}|\Pi_{\beta r}(x)|+\sum_{\gamma\colon |\gamma|\le|\beta|}
\mathbb{E}^\frac{1}{2}|(\Gamma_x^*)_\beta^\gamma|^2 \,
\mathbb{E}^\frac{1}{2}|\Pi_{x\gamma r}(x)|^2.
\end{align*}
By the argument at the end of Subsection \ref{ss:trans} we obtain
$\limsup_{r\uparrow\infty}$ $r^{-|\beta|}$ $\mathbb{E}|\Pi_{\beta r}(x)|<\infty$ from (\ref{ck52}).
Strengthening this postulate (\ref{ck52}) to holds with $\mathbb{E}|\cdot|$
replaced by $\mathbb{E}^\frac{1}{2}|\cdot|^2$, we obtain by (\ref{temp30}) that also
$\limsup_{r\uparrow\infty}r^{-|\gamma|}\mathbb{E}^\frac{1}{2}|\Pi_{x\gamma r}(x)|^2<\infty$.
Combining this with the purely qualitatively postulate that $\mathbb{E}^\frac{1}{2}
|(\Gamma^*_x)_\beta^\gamma|^2<\infty$ we thus obtain
$\limsup_{r\uparrow\infty}$ $r^{-|\beta|}$ $\mathbb{E}|p_{r}(x)|<\infty$.
By the argument at the end of Subsection \ref{ss:unique} we conclude that $p$ has
degree $\le|\beta|$ to the effect of
\begin{align}\label{ap11}
(\Gamma_x^*)_\beta^{\delta_{\bf n}}=0\quad\mbox{unless}\;|{\bf n}|\le|\beta|.
\end{align}

\medskip

In the case of equality $|{\bf n}|=|\beta|$ in (\ref{ap11}),
we note that by (\ref{cw46}), the l.~h.~s.~vanishes unless $[\beta]\ge 0$ or $\beta$ is pp.
Hence by (\ref{temp04bis}) we must have that $\beta$ is pp.
We then appeal to (\ref{temp26}) to learn that (\ref{ap11}) upgrades to
\begin{align*}
(\Gamma_x^*-{\rm id})_\beta^{\delta_{\bf n}}=0\quad\mbox{unless}\;|{\bf n}|<|\beta|.
\end{align*}
In view of (\ref{temp04}), we thus established (\ref{cw44}) in the remaining case of $\gamma=\mbox{pp}$.
\end{proof}


\section{Main result and sketch of proof} \label{section:main}

\subsection{What estimates to expect?}\label{ss:expect}

We now continue exploring
the consequences of combining the shift and scaling invariances (\ref{ck19}) and (\ref{ck22}) of $\xi$.
Recall that we argued in Subsection~\ref{ss:unique} that 
in the limit $\rho\downarrow0$, 
the laws of $r^{-|\beta|}\Pi_{\beta r}(0)$ and $r^{2-|\beta|}\Pi_{\beta r}^{-}(0)$ 
are independent of $r$, see (\ref{ck27bis}).

\medskip 

Likewise, we now discuss an analogous scaling invariance for $\Gamma^*$.
Clear\-ly, the transformation (\ref{ck66}) depends on $\xi$
so that $\Gamma_x$ is a random endomorphism of $\mathbb{R}[\mathsf{z}_3,\mathsf{z}_{\bf n}]$,
and thus the matrix elements $(\Gamma^*_x)_\beta^\gamma$ are random numbers.
We claim that the invariances in law (\ref{ck27}) and (\ref{temp30}) of $\Pi_x$
are consistent with
\begin{align}\label{ck76}
(\Gamma_{Rx}^*)_\beta^\gamma=_{law}r^{|\beta|-|\gamma|}(\Gamma_x^*)_\beta^\gamma
\quad\mbox{in the limit}\;\rho\downarrow 0.
\end{align}
Note that for purely polynomial $\beta$, 
this is consistent with (\ref{temp26})
in view of (\ref{temp04}).

\medskip

Let us argue that \eqref{ck76} is reasonable in view of (\ref{ck27}) and (\ref{temp30}).
Indeed, by the component-wise version of (\ref{temp02})
\begin{align}\label{ck75}
\Pi_\beta=\sum_{\gamma}(\Gamma_x^*)_\beta^\gamma\Pi_{x\gamma},
\end{align}
which we first use with $x$ replaced by $Rx$, and evaluated at $R(\cdot+x)$
\begin{align*}
\Pi_\beta(R(\cdot+x))=\sum_{\gamma}(\Gamma_{Rx}^*)_\beta^\gamma\Pi_{Rx\;\gamma}(R(\cdot+x)).
\end{align*}
Combining (\ref{ck27}) and (\ref{temp30}) we have
\begin{align*}
\Pi_{Rx\;\gamma}(R(\cdot+x))=_{law}\Pi_{\gamma}(R\cdot)=_{law}r^{|\gamma|}\Pi_{\gamma}
=_{law}r^{|\gamma|}\Pi_{x\;\gamma}(\cdot+x).
\end{align*}
Since we think of this and (\ref{ck76}) as holding jointly, we thus learn
\begin{align*}
\Pi_\beta(R(\cdot+x))=_{law}
r^{|\beta|}\sum_{\gamma}(\Gamma_{x}^*)_\beta^\gamma\Pi_{x\;\gamma}(\cdot+x)
\stackrel{\eqref{ck75}}{=}r^{|\beta|}\Pi_\beta(\cdot+x),
\end{align*}
which in turn is consistent with (\ref{ck27}).

\medskip

It is therefore reasonable to expect by \eqref{ck27bis} that even before 
passing to the limit, we have the estimates
\begin{align*}
\mathbb{E}^\frac{1}{p}|\Pi_{\beta r}(0)|^p\lesssim r^{|\beta|}\quad\mbox{and}\quad
\mathbb{E}^\frac{1}{p}|\Pi_{\beta r}^{-}(0)|^p\lesssim r^{|\beta|-2},
\end{align*}
where the implicit constant depends on $\beta$, the arbitrary exponent $p$ $<\infty$,
and the semi-norms (\ref{ck52a}) of the kernel $\psi$ in the definition (\ref{ck30})
of the mollification operator $(\cdot)_r$ -- but not on $\rho$. By (\ref{temp30}),
this extends to 
\begin{align}\label{ck77}
\mathbb{E}^\frac{1}{p}|\Pi_{x\beta r}    (x)|^p\lesssim r^{|\beta|}\quad\mbox{and}\quad
\mathbb{E}^\frac{1}{p}|\Pi_{x\beta r}^{-}(x)|^p\lesssim r^{|\beta|-2}.
\end{align}

\medskip

Similarly, in view of (\ref{ck76})
the law of $r^{|\gamma|-|\beta|}(\Gamma_{Rx}^*)_\beta^\gamma$
is independent of $r$. This suggests that we may hope for the estimate
\begin{align}\label{ck78}
\mathbb{E}^\frac{1}{p}|(\Gamma_{x}^*)_\beta^\gamma|^p\lesssim|x|^{|\beta|-|\gamma|},
\end{align}
where
\begin{align}\label{cw17bis}
|x|:=\sqrt[4]{x_0^2+({\textstyle\sum_{i=1}^d}x_i^2)^2},
\end{align}
which behaves as the parabolic Carnot-Carath\'eodory distance $\sqrt{|x_0|}$ 
$+\sqrt{\sum_{i=1}^dx_i^2}$ between $x$ and the origin.

\medskip

We now motivate an assumption on ensembles $\xi$ for which we are able to establish the estimates \eqref{ck77} \& (\ref{ck78}).
Let us start by discussing an example which satisfies \eqref{ck22} as well as (\ref{ck19}) - (\ref{ck21}).
It is given by the centered Gaussian ensemble on Schwartz distributions characterized by
\begin{align}\label{ck81}
\mathbb{E}(\xi,\zeta)(\xi,\zeta')=\mbox{$\dot H^{-s}$-inner product of $\zeta,\zeta'$},
\end{align}
where the homogeneous Sobolev space $\dot H^s$ of (fractional and possibly negative parabolic)
order $s$ of $\zeta$ is conveniently defined in terms of the Fourier transform
$({\mathcal F}\zeta)(q)$ $=\int\frac{dx}{(2\pi)^{1+d}} e^{-ix\cdot q}\zeta(x)$ via
\begin{align}\label{ck81ter}
\mbox{$\dot H^{s}$-inner product of $\zeta,\zeta'$}
=\int dq \, |q|^s{\mathcal F}\zeta \, \overline{|q|^s{\mathcal F}\zeta'},
\end{align}
where, in line with (\ref{cw17bis}),
\begin{align}\label{ap37}
|q|:=\sqrt[4]{q_0^2+({\textstyle\sum_{i=1}^d}q_i^2)^2}
\end{align}
is a size measure on wave vectors $q$ that scales as $r^{-1}$ under the parabolic
rescaling $R^{-1}q$, thus ensuring that (\ref{ck81ter}) is of (parabolic) order $s$.
In case of $d=0$, this Gaussian ensemble coincides with $\partial_0B_{x_0}$
where $B_{x_0}$ is a fractional Brownian motion of Hurst exponent\footnote{the factor 2 here is due to the parabolic nature of \eqref{ap37} which reduces to $| q | = \sqrt[4]{q^2}$ when $d = 0$} $2(\frac{1}{2}+s)$.
Note that $|q|^4$ is the symbol of the fourth-order operator $LL^*$.
For $s>0$, the $\dot H^{-s}$-norm is obviously weaker than the $L^2$-norm on small scales,
so that (\ref{ck81}) implies that the realizations of $\xi$ are less rough than white
noise, while for $s<0$ this implies that the realizations of $\xi$ are rougher than white noise.
\medskip

In fact, our main result establishes the estimates \eqref{ck77} \& (\ref{ck78}) under
the inequality version of (\ref{ck81}), that is
\begin{align}\label{ck81bis}
\mbox{variance of}\;(\xi,\zeta)
\le\big(\mbox{$\dot H^{-s}$ norm of $\zeta$}\big)^2,
\end{align}
cf.~(\ref{ck81ter}). This means that the Cameron-Martin space of this Gaussian
ensemble, interpreted as a space of functions rather than their dual (i.~e.~~distributions), 
is contained in the $L^2$-dual of $\dot H^{-s}$, which is $\dot H^s$. 
In the course of the proof, we will relax the assumption of Gaussianity.


\subsection{Statement of main result}\label{ss:main}

As long as $\rho>0$, the above construction of $\Pi_x$, $\Pi^{-}_x$,
and $\Gamma^*_x$ can actually be carried out rigorously.
Our main result is that these objects are estimated uniformly
in $\rho\downarrow 0$, while $c$ typically diverges.

\begin{theorem}\label{th:main}
Suppose that the centered Gaussian ensemble on Schwartz distributions on $\mathbb{R}^{1+d}$
satisfies the symmetries \eqref{ck19} -- \eqref{ck21},
and satisfies \eqref{ck81bis} for some $s$ with \eqref{ck53} and \eqref{ck28}.
Moreover, we assume 
$D \geq 3$ and
\begin{align}\label{ap43}
3\alpha+\frac{D}{2}>0.
\end{align}
Then we have 
\begin{align} 
\mathbb{E}^{\frac{1}{p}}|\Pi_{x \beta r}    (x)|^p 
&\lesssim r^{|\beta|} , 
\label{ao65impov}\\
\mathbb{E}^{\frac{1}{p}}|\Pi_{x \beta r}^{-}(x)|^p
&\lesssim r^{|\beta|-2} ,
\label{ao65bisimpov}\\
\mathbb{E}^{\frac{1}{p}}|(\Gamma_{x}^*)_{\beta}^{\gamma}|^p&\lesssim|x|^{|\beta|-|\gamma|}
\quad\mbox{for all populated $\gamma$}.
\label{ao65ter}
\end{align}
The implicit multiplicative constants only depend on $d$, $s$, $\beta$, $\gamma$, $p$,
and on $\psi$ only through the semi-norms, but not on $x$, $r$, and $\rho$.
\end{theorem}

The motivation for assumption (\ref{ap43}) will be given in the next Subsection
\ref{ss:usage}, see (\ref{ap44}); it is empty for $D\ge 6$.
It can also be shown that 
$\Pi$, $\Pi^{-}$, and $\Gamma$ converge as $\rho\downarrow0$
to a uniquely characterized limit, see \cite{Tem}.
This unique characterization involves the Malliavin derivative, which will be
introduced in Subsection~\ref{ss:usage}.

\medskip

\begin{remark}
The estimates \eqref{ao65impov} and \eqref{ao65bisimpov} are still impoverished, 
because the center of the average
agrees with the base-point. 
Here $\Gamma^*$ comes to help, 
which allows to post-process 
\eqref{ao65impov} into the stronger estimate
\begin{align}\label{ao65}
\mathbb{E}^{\frac{1}{p}} | \Pi_{\beta r} ( x ) |^p 
\lesssim r^{\alpha} ( r + | x | )^{| \beta | - \alpha} .
\end{align}
Similarly, \eqref{ao65bisimpov} can be upgraded to the stronger estimate
\begin{align}\label{ao65bisimproved}
\mathbb{E}^{\frac{1}{p}} | \Pi_{\beta r}^- ( x ) |^p 
\lesssim r^{3 \alpha} ( r + | x | )^{| \beta | - 2 - 3 \alpha} , 
\quad \text{ provided $\beta \neq 0$},
\end{align}
with the understanding that the l.~h.~s.~ vanishes unless $|\beta|\geq2+3\alpha$. 
Note that for $\beta=0$ we have $\Pi^-_{x \, \beta=0} = \xi$, 
which is independent of $x$ anyway, see \eqref{ck62}. 
\end{remark}

We note that \eqref{ao65} contains several pieces of information:
The local degree of regularity of $\Pi_{\beta}$ is of the (negative) order $\alpha$;
however, in $x=0$ $\Pi_{\beta}$ (on average) vanishes to order $|\beta|\geq\alpha$;
finally, at infinity $\Pi_{\beta}$ grows (on average) at order $|\beta|-\alpha$.
The first exponent in \eqref{ao65bisimproved} is expected, 
since on a heuristic level the cubic terms in the hierarchy \eqref{ao44} have regularity $3\alpha$.
We note that the appearance of the exponent $3\alpha$ already points towards \eqref{ap43}.

\begin{proof}[Proof of \eqref{ao65} \textnormal{\&} \eqref{ao65bisimproved}]
We give the proof of \eqref{ao65bisimproved}, 
the proof of \eqref{ao65} proceeds analogously.
We may appeal to \eqref{temp01} in its component-wise version, 
to which we apply the convolution operator
$(\cdot)_r$ from \eqref{ck30} and which we evaluate at $x$; 
using the triangle inequality w.~r.~t.~the norm 
$\mathbb{E}^\frac{1}{p}|\cdot|^p$ and then H\"older's inequality we obtain
\begin{align*}
\mathbb{E}^\frac{1}{p}|\Pi^-_{\beta r}(x)|^p
\le\sum_\gamma\mathbb{E}^\frac{1}{2p}|(\Gamma_x^*)_\beta^\gamma|^{2p}
\, \mathbb{E}^\frac{1}{2p}|\Pi^-_{x\gamma r}(x)|^{2p}.
\end{align*}
We now may paste in \eqref{ao65bisimpov} and \eqref{ao65ter} with $p$ replaced by $2p$,
and appeal to (\ref{cw44}) in order to obtain
$\mathbb{E}^\frac{1}{p}|\Pi^-_{\beta r}(x)|^p$
$\lesssim\sum_{|\gamma|\le|\beta|}|x|^{|\beta|-|\gamma|}r^{|\gamma|-2}$, 
which by \eqref{temp15} 
can be consolidated to \eqref{ao65bisimproved} 
provided $|\gamma|\geq2+3\alpha$.
Indeed, from \eqref{temp05} and $\beta\neq0$ we deduce that 
$(\Gamma^*_x)_\beta^0=0$, 
and thus in the above expansion effectively $\gamma\neq0$. 
From the population condition \eqref{temp17} on $\Pi^{-}$, $\gamma$ cannot be purely polynomial either.
Thus $\gamma (3) \geq 1$ and in turn from \eqref{ao07} we deduce $| \gamma | \geq 2 + 3 \alpha$ as claimed.
\end{proof}


\subsection{Population of \texorpdfstring{$\Gamma^*$}{Gamma} and estimates of \texorpdfstring{$\Pi_\beta$}{Pi beta} for \texorpdfstring{$\beta=\delta_3+\delta_{\bf 0},2\delta_3+\delta_{\bf 0}$}{beta = delta 3 + delta 0, 2 delta 3 + delta 0} revisited}\label{ss:improve}

For later purpose we remark that \eqref{ao65bisimproved} is still too pessimistic for certain classes of 
multi-indices, the simplest among those are $\beta=\delta_3+\delta_{\bf 0},2\delta_3+\delta_{\bf 0}$. 
To this aim, we first note that for any $k\ge 0$
\begin{align}
\{\,\gamma\;\mbox{populated}\,&|\,(\Gamma_x^*-{\rm id})_{\delta_3+k\delta_{\bf 0}}^\gamma\not=0\,\}  \label{ap35} \\
&\subset\{\mbox{purely polynomial}\} , \nonumber \\
\{\,\gamma\;\mbox{populated}\,&|\,(\Gamma_x^*-{\rm id})_{2\delta_3+\delta_{\bf 0}}^\gamma\not=0\,\} \label{ap35b} \\
&\subset \{\mbox{purely polynomial}\} \nonumber \\
&\,\quad\cup \{\delta_3+\mbox{purely polynomial}\} \nonumber\\
&\,\quad\cup \{\delta_3+\delta_{\bf 0} + \mbox{purely polynomial}\} . \nonumber
\end{align}
This is reminiscent of the notion of sector in regularity structures, which is related to the ``bare'' regularity of $\Pi_{\beta}$, see e.~g.~ \cite[Corollary~3.16]{Hai14}. 
\begin{proof}[Proof of \eqref{ap35} \textnormal{\&} \eqref{ap35b}]
We start with the proof of \eqref{ap35} and 
first note that from (\ref{temp05}) and (\ref{temp06}) we obtain $\gamma(3)\le 1$.
If $\gamma(3)=0$, we learn from (\ref{ck54}) that $\gamma$ 
must be pp or $\gamma=0$, but $(\Gamma^*_x-{\rm id})_\beta^0=0$ for all $\beta$ by the 
second item in \eqref{temp05}; thus $\gamma$
must be pp, as desired.
If $\gamma(3)=1$, that is, $\gamma=\delta_3+\sum_{j=1}^l\delta_{{\bf n}_j}$
we infer from (\ref{cw44}) that $l\le k$. Hence 
we obtain from (\ref{temp05}) \& \eqref{temp06} that 
$(\Gamma_x^*)_{\delta_3+k\delta_{\bf 0}}^\gamma$ 
$=\prod_{j=1}^l(\Gamma_x^*)_{\beta_j}^{\delta_{{\bf n}_j}}$ with
$\sum_{j=1}^l\beta_j=k\delta_{\bf 0}$. 
Since $k\delta_{\bf 0}$ is not populated unless $k\le 1$ we learn
from (\ref{cw46}) that $k=l$ and $\beta_j=\delta_{\bf 0}$. By (\ref{temp26}) this
yields ${\bf n}_j$ $={\bf 0}$. Hence we obtain $\gamma=\delta_3+k\delta_{\bf 0}$, 
and we conclude with \eqref{cw44}.

\medskip

We turn to the proof of \eqref{ap35b}.
From the second item of \eqref{temp05} and \eqref{temp06} we infer $\gamma\neq0,\delta_3$;
In view of \eqref{ck54} we thus may assume that $\gamma=\gamma_1+\gamma_2$ with $\gamma_{1,2}\neq0$.
From the first item of \eqref{temp05}, we obtain 
\begin{align*}
(\Gamma^*_x)_{2\delta_3+\delta_{\bf 0}}^\gamma
= \sum_{\beta_1+\beta_2=2\delta_3+\delta_{\bf 0}}
(\Gamma^*_x)_{\beta_1}^{\gamma_1}
(\Gamma^*_x)_{\beta_2}^{\gamma_2},
\end{align*}
where by \eqref{ap68} we have $\beta_{1,2}\neq0$. 
We first treat the case $\beta_1=\delta_{\bf 0}$ and $\beta_2=2\delta_3$. 
If this summand is non vanishing, 
then it follows from \eqref{temp26} that $\gamma_1=\delta_{\bf 0}$, 
and from \eqref{temp05} and \eqref{temp06} that $\gamma_2=2\delta_3$;
thus $\gamma=2\delta_3+\delta_{\bf 0}$ and we conclude by \eqref{cw44}.
We now treat the case $\beta_1=\delta_3$ and $\beta_2=\delta_3+\delta_{\bf 0}$,
the remaining cases are dealt with by symmetry. 
If this summand is non vanishing, 
it then follows from \eqref{ap35} 
that $\gamma_1=\delta_3$ or $\gamma_1$ is purely polynomial, 
and that $\gamma_2=\delta_3+\delta_{\bf 0}$ 
or $\gamma_2$ is purely polynomial.
If $\gamma_1=\delta_3$ and $\gamma_2=\delta_3+\delta_{\bf 0}$ we conclude again by \eqref{cw44}. 
If $\gamma_1$ and $\gamma_2$ are purely polynomial, then $\gamma$ is not populated. 
In the remaining two cases $\gamma$ is an element of the r.~h.~s.~of \eqref{ap35b}. 
\end{proof}

\medskip

Because of \eqref{cw44} and (\ref{ap35}), 
in the case of $\alpha<-1/2$, 
(\ref{temp02}) and (\ref{temp01}) collapse to
\begin{align}\label{ap36}
\Pi_{\delta_3+\delta_{\bf 0}}=\Pi_{x\;\delta_3+\delta_{\bf 0}}
+(\Gamma_x^*)_{\delta_3+\delta_{\bf 0}}^{\delta_{\bf 0}}\quad\mbox{and}\quad
\Pi_{\delta_3+\delta_{\bf 0}}^{-}=\Pi_{x\;\delta_3+\delta_{\bf 0}}^{-}.
\end{align}
Note that the second item in (\ref{ap36}) is of the same type as (\ref{ck62});
in view of (\ref{temp30}) it shows that both $\Pi_0$ and 
$\Pi_{\delta_3+\delta_{\bf 0}}^{-}$ are stationary.
From the second item in (\ref{ap36}) we learn that \eqref{ao65bisimpov} yields
\begin{align}\label{ck79bis}
\mathbb{E}^\frac{1}{p}|\Pi_{\delta_3+\delta_{\bf 0}\;r}^{-}(x)|^p
\lesssim r^{|\delta_3+\delta_{\bf 0}|-2},
\end{align}
which is an improvement over \eqref{ao65bisimproved} in the sense that it eliminates $x$
from the r.~h.~s..
Similarly, we can exploit \eqref{ap35b} in the argument that led to \eqref{ao65bisimproved} to obtain
\begin{equation}\label{t03}
\mathbb{E}^{\frac{1}{p}} | \Pi^-_{2 \delta_3 + \delta_{\bf 0} \, r} (x) |^p
\lesssim r^{2\alpha}(r+|x|)^{|2\delta_3+\delta_{\bf 0}|-2-2\alpha}, 
\end{equation}
which is an improvement over \eqref{ao65bisimproved} because of $\alpha\leq0$.
Now we are in a position to prove \eqref{ap30}.

\medskip
\begin{proof}[Proof of \eqref{ap30}]
From the first item in (\ref{ap36}) we actually learn that we cannot expect
that $(\Gamma_x^*)_{\delta_3+\delta_{\bf 0}}^{\delta_{\bf 0}}$ vanishes for
all $x$ and all realizations, which is the argument for (\ref{ap34}) and thus, working our
way back, for (\ref{ap31}) and ultimately (\ref{ap30}). Indeed, if it were vanishing,
(\ref{ao65impov}) would translate by (\ref{ap36}) into
$\mathbb{E}^\frac{1}{p}|\Pi_{\delta_3+\delta_{\bf 0}\,r}(x)|^p\lesssim r^{|\delta_3+\delta_{\bf 0}|=2\alpha+2}$
for all $x$, from which we learn by $r\downarrow 0$ that $\Pi_{\delta_3+\delta_{\bf 0}}$
vanishes, which would imply that $\Pi_{\delta_3+\delta_{\bf 0}}^{-}$ vanishes
by (\ref{ao43ter}). In view of (\ref{ck68}) this would imply that $\Pi_0^2$ is constant.
Since in view of (\ref{ao43ter}), $\Pi_0$ inherits from $\xi$ that it is a
centered Gaussian, $\Pi_0$ must vanish, and thus also $\xi_\rho$. 
\end{proof}


\subsection{Usage of the (directional) Malliavin derivative\texorpdfstring{ $\delta$}{}, first attempt}\label{ss:usage}

We will take the liberty to use the Malliavin derivative as a conceptual tool here in an informal fashion. 
The intuition best suited for our purposes is that
the Malliavin derivative of a random variable $F$ amounts to a Fr\'echet derivative 
of $F$ considered as a functional $F=F[\xi]$ of $\xi$. In fact, we will 
work with the directional Malliavin derivative: Given an element $\delta\xi$ of the 
Cameron-Martin space, which we think of as an infinitesimal
variation of $\xi$ (thus the notation), we consider the random variable 
\begin{align*}
\delta F[\xi]:=\frac{d}{dt}_{|t=0}F[\xi+t\delta\xi],
\end{align*}
which is the infinitesimal variation\footnote{not to be confused with the notion of divergence
in Malliavin calculus} of $F$ generated by $\delta\xi$.

\medskip
 
We will apply the derivation $\delta$ to $F=\Pi_{\beta r}^{-}(x)$,
which in view of (\ref{ao44}) and (\ref{ao43ter}) arises from $\xi$ 
by a sequence of operations that correspond to taking products and inverting\footnote{and 
thus is a multi-linear function in $\xi$, a property we make no explicit use of in this work, neither on the heuristic nor on the rigorous level} $L$.
Hence applying $\delta$ amounts to a linearization of 
these operations around a given $\xi$. Loosely speaking,
it monitors how the solution $\phi$ at fixed parameter $(\lambda,p)$
depends on the r.~h.~s.~$\xi$.
One aspect of Malliavin calculus proper that we will appeal to
in this heuristic discussion is that 
\begin{align}\label{ck72}
\delta\xi\in\dot H^s, 
\end{align}
see the discussion after (\ref{ck81bis}).

\medskip

The main challenge of renormalization is that (\ref{ao43}) encodes the map $\Pi\mapsto\Pi^{-}$ 
in a non-robust way as $\rho\downarrow 0$ since it contains the divergent $c$ and 
a singular (triple) product.
Applying $\delta$ to (\ref{ao43}) seems promising:
\begin{itemize}
\item Since $c$ is deterministic we have $\delta c=0$.
\item While $\xi$ has regularity of order $s-\frac{D}{2}$, cf.~(\ref{ck22}),
$\delta\xi$ has regularity of order $s$, cf.~(\ref{ck72}), an improvement 
by $\frac{D}{2}$ units.\footnote{Only in Subsection \ref{ss:Besov}, we will
start to worry that these orders of regularity are measured in different norms,
namely uniform/stationary vs.~space-time $L^2$.}
\item The control of the Malliavin derivative $\delta\Pi^{-}$, 
which captures the fluctuations\footnote{this is quantified by the variance control via the spectral gap inequality \eqref{ap18}} of the random variable $\Pi^{-}$, naturally complements 
the BPHZ-choice of renormalization (\ref{ck52ter}), which takes care of the expectation.
\end{itemize}
Indeed, in Subsection \ref{ss:robust}, we argue that there exists a robust
map $\delta\Pi\mapsto\delta\Pi^{-}$. 

\medskip

Unfortunately, things are a bit more complicated:
Applying $\delta$ to (\ref{ao43}) does not eliminate $c$, since we obtain from Leibniz' rule
\begin{align}\label{ck69}
\delta\Pi^{-}=(3\mathsf{z}_3\Pi^2+c)\delta\Pi+\delta\xi_\rho\mathsf{1};
\end{align}
applying $\delta$, which commutes with $L$, to (\ref{ao43ter}) we obtain
\begin{align}\label{ck70bis}
L\delta\Pi_\beta=\delta\Pi_\beta^{-} \quad\mbox{ mod polynomials of degree $\le|\beta|-2$}.
\end{align}

\medskip

To make things worse, the product in (\ref{ck69}) 
still is not robust in the limit $\rho\downarrow0$:
Consider the multi-indices $\beta=0,\delta_3,2\delta_3$; 
in view of $c_0=0$ and (\ref{ck68}), (\ref{ck69}) and (\ref{ck70bis}) yield for
the corresponding components
\begin{align*}
&L\delta\Pi_0=\delta\xi_\rho,\quad
L\delta\Pi_{\delta_3}=\Pi_{\delta_3+\delta_{\bf 0}}^{-}\delta\Pi_0,\quad\mbox{and}\\
&\delta\Pi_{2\delta_3}^{-}=\Pi_{\delta_3+\delta_{\bf 0}}^{-}\delta\Pi_{\delta_3}
+(6\Pi_0\Pi_{\delta_3}+c_2)\delta\Pi_{0}.
\end{align*}
From the first item we learn that (\ref{ck72}) translates into $\delta\Pi_0\in\dot H^{s+2}$.
Since $|\delta_3+\delta_{\bf 0}|$ $=2(1+\alpha)$, cf.~(\ref{ap42}),
we learn from (\ref{ck79bis}) that $\Pi_{\delta_3+\delta_{\bf 0}}^{-}$
has (limiting) regularity $2\alpha$. 
The sum of the regularities of the two factors is given by
\begin{align}\label{ap44}
2\alpha+(s+2)\stackrel{\eqref{ap40}}{=}3\alpha+\frac{D}{2}.
\end{align}
Since by assumption (\ref{ap43}), this sum is positive, it is well-known that the product 
$\Pi_{\delta_3+\delta_{\bf 0}}^{-}\delta\Pi_0$ converges in the limit $\rho\downarrow0$.
However, its regularity is obviously dominated by the worst factor, 
so that $\delta\Pi_{\delta_3}\in \dot H^{2\alpha+2}$ and no better. 

\medskip

We now turn to the product
$\Pi_{\delta_3+\delta_{\bf 0}}^{-}\delta\Pi_{\delta_3}$ .
According to the above, the sum of the regularities is given by 
$2\alpha+(2\alpha+2)$ $=4\alpha+2$. This expression is negative for $\alpha<-\frac{1}{2}$, 
and thus the product not robust in the limit $\rho\downarrow 0$. 
This is mirrored by the fact that
the other factor $6\Pi_0\Pi_{\delta_3}+c_2$ does not quite agree with
$\Pi^{-}_{2\delta_3+\delta_{\bf 0}}$, cf.~(\ref{ck74bis}), which would be controlled\footnote{
In fact, by (\ref{ck74bis})
this first factor is
$\Pi_{2\delta_3+\delta_{\bf 0}}^{-}-\Pi_{\delta_3+\delta_{\bf 0}}^{-}\Pi_{\delta_3+\delta_{\bf 0}}$;
while the first summand $\Pi_{2\delta_3+\delta_{\bf 0}}^{-}$
(and by \eqref{t03} and \eqref{ap43} even its product with $\delta\Pi_0$) stays finite 
as $\rho\downarrow 0$, the second term 
$\Pi_{\delta_3+\delta_{\bf 0}}^{-}$ $\Pi_{\delta_3+\delta_{\bf 0}}$
does not, since it can be rewritten as $L\frac{1}{2}\Pi_{\delta_3+\delta_{\bf 0}}^2$
$+\sum_{i=1}^d(\partial_i\Pi_{\delta_3+\delta_{\bf 0}})^2$: each of the 
$(\partial_i\Pi_{\delta_3+\delta_{\bf 0}})^2$ diverges since in view of (\ref{ck79bis}) the
bare regularity of
$\Pi_{\delta_3+\delta_{\bf 0}}$ is given by its homogeneity,
namely $|\delta_3+\delta_{\bf 0}|=2+2\alpha$, which is
$<1$ for $\alpha<-\frac{1}{2}$, again.}.
Hence we need to better explore the structure of $\delta\Pi$ and $\delta\Pi^{-}$.


\subsection{Tangent space to solution manifold}\label{ss:tangent}

To understand the structure of $\delta\Pi$ and $\delta\Pi^-$ in more detail in Subsections \ref{ss:modelling} and \ref{ss:robust},
we return to the informal discussion of the solution manifold of Subsection \ref{ss:para}.
Suppose we are given a curve $\mathbb{R}\ni t\mapsto\phi(t)$ on the space of
all $\phi$'s satisfying (\ref{ao40bis}) at fixed $\lambda$,
that passes through the ``point'' $\phi$
at time $t=0$. Then $\dot\phi:=\frac{d}{dt}_{|t=0}\phi(t)$ is a generic ``tangent vector''.
In view of (\ref{ao40bis}), it is characterized by
\begin{align}\label{ck88bis}
L\dot\phi-(3\lambda\phi^2+h^{(\rho)})\dot\phi=0 \quad\mbox{ mod analytic functions}.
\end{align}
Suppose that $\phi$ is parameterized by $(\lambda,p)$ via (\ref{ao41}). Then we
informally claim that there exist $\{\dot\pi^{(\bf n)}\}\subset\mathbb{R}$ such that
\begin{align}\label{ck86}
\dot\phi=\sum_{\beta}\big(\sum_{\bf n}\dot\pi^{({\bf n})}
(\partial_{\mathsf{z}_{\bf n}}\mathsf{z}^{\beta})[\lambda,p]\big)\Pi_\beta.
\end{align}
Hence (\ref{ck86}) parameterizes the tangent space of the solution manifold in 
the configuration $\phi$, in terms of $\{\dot\pi^{(\bf n)}\}\subset\mathbb{R}$. 

\medskip

\begin{sloppypar} 
Here $\partial_{\mathsf{z}_{\bf n}}$ informally denotes the partial derivative 
with respect to $\mathsf{z}_{\bf n}[p]$ $=\frac{1}{{\bf n}!}\partial^{\bf n}p(0)$;
it can be rigorously defined 
as an endomorphism on $\mathbb{R}[[\mathsf{z}_3,\mathsf{z}_{\bf n}]]$ via its matrix entries
\begin{align*}
(\partial_{\mathsf{z}_{\bf n}})_\beta^\gamma=\left\{\begin{array}{cc}
\gamma({\bf n})&\mbox{for}\;\gamma=\beta+\delta_{\bf n}\\
0&\mbox{else}
\end{array}\right\},
\end{align*}
\end{sloppypar}
which just encodes the desired action on monomials 
$\partial_{\mathsf{z}_{\bf n}}\mathsf{z}^{\gamma}$ 
$=\gamma({\bf n})\mathsf{z}^{\gamma-\delta_{\bf n}}$,
and automatically satisfies the finiteness properties (\ref{ck88}). In fact, it is a
derivation on the algebra $\mathbb{R}[[\mathsf{z}_3,\mathsf{z}_{\bf n}]]$,
by which the algebraists understand a linear endomorphism that satisfies Leibniz' rule
\begin{align}\label{ck92}
\partial_{\mathsf{z}_{\bf n}}\pi\pi'=(\partial_{\mathsf{z}_{\bf n}}\pi)\pi'
+\pi(\partial_{\mathsf{z}_{\bf n}}\pi')\quad\mbox{and}\quad
\partial_{\mathsf{z}_{\bf n}}\mathsf{1}=0
\end{align}
for all $\pi,\pi'\in\mathbb{R}[[\mathsf{z}_3,\mathsf{z}_{\bf n}]]$.
In view of the finiteness property (\ref{ck88}), such a derivation is characterized by
imposing its value on the coordinates; here $\partial_{\mathsf{z}_{\bf n}}\mathsf{z}_{3}=0$
and $\partial_{\mathsf{z}_{\bf n}}\mathsf{z}_{\bf m}=\delta_{\bf n}^{\bf m}$.

\medskip

The argument for (\ref{ck86}) is almost tautological: By definition of $\dot\phi$,
there exists a curve $t\mapsto\phi(t)$ with $\phi(t=0)=\phi$ and $\frac{d}{dt}_{|t=0}\phi(t)$
$=\dot\phi$; in view of (\ref{ao41}) it lifts to a curve $t\mapsto p(t)$ in parameter
space with $p(t=0)=p$ and $\phi(t)$ $=\sum_{\beta}\mathsf{z}^\beta[\lambda,p(t)]\Pi_\beta$.
Applying $\frac{d}{dt}_{|t=0}$ to this identity yields (\ref{ck86}) by the chain rule,
where 
\begin{equation*}
\dot\pi^{({\bf n})} := \mathsf{z}_{\bf n}\big[\frac{d}{dt}_{|t=0}p(t)\big] 
\end{equation*}
are the inner derivatives.

\medskip

We now algebrize (\ref{ck86}) by considering a $\dot\Pi\in X[[\mathsf{z}_3,\mathsf{z}_{\bf n}]]$
with
\begin{align}\label{ck87}
L\dot\Pi-(3\mathsf{z}_3\Pi^2+c)\dot\Pi=0 \quad\mbox{ mod analytic functions,}
\end{align}
and informally claim that this implies the representation
\begin{align}\label{ck90}
\dot\Pi=\sum_{\bf n}\dot\pi^{({\bf n})}\partial_{\mathsf{z}_{\bf n}}\Pi\quad\mbox{for some}\;
\{\dot\pi^{({\bf n})}\}\subset\mathbb{R}[[\mathsf{z}_3,\mathsf{z}_{\bf n}]].
\end{align}
Indeed, for arbitrary parameter $(\lambda,p)$ consider $\phi$ and $h^{(\rho)}$ given
by (\ref{ao41}), and $\dot\phi$ informally defined through
\begin{align}\label{ck89}
\dot\phi :=\sum_\beta\mathsf{z}^\beta[\lambda,p]\dot\Pi_\beta.
\end{align}
Then as (\ref{ao40bis})
did informally translate into (\ref{ao43bis}) \& (\ref{ao43}), so does (\ref{ck87}) 
translate back into (\ref{ck88bis}). Hence we may apply (\ref{ck86}); 
now the $\dot\pi^{({\bf n})}$'s implicitly depend on $(\lambda,p)$ 
and thus can be (informally) interpreted as 
elements of $\mathbb{R}[[\mathsf{z}_3,\mathsf{z}_{\bf n}]]$. Equating this representation
with (\ref{ck89}), and using that $(\lambda,p)$ was arbitrary, we obtain
\begin{align*}
\sum_\beta\sum_{\bf n}\dot\pi^{({\bf n})}(\partial_{\mathsf{z}_{\bf n}}\mathsf{z}^\beta)\Pi_\beta
=\sum_{\beta}\mathsf{z}^\beta\dot\Pi_\beta
\end{align*}
as an identity in $X[[\mathsf{z}_3,\mathsf{z}_{\bf n}]]$,
which amounts to (\ref{ck90}). 

\medskip

Finally, we claim that (\ref{ck90}) holds for an arbitrary base-point $x$, which for
variety we logically reverse and formulate as a rigorous version:
Provided the sequence $\{\dot\pi^{({\bf n})}_x\}$ 
$\subset\mathbb{R}[[\mathsf{z}_3,\mathsf{z}_{\bf n}]]$ is finite, then
\begin{align}\label{ck91}
\lefteqn{\dot\Pi=\sum_{\bf n}\dot\pi^{({\bf n})}_x\Gamma_x^*\partial_{\mathsf{z}_{\bf n}}\Pi_x}\nonumber\\
&\implies\quad
L\dot\Pi=(3\mathsf{z}_3\Pi^2+c)\dot\Pi \quad\mbox{ mod analytic functions}.
\end{align}
\begin{proof}[Proof of \eqref{ck91}]
By linearity, it is enough to consider $\dot\Pi$ 
$=\Gamma_x^*\partial_{\mathsf{z}_{\bf n}}\Pi_x$. Since $L$ commutes with 
$\Gamma_x^*\partial_{\mathsf{z}_{\bf n}}$, we obtain
$L\dot\Pi$ $=\Gamma_x^*\partial_{\mathsf{z}_{\bf n}}\Pi_x^{-}$mod analytic functions from (\ref{ck93}).
Applying the derivation $\partial_{\mathsf{z}_{\bf n}}$ to (\ref{ck84}) we obtain 
using Leibniz' rule (\ref{ck92}) that $\partial_{\mathsf{z}_{\bf n}}\Pi_x^{-}$
$=(3\mathsf{z}_3\Pi_x^2+c)\partial_{\mathsf{z}_{\bf n}}\Pi_x$.
We now apply $\Gamma_x^*$ to this identity;
by its multiplicativity (\ref{temp05}) followed by (\ref{temp06}), (\ref{temp02}), and 
(\ref{ck94}) this yields the desired r.~h.~s.~$(3\mathsf{z}_3\Pi^2+c)$ 
$\Gamma_x^*\partial_{\mathsf{z}_{\bf n}}\Pi_x$.
\end{proof}

\begin{figure}[t]
\includegraphics[width=0.6\textwidth]{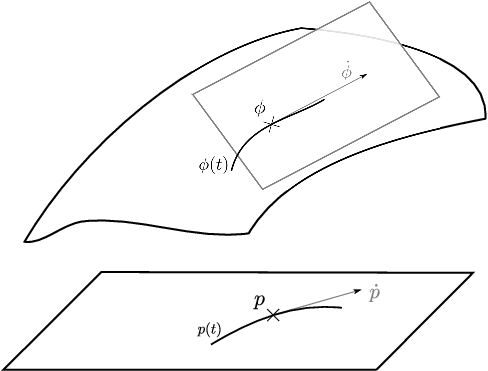}
\centering
\caption{
Heuristic visualization of the tangent space to the solution manifold.}
\label{fig:tangent}
\end{figure}


\subsection{Modelling the Malliavin derivative\texorpdfstring{ $\delta\Pi$ by ${\rm d}\Gamma_x^*$}{}}\label{ss:modelling}

We note that (\ref{ck69}) and (\ref{ck70bis}) combine to
\begin{align*}
L\delta\Pi=(3\mathsf{z}_3\Pi^2+c)\delta\Pi+\delta\xi_\rho \mathsf{1} \quad\mbox{ mod analytic functions};
\end{align*}
since $\delta\xi$ is not analytic, this is not exactly of the form \eqref{ck87}. 
Hence we cannot hope that the left statement of \eqref{ck91} holds for $\dot\Pi=\delta\Pi$. However, since
$\delta\xi\in\dot H^s$ has some regularity, we expect that it
holds approximately. More precisely, since $L$ is of second order and in view of (\ref{ck53}), 
we expect that the left statement of \eqref{ck91} holds up to an error of order $2+s$, i.~e.~there exist random $\{{\rm d}\pi_x^{({\bf n})}\}_{|{\bf n}|<2+s}$ 
$\subset\mathbb{R}[[\mathsf{z}_3,\mathsf{z}_{\bf n}]]$ such that
\begin{align}\label{ap29}
\delta\Pi-\sum_{|{\bf n}|<2+s}{\rm d}\pi_x^{({\bf n})}\Gamma_x^*
\partial_{\mathsf{z}_{\bf n}}\Pi_x=O(|\cdot-x|^{2+s}).
\end{align}
Note that by definition (\ref{ap40}), $2+s$ $=\alpha+\frac{D}{2}$,
so that in view of (\ref{ck28}) and (\ref{ap43})
we always have that $2+s>0$.
Property (\ref{ap29}) motivates to introduce the random endomorphism of
$\mathbb{R}[[\mathsf{z}_3,\mathsf{z}_{\bf n}]]$ via
\begin{align}\label{ck95}
{\rm d}\Gamma_x^*:=\sum_{|{\bf n}|<2+s}{\rm d}\pi_x^{({\bf n})}
\Gamma_x^*\partial_{\mathsf{z}_{\bf n}}.
\end{align}
We observe in passing that ${\rm d}\Gamma_x^*$ has the finiteness property (\ref{ck88}),
as a sum of products of operators that have this property; the latter
was already established for $\Gamma_x^*$ and $\partial_{\mathsf{z}_{\bf n}}$, and
is obvious for the operator $M$ of multiplication with an
$\pi\in\mathbb{R}[[\mathsf{z}_3,\mathsf{z}_{\bf n}]]$, 
which has the coordinate representation $M^\gamma_\beta$ $=\pi_{\beta-\gamma}$ 
with the understanding that this matrix element vanishes unless $\gamma\le\beta$
(coordinate-wise).

\medskip

We shall indeed establish a Schwartz-distributional version of (\ref{ap29}), 
see (\ref{ap61}) below, namely
\begin{align}\label{ck96}
(\delta\Pi-{\rm d}\Gamma_x^*\Pi_x)_r(x)=O(r^{2+s})\quad\mbox{as}\;r\downarrow0\;
\mbox{for all}\;x\in\mathbb{R}^{1+d}.
\end{align}
In this sense, $\delta\Pi_\beta$ is described (``modelled'' in the jargon of regularity structures) in terms
of $\{\Pi_{x\gamma}\}$ to order $2+s$; the coefficients are given by 
$\{({\rm d}\Gamma_x^*)_\beta^\gamma\}$ (they combine to a
``modelled distribution'' $\{({\rm d}\Gamma_x^*)_\beta^\gamma\}_{\gamma, x}$).
The statement (\ref{ck96}) is a multi-dimensional version of Gubinelli's controlled
rough-path condition.
We note that by the (qualitative) smoothness of $\delta\Pi$ and
$\Pi_x$ for $\rho>0$, (\ref{ck96}) implies
\begin{align}\label{ck96bis}
\partial^{\bf n}(\delta\Pi-{\rm d}\Gamma_x^*\Pi_x)(x)=0\quad\mbox{for}\;|{\bf n}|<2+s,
\end{align}
and shall argue in the next Subsection \ref{ss:struct} that this determines ${\rm d}\Gamma_x^*$.


\subsection{Uniqueness of \texorpdfstring{${\rm d}\Gamma_x^*$}{dGamma}}\label{ss:struct}

By definition (\ref{ck95}), we see that (\ref{temp05}) and (\ref{ck92}) translate into
\begin{align}
{\rm d}\Gamma^*_x \pi\pi'=({\rm d}\Gamma^*_x \pi)(\Gamma^*_x\pi')
+(\Gamma^*_x\pi)({\rm d}\Gamma^*_x\pi')
\quad\mbox{and}\quad{\rm d}\Gamma^*_x\mathsf{1}=0,\label{temp10}\\
{\rm d}\Gamma^*_x\mathsf{z}_3=0.\label{temp09}
\end{align}
For later purpose, we note that (\ref{temp09}) yields the counterpart of (\ref{ck94})
\begin{align}\label{ck94bis}
{\rm d}\Gamma_x^*c=0.
\end{align}
The three above properties motivate our notation of ${\rm d}$: ${\rm d}\Gamma_x^*$
like the Malliavin derivative $\delta\Gamma_x^*$ (which will not play a role
in these notes\footnote{Indeed, on the one hand $\delta\Gamma^*$ is too impoverished to model $\delta\Pi$ to order $2+s$ in the sense of \eqref{ck96}, which is the reason we introduce $\mathrm{d}\Gamma^*$;
on the other hand it is also not used to estimate $\Gamma^*$ via the spectral gap inequality. 
However, some estimates such as \eqref{ap75} below could be improved by using $\delta\Gamma^*$ as was done in \cite{LOTT21}.}) can be considered as a tangent vector to the
group of automorphisms of the algebra $\mathbb{R}[[\mathsf{z}_3,\mathsf{z}_{\bf n}]]$
in the group element $\Gamma_x^*$.

\medskip

However, ${\rm d}\Gamma_x^*$ does not have the population properties of a tangent vector to the 
structure group in $\Gamma_x^*$ (like $\delta\Gamma_x^*$ is):
We immediately obtain from (\ref{ck95}) the population condition
\begin{align}\label{cw59}
{\rm d}\Gamma^*_x\mathsf{z}_{\bf n}=0\quad\mbox{unless}\;|{\bf n}|<2+s.
\end{align}
As a consequence of \eqref{temp10} \& \eqref{temp09} in conjunction with
(\ref{cw59}), and due to its finiteness property, ${\rm d}\Gamma^*_x$ is 
determined by its values on the finitely many coordinates $\{\mathsf{z}_{\bf n}\}_{|{\bf n}|<2+s}$.

\medskip

We now use these algebraic properties to argue, 
mimicking the uniqueness argument for $\Gamma^*_x$ of Subsection \ref{ss:oneGamma},
that via (\ref{ck96bis}),
\begin{align}\label{ap00}
{\rm d}\Gamma_x^*\quad\mbox{is determined by}\;\delta\Pi,\;\mbox{next to}\;\Pi_x,\Gamma_x^*.
\end{align}
\begin{proof}[Proof of \eqref{ap00}]
Indeed, by (\ref{ck96bis}),
\begin{align}\label{ap02}
{\rm d}\Gamma_x^*\partial^{\bf n}\Pi_x(x)
\quad\mbox{is determined for $|{\bf n}|<2+s$}.
\end{align}
We now argue by induction in $k=|{\bf n}|<2+s$ that ${\rm d}\Gamma_x^*\mathsf{z}_{\bf n}$
is determined, which by the above
determines ${\rm d}\Gamma_x^*$.
The base case $k=0$, which just contains ${\bf n}={\bf 0}$, follows from
(\ref{ap02}), appealing to (\ref{ap01}) and (\ref{temp09}).
For the induction step $k-1\leadsto k$ we give ourselves an ${\bf n}$
with $|{\bf n}|=k$. By induction hypothesis and (\ref{temp10}) \& (\ref{temp09}),
we already identified ${\rm d}\Gamma_x^*$ on 
$\mathbb{R}[[\mathsf{z}_3,\{\mathsf{z}_{\bf m}\}_{|{\bf m}|< k}]]$.
Hence via (\ref{ap01}) we learn from this and (\ref{ap02}) that
${\rm d}\Gamma_x^*\mathsf{z}_{\bf n}$ is determined.
\end{proof}


\subsection{A robust relation \texorpdfstring{$\delta\Pi\mapsto\delta\Pi^{-}$}{between deltaPi and deltaPi-}}\label{ss:robust}

We now claim that (\ref{ck96bis}) implies
\begin{align}\label{ck97}
(\delta\Pi^{-}-\delta\xi_\rho\mathsf{1}-{\rm d}\Gamma_x^*\Pi_x^{-})(x)=0.
\end{align}
The combination of (\ref{ck96bis}) and (\ref{ck97}) provides the desired robust 
map $\delta\Pi\mapsto\delta\Pi^{-}$ that substitutes the non-robust $\Pi\mapsto\Pi^{-}$
given by (\ref{ao43}); in the sense that it bypasses the divergent $c$: 
In view of (\ref{ap00}), ${\rm d}\Gamma_x^*$ 
is uniquely determined by (\ref{ck96}) in terms of $\delta\Pi$ (at given $\Pi_x,\Gamma_x^*$),
so that (\ref{ck97}) determines $\delta\Pi^{-}$ (at given $\Pi_x,\Pi_x^{-},\Gamma_x^*$
and of course $\delta\xi$). 
Hence the mission
of Subsection \ref{ss:usage} is accomplished.

\medskip

\begin{proof}[Proof of \eqref{ck97}]
We start the argument for (\ref{ck97}) by noting that (\ref{ck96bis}) implies in particular
\begin{align*}
(\delta\Pi-{\rm d}\Gamma_x^*\Pi_x)(x)=0.
\end{align*}
In view of (\ref{ck84}) and (\ref{ck69}), we may pass from this to
(\ref{ck97}) based on the identity
\begin{align*}
{\rm d}\Gamma_x^*(\mathsf{z}_3\Pi_x^3+c\Pi_x+\xi_\rho\mathsf{1})
=(3\mathsf{z}_3\Pi^2+c){\rm d}\Gamma_x^*\Pi_x,
\end{align*}
which itself follows from the rules (\ref{temp10}) \& (\ref{temp09}) \& (\ref{ck94bis}),
the rules (\ref{temp05}) \& (\ref{temp06}) they are based on,
and (\ref{temp02}).
\end{proof}


\subsection{Population of \texorpdfstring{${\rm d}\Gamma_x^*$}{dGamma}}\label{ss:order1}

We claim that in analogy to (\ref{cw46}) \& (\ref{cw46ter}) we have the population property
\begin{align}\label{ap17}
\mbox{for all populated $\gamma$}:\quad
({\rm d}\Gamma_x^*)_\beta^\gamma=0\quad\mbox{unless}\;[\beta]\ge 0.
\end{align}
For $\beta=0$, we claim the more precise information
\begin{align}\label{ap63}
({\rm d}\Gamma_x^*)_0^\gamma=0\quad\mbox{unless}\;\gamma=\delta_{\bf n}\;
\mbox{with}\;|{\bf n}|<2+s.
\end{align}

\medskip

\begin{proof}[Proof of \eqref{ap17} \textnormal{\&} \eqref{ap63}]
We follow Subsection \ref{ss:triang} and establish (\ref{ap17}) by induction in $|\beta|$.
We follow the argument for (\ref{cw46}) \& (\ref{cw46ter}), just indicating the changes.
By the last item in (\ref{temp10}), the case of $\gamma=0$ is automatically satisfied.
As in Subsection \ref{ss:triang}, we start with the cases (\ref{cw71}) \& (\ref{cw71ter}).
Since the Ansatz (\ref{cw71bis}) translates
into $\mathsf{z}^\gamma$ $=\mathsf{z}_3^k\prod_{j=1}^l\mathsf{z}_{{\bf n}_j}$,
we learn from (\ref{temp05}) \& (\ref{temp06}) and from (\ref{temp10}) \& (\ref{temp09})
that ${\rm d}\Gamma_x^*\mathsf{z}^\gamma$ is a (finite) linear combination of
terms of the form $\mathsf{z}_3^k$ $({\rm d}\Gamma_x^*\mathsf{z}_{{\bf n}_1})$
$\prod_{j=2}^l\Gamma_x^*\mathsf{z}_{{\bf n}_j}$. Hence in view of (\ref{ck74})
we obtain the analogue of (\ref{ap09}):
\begin{align}\label{ap19}
&({\rm d}\Gamma_x^*)^\gamma_\beta=\mbox{linear combination of}\quad
({\rm d}\Gamma_x^*)^{\delta_{{\bf n}_1}}_{\beta_1}
(\Gamma_x^*)^{\delta_{{\bf n}_2}}_{\beta_2}\cdots(\Gamma_x^*)^{\delta_{{\bf n}_l}}_{\beta_l}\nonumber\\
&\mbox{where}\;k\delta_3+\beta_1+\cdots+\beta_l=\beta.
\end{align}
As in Subsection \ref{ss:triang} we argue that $|\beta_j|<|\beta|$ for $j=1,\dots,l$.
On the one hand, for $j=1$ we learn from the induction hypothesis that
the r.~h.~s.~term in (\ref{ap19}) vanishes unless $[\beta_1]\ge 0$.
For $j=2,\dots,l$ we infer from
\eqref{cw47} that it vanishes unless $\beta_j$ is
populated and thus in particular $[\beta_j]\ge -1$.
Hence by additivity of $[\cdot]$ we learn that $[\beta]{\ge2k-l+1}$.
On the other hand, by (\ref{cw71bis}) we have $[\gamma]$ $=2k-l$,
and since $\gamma$ is assumed to be populated so that $[\gamma]+1\ge 0$
we have $2k-l+1\ge 0$.
In combination, we obtain the desired $[\beta]\ge 0$.

\medskip

We conclude the induction step for (\ref{ap17}) by treating the case of (\ref{ap07}),
i.~e.~ $\gamma=\delta_{\bf n}$.
According to (\ref{cw59}) we may restrict to the case ${|{\bf n}|<2+s}$.
We rewrite (\ref{ck96bis}) component-wise
$\partial^{\bf n}\delta\Pi_\beta(x)=\sum_{\gamma}({\rm d}\Gamma_x^*)_\beta^\gamma
\partial^{\bf n}\Pi_{x\gamma}(x)$.
Splitting the sum into $\gamma$ that are pp, on which we use (\ref{ck38}),
and those that are not, we obtain the representation
\begin{align*}
({\rm d}\Gamma_x^*)_\beta^{\delta_{\bf n}}=
\tfrac{1}{{\bf n}!}\partial^{\bf n}\delta\Pi_\beta(x)
-\sum_{\gamma\textnormal{ not pp}}({\rm d}\Gamma_x^*)_\beta^\gamma\tfrac{1}{{\bf n}!}\partial^{\bf n}\Pi_{x\gamma}(x).
\end{align*}
In view of (\ref{ao08}), the second factor $\partial^{\bf n}\Pi_{x\gamma}(x)$
vanishes unless $\gamma$ is populated; if $\gamma$ is populated,
in view of what we showed in the previous paragraph, the first factor
$({\rm d}\Gamma_x^*)_\beta^\gamma$ vanishes unless $[\beta]\ge 0$.
In view of again (\ref{ao08}) and its Malliavin derivative,
the first term $\partial^{\bf n}\delta\Pi_\beta(x)$ vanishes unless $[\beta]\ge 0$.
This concludes the induction step and thus the induction argument for (\ref{ap17}).

\medskip

We finally tackle (\ref{ap63}).
Writing $\gamma$ as in (\ref{cw71bis}), we obtain (\ref{ap19})
with $\beta_j=0$ for all $j=1,\dots,l$; it vanishes unless $k=0$.
Since by (\ref{temp04}) and (\ref{ck28}) we have
$|\delta_{{\bf n}_j}|$ $=|{\bf n}_j|$ $\ge 0$ $>|0|$, it follows from
(\ref{cw44}) that $(\Gamma_x^*)_0^{\delta_{{\bf n}_j}}$ $=0$,
so that (\ref{ap19}) vanishes unless $l=1$. Thus $\gamma$ is pp, i.~e.~$\gamma$
$=\delta_{\bf n}$; the remainder of (\ref{ap63}) thus follows from (\ref{cw59}).
\end{proof}


\subsection{Strict triangularity of \texorpdfstring{${\rm d}\Gamma_x^*$}{dGamma},
order\texorpdfstring{ $\prec$}{} of induction}\label{ss:order}

Theorem~\ref{th:main} is established inductively in $\beta$;
in view of (\ref{ao08}), (\ref{temp17}), (\ref{temp26}), (\ref{cw46}),
and (\ref{cw46bis}) it is sufficient to treat $\beta$ with $[\beta]\ge 0$.
The inductive proof relies on triangular properties, in particular those of ${\rm d}\Gamma_x^*$.  
While ${\rm d}\Gamma_x^*$ has the same algebraic properties as $\delta\Gamma_x^*$,
namely (\ref{temp10}) \& (\ref{temp09}),
its population pattern is quite different, cf.~(\ref{cw59}). 
As a consequence, ${\rm d}\Gamma_x^*$ is
not strictly triangular w.~r.~t.~$|\cdot|$; in fact we just have
\begin{align}\label{ap70}
({\rm d}\Gamma_x^*)_\beta^\gamma=0\quad\mbox{unless}\;|\gamma|<|\beta|+\frac{D}{2},
\end{align}
as we shall argue now by induction in $|\beta|$: 

\begin{proof}[Proof of \eqref{ap70}]
The base case follows from (\ref{ap63}),
where we use (\ref{ao12}), which implies $|0|$ $=s-\frac{D}{2}+2$, and (\ref{temp04}).
For the induction step we distinguish two cases.
If $\gamma$ is of the form \eqref{ap07}, 
i.~e.~$\gamma=\delta_{\bf n}$, 
then we conclude by \eqref{cw59} which yields $({\rm d}\Gamma^*_x)_\beta^\gamma=0$ unless
\begin{align*}
|\gamma|=|\delta_{\bf n}|
\stackrel{\eqref{temp04}}{=} |{\bf n}|
\stackrel{\eqref{cw59}}{<}2+s 
\stackrel{\eqref{ap40}}{=}\alpha+\frac{D}{2}
\stackrel{\eqref{temp08}}{\leq} |\beta|+\frac{D}{2} .
\end{align*}
If $\gamma$ is of the form \eqref{cw71} or \eqref{cw71ter}, 
we appeal to (\ref{ap19}); by (\ref{temp08}) we have
$|\beta|-|0|$ $=(|k\delta_3|-|0|)$ $+\sum_{j=1}^l(|\beta_j|-|0|)$,
by the induction hypothesis and (\ref{cw44}), the term only contributes if
$|\beta|-|0|$ $>(|k\delta_3|-|0|)$ $+\sum_{j=1}^l(|\delta_{{\bf n}_j}|-|0|)$ $-\frac{D}{2}$,
which by (\ref{cw71bis}) and once more by (\ref{temp08}) amounts to the desired
$|\beta|-|0|$ $>(|\gamma|-|0|)$ $-\frac{D}{2}$.
\end{proof}

\medskip

For the proof of Theorem \ref{th:main}, we need to find an order
$\prec$ on the multi-indices $\beta$ with $[\beta]\ge0$ w.~r.~t.~which both $\Gamma_x^*$
and ${\rm d}\Gamma_x^*$ are strictly triangular. Here are the three properties 
of the ordinal $|\beta|_{\prec}$ we need: On the subset of multi-indices the induction takes place,
the basic feature (\ref{temp08}) is required:
\begin{align}\label{ap12}
&\mbox{for $\beta$ with $[\beta]\ge 0$}:\nonumber\\
&|\cdot|_{\prec}-|0|_{\prec}\;\mbox{is additive},\quad\ge 0,\quad\mbox{and $=0$ only for 
$\beta=0$}.
\end{align}
The ordinal needs to be comparable to $|\cdot|$ in the sense of
\begin{align}\label{ap13}
\mbox{for populated $\beta$}:\quad
|\beta|_{\prec}\ge|\beta|\quad\mbox{and}\quad|\delta_{\bf n}|_{\prec}=|\delta_{\bf n}|
\end{align}
(which also ensures that coercivity, cf.~(\ref{temp15}), is preserved), 
and the ordinal needs to dominate the truncation order of ${\rm d}\Gamma_x^*$, see (\ref{cw59}):
\begin{align}\label{ap14}
\mbox{for $\beta$ with $[\beta]\ge 0$}:\quad|\beta|_{\prec}\ge2+s.
\end{align}
Since $[\beta]+1\ge 0$ for populated $\beta$ and $[\beta]+1=0$ for
purely polynomial $\beta$, cf.~\eqref{temp24}, it follows from (\ref{ck28}) and (\ref{temp08})
in conjunction with $2+s$ $=\alpha+\frac{D}{2}$, cf.~(\ref{ap40}),
that these postulates are satisfied by
\begin{align}\label{ap64}
|\beta|_\prec:=|\beta|+\frac{D}{2}([\beta]+1).
\end{align}
We note that $|\beta|_{\prec}$ differs from $|\beta|$, in its representation 
(\ref{ao12}), by replacing the pre-factor $s-\frac{D}{2}$ of the noise homogeneity $[\beta]+1$
by $s$.

\medskip

We claim the strict triangularities 
\begin{align}
\mbox{for all $\gamma$}:&\quad
(\Gamma_x^*-{\rm id})_\beta^\gamma=0\quad\mbox{unless}\;|\gamma|_\prec<|\beta|_\prec,\label{ap15}\\
\mbox{for all $\gamma$}:&\quad
({\rm d}\Gamma_x^*)_\beta^\gamma=0\quad\mbox{unless}\;|\gamma|_\prec<|\beta|_\prec.\label{ap16}
\end{align}

\medskip

\begin{proof}[Proof of \eqref{ap15} \textnormal{\&} \eqref{ap16}]
The strict triangularity (\ref{ap15}) is established closely following the argument
for (\ref{cw44}) in Subsection \ref{ss:triang2}, which is based on induction in $|\beta|$.
We start with the base case of $\beta=0$; by (\ref{temp08}),
(\ref{cw44}) assumes the form that $(\Gamma_x^*-{\rm id})_0^\gamma=0$ for all $\gamma$,
which trivially implies (\ref{ap15}). In the induction step, we  
distinguish whether $\gamma$ is purely polynomial or not, following (\ref{ap69}).
In case of $\gamma\not=$ pp, 
we may reproduce the argument for \eqref{cw44} with $|\cdot|$ replaced by $|\cdot|_\prec$ 
since the inspection of it reveals that it only relies on \eqref{ap12}.
We now turn to the case of $\gamma=\delta_{\bf n}$; by (\ref{cw46}) we may
assume that $\beta$ is populated. Hence our requirement (\ref{ap13}) ensures
that we may pass from (\ref{cw44}) to (\ref{ap15}).
This concludes the proof of (\ref{ap15}).

\medskip

We now turn to (\ref{ap16}) which we establish by induction in the
plain length $\gamma(3)+\sum_{\bf n}\gamma({\bf n})$.
For $\gamma=0$ and $\gamma=\delta_3$
we have $({\rm d}\Gamma^*_x)_\beta^\gamma=0$ by the second item in (\ref{temp10})
and by (\ref{temp09}), respectively.
For $\gamma=\delta_{\bf n}$ we note that by (\ref{cw59}) we have on the one hand
$({\rm d}\Gamma_x^*\mathsf{z}_{\bf n})_\beta=0$ unless $|{\bf n}|<2+s$,
which by (\ref{temp04}), the second item in (\ref{ck74}), and (\ref{ap13}), translates into
$({\rm d}\Gamma_x^*)_\beta^{\delta_{\bf n}}=0$ unless $|\delta_{\bf n}|_{\prec}<2+s$.
On the other hand, since the purely polynomial $\gamma=\delta_{\bf n}$ is in
particular populated, we may by (\ref{ap17}) restrict to $[\beta]\ge 0$,
so that by postulate (\ref{ap14}) we have $|\beta|_{\prec}\ge2+s$. 
Both statements combine into the desired statement
$({\rm d}\Gamma_x^*)_\beta^{\delta_{\bf n}}=0$ unless $|\delta_{\bf n}|_\prec<|\beta|_{\prec}$.

\medskip

Turning to the induction step we are given a $\gamma$ of plain length $\ge 2$,
which we write as $\gamma=\gamma_1+\gamma_2$
with $\gamma_1$ and $\gamma_2$ of plain length strictly less than that of $\gamma$.
As for (\ref{ap19}), this implies
\begin{align}\label{cw57bis}
({\rm d}\Gamma^*_x)_\beta^\gamma=\sum_{\beta_1+\beta_2=\beta} \big(
( {\rm d}\Gamma^*_x)_{\beta_1}^{\gamma_1} (\Gamma^*_x)_{\beta_2}^{\gamma_2}
+(\Gamma^*_x)_{\beta_1}^{\gamma_1}({\rm d}\Gamma^*_x)_{\beta_2}^{\gamma_2} \big).
\end{align}
According to additivity of $|\cdot|_{\prec}-|0|_{\prec}$, 
$|\gamma|_\prec\ge|\beta|_\prec$ would imply
$|\gamma_1|_\prec$ $+|\gamma_2|_\prec$ $\ge|\beta_1|_\prec$ $+|\beta_2|_\prec$,
and thus 
$|\gamma_1|_\prec\ge|\beta_1|_\prec$ or $|\gamma_2|_\prec\ge|\beta_2|_\prec$; 
w.~l.~o.~g.~we may restrict to the former.
In this case the first term in (\ref{cw57bis}) vanishes 
because its first factor vanishes by induction hypothesis.
By (\ref{ap15}), the second term vanishes unless $\gamma_1=\beta_1$ which implies
$\gamma_2=\beta_2$, so that also this term vanishes by induction hypothesis.
\end{proof}


\subsection{Usage of the spectral gap inequality}\label{ss:spectral_gap}

It is well-known, see \cite[Theorem~5.5.1]{Bog98}, that for a (centered) Gaussian ensemble that satisfies (\ref{ck81bis}), 
the variance of a random variable $F$
is dominated by the expectation of the carr\'e-du-champs of its Malliavin derivative
\begin{align}\label{ap19bis}
\Big\|\frac{\partial F}{\partial\xi}\Big\|_{\dot H^{-s}}^2
=\Big(\sup_{\delta\xi}\frac{\delta F}{\|\delta\xi\|_{\dot H^s}}\Big)^2
\end{align}
in the sense of
\begin{align}\label{ap18}
\mathbb{E}(F-\mathbb{E}F)^2\le\mathbb{E}\Big\|\frac{\partial F}{\partial\xi}\Big\|_{\dot H^{-s}}^2
\end{align}
for any (reasonable) random variable $F$ -- we continue to be informal. Note that (\ref{ck81bis})
is a special case of (\ref{ap18}) for the simple cylinder function(al)
$F[\xi]=(\xi,\zeta)$. The inequality (\ref{ap18}) 
can be seen as a Poincar\'e inequality in probability; it bounds the spectral gap
of the generator of the stochastic process that is defined on the basis of
the Dirichlet form on the r.~h.~s~of (\ref{ap18}) and has the ensemble at
hand as a stationary measure. Our assumption of Gaussianity in conjunction
with (\ref{ck81bis}) can be replaced by directly assuming (\ref{ap18}) (and the closability
of the Malliavin derivative).

\medskip

An easy argument based on Leibniz' rule and H\"older's inequality (see e.~g.\ \cite[Proposition~5.1]{IORT23}) shows that (\ref{ap18})
can be upgraded to the following $\mathbb{E}^\frac{1}{p}|\cdot|^p$-version for any $2\le p<\infty$
\begin{align}\label{ap20bis}
\mathbb{E}^\frac{1}{p}|F-\mathbb{E}F|^p\lesssim
\mathbb{E}^\frac{1}{p}\Big\|\frac{\partial F}{\partial\xi}\Big\|_{\dot H^{-s}}^p,
\end{align}
where the implicit constant depends on $p$.
In view of (\ref{ap19bis}), (\ref{ap20bis}) can be reformulated as
\begin{align*}
\mathbb{E}^\frac{1}{p}|F-\mathbb{E}F|^p&\lesssim\sup_{\delta\xi\;\textnormal{random}}
\frac{\mathbb{E}\delta F}{\mathbb{E}^\frac{1}{p^*}\|\delta\xi\|_{\dot H^s}^{p^*}}\nonumber\\
&=\sup\{\,\mathbb{E}\delta F\,|\,
\delta\xi\;\mbox{random with}\;\mathbb{E}\|\delta\xi\|_{\dot H^s}^{p^*}\le1\,\},
\end{align*}
where $1<p^*\le 2$ is the dual exponent to $p$, i.~e.~$\frac{1}{p}+\frac{1}{p^*}=1$.
Hence 
our task at hand is to estimate $\delta F$ for random variables
of interest like $F=\Pi_{\beta r}^{-}(0)$. In fact, for the induction argument,
it will be important to monitor a norm instead of its expectation, namely
\begin{align}\label{ap20.4}
\mathbb{E}^\frac{1}{q}|\delta\Pi_{\beta r}^{-}(0)|^{q}\quad
\mbox{for all}\;1\le q<p^*\le 2 \le p.
\end{align}
That is, we will use the spectral gap assumption in the form
	\begin{align} 
		\mathbb{E}^\frac{1}{p}|\Pi_{\beta r}^-(0)|^p
		&\lesssim |\mathbb{E}\Pi_{\beta r}^-(0)| \nonumber\\
		& \quad +\sup\{\,\mathbb{E}^\frac{1}{q}|\delta\Pi_{\beta r}^{-}(0)|^{q}\,|\,
\delta\xi\;\mbox{random with}\;\mathbb{E}\|\delta\xi\|_{\dot H^s}^{p^*}\le1\,\}. \label{ap20ter}
	\end{align}
In what follows, all estimates will be proved simultaneously for all $p, q$ in the range \eqref{ap20.4} (with implicit constants depending on $p, p^*, q$), so that it will not be a problem to recursively appeal to the same estimates with exponents $q^{\prime} \in (q, p^*)$ or $p^{\prime} > p$, as happens e.~g.\ when applying H\"older's inequality in probability.
In practice: Within the induction we \emph{fix} $p, q$ and a random $\delta\xi\in \dot H^s$ with
\begin{align}\label{ap71}
\mathbb{E}\|\delta\xi\|_{\dot H^s}^{p^*}\le1,
\end{align}
and estimate Malliavin derivatives in the direction $\delta \xi$.

\medskip

In order to complete the argument \eqref{ap20ter}, also the expectation $| \mathbb{E}\Pi_{\beta r}^{-}(0)|$ has to be bounded.
This will follow from \eqref{ck52ter}, see Subsection~\ref{ss:expectation}.
This relies on the BPHZ-choice of renormalization, so that this requires the restriction to $|\beta|<2$.


\subsection{Besov-type norms and base case}\label{ss:Besov}

At the core of the proof is a quantification of (\ref{ck96}). From now onwards,
we restrict ourselves to $\beta$ with $[\beta]\ge 0$ and $|\beta|<2$ (the case\footnote{recall that the case $|\beta|=2$ is trivial by \eqref{temp17} and \eqref{temp04bis}} $|\beta|>2$ is treated differently, see Subsection~\ref{ss:logic}). We start by noting
\begin{align}\label{ap27}
\lefteqn{L(\delta\Pi-{\rm d}\Gamma_x^*\Pi_x)_{\beta} 
- (\delta\Pi^{-}-{\rm d}\Gamma_x^*\Pi_x^{-})_{\beta}
}\nonumber\\
&=0
\quad\mbox{ mod polynomials of degree $<|\beta|+\frac{D}{2}-2$}.
\end{align}
Indeed, this follows from (\ref{ck93}), (\ref{ck70bis}), and (\ref{ap70}).

\medskip

We momentarily consider the base case $\beta=0$; in view of (\ref{ao12}),
the polynomial in (\ref{ap27}) is of degree $<s$, and 
according to (\ref{ao46}) and (\ref{ap63}), (\ref{ap27}) thus collapses to 
\begin{align*}
L\big(\delta\Pi_{0}
-\sum_{|{\bf n}|<2+s}({\rm d}\Gamma_x^*)_{0}^{\delta_{\bf n}}(\cdot-x)^{\bf n}\big)
=\delta\xi_\rho
\mbox{ mod polynomials of degree}<s.
\end{align*}
This shows that the order of truncation of what now is a Taylor polynomial
of order $<2+s$ of $\delta\Pi_0$
is consistent with the order of regularity $s$ of $\delta\xi$ and the order $2$ of $L$.
It also shows that in its pointwise-in-$x$ form, statement (\ref{ck96}) is too strong;
in fact, because of the $L^2$-based nature of $\dot H^s$ we only have
\begin{align*}
\big(\int dx|\delta\Pi_{0}(x+y)
-\sum_{|{\bf n}|<2+s}({\rm d}\Gamma_x^*)_{0}^{\delta_{\bf n}}y^{\bf n}|^2
\big)^\frac{1}{2}
& \lesssim |y|^{2+s}\|\delta\Pi_0\|_{\dot{H}^{2+s}} \\
& \lesssim |y|^{2+s}\|\delta\xi\|_{\dot H^s} ,
\end{align*}
which is easily seen by Fourier transformation.
This amounts to an estimate of $\delta\Pi_0$ in the Besov space $\dot B^{2+s}_{2,\infty}$.
Because of Minkowski's inequality (recall $q\leq 2$) 
and the normalization (\ref{ap71}), it implies the ``annealed'' estimate
\begin{align}\label{ap72}
\big(\int dx\mathbb{E}^\frac{2}{q}|\delta\Pi_{0}(x+y)
-\sum_{|{\bf n}|<2+s}({\rm d}\Gamma_x^*)_{0}^{\delta_{\bf n}}y^{\bf n}|^q
\big)^\frac{1}{2}\lesssim |y|^{2+s}.
\end{align}

\medskip

What can we expect for $(\delta\Pi-{\rm d}\Gamma_x^*\Pi_x)_\beta$ with
$\beta\not=0$? First of all, as discussed in Subsection \ref{ss:usage}
in case of $\beta=\delta_3$, 
we cannot expect to estimate $\delta\Pi_\beta$ as a (square-integrable) function,
but just as a distribution. Hence in anticipation, we relax (\ref{ap72}) to
\begin{align}\label{ap73}
\big(\int dx\mathbb{E}^\frac{2}{q}|\big(\delta\Pi-{\rm d}\Gamma_x^*\Pi_x)_{0r}(x)|^q
\big)^\frac{1}{2}\lesssim r^{2+s}\stackrel{\eqref{ap40}}{=}r^{\alpha+\frac{D}{2}},
\end{align}
where the constant now implicitly depends on the control \eqref{ck52a}.
Moreover, in view of its structure (\ref{ap19}), ${\rm d}\Gamma_x^*$ will acquire (some of) the growth
(\ref{ck76}) of $\Gamma_x^*$ in $x$, while according to (\ref{temp30}),
the law of $\Pi_{x\gamma r}(x)$ is independent of $x$, so that one cannot expect
$\mathbb{E}^\frac{1}{q}|({\rm d}\Gamma_x^*\Pi_x)_{\beta r}(x)|^q$ 
to be square integrable (at infinity) in $x$.
Hence for $\beta\not=0$, we need to relax (\ref{ap73}) by restricting the $x$-integral
to a (parabolic) ball ${B_R:=\{|x|<R\}}$. So next to the mollification length scale $r$,
we acquire a second length scale, the localization scale $R$. 
This motivates the form of the l.~h.~s.~of
\begin{align}\label{ap61}
\big(\int_{B_R}dx\mathbb{E}^\frac{2}{q}&
|(\delta\Pi-{\rm d}\Gamma_x^*\Pi_x)_{\beta r}(x)|^{q}\big)^\frac{1}{2}
\lesssim r^{\alpha+\frac{D}{2}}(r+R)^{|\beta|-\alpha},
\end{align}
which is indeed what we shall establish -- and just have established for $\beta=0$
in view of (\ref{ck28}).

\medskip

Note that compared to (\ref{ao65}), 
passing from $\Pi$ to $\delta\Pi-{\rm d}\Gamma_x^*\Pi_x$ comes
with a (beneficial) factor of $r^\frac{D}{2}$, and $R$ plays the role of $|x|$.
There are two consistency checks for the two exponents in (\ref{ap61}):
1) For $r\ll R$, (\ref{ap61}) contains the expected $O(r^{\alpha+\frac{D}{2}=2+s})$-behavior
of (\ref{ck96}).
2) The dimension in terms of length of (\ref{ap61}) is
consistent with the one of (\ref{ao65}) since the $L^2(B_R)$-norm contributes
$\frac{D}{2}$ dimensions of length.

\medskip

While we have established (\ref{ap61}) in the base case of $\beta=0$, for
$\beta\not=0$, we shall derive it by ``integration'' of (\ref{ap27}) from 
\begin{align}\label{ap21bis}
\big(\int_{B_R}dx\mathbb{E}^\frac{2}{q}
|(\delta\Pi^{-}-{\rm d}\Gamma_x^*\Pi_x^{-})_{\beta r}(x)|^{q}\big)^\frac{1}{2}
\lesssim r^{\alpha-2+\frac{D}{2}}(r+R)^{|\beta|-\alpha}.
\end{align}
In fact we shall establish the stronger
\begin{align}\label{ap21ter}
\big(\int_{B_R}dx\mathbb{E}^\frac{2}{q}
|(\delta\Pi^{-}-{\rm d}\Gamma_x^*\Pi_x^{-})_{\beta r}(x)|^{q}\big)^\frac{1}{2}
\lesssim r^{3\alpha+\frac{D}{2}}(r+R)^{|\beta|-2-3\alpha},
\end{align}
which, as a quantitative version of \eqref{ck97} (recall $3\alpha+\frac{D}{2}>0$ by assumption \eqref{ap43}), is the estimate required for the characterization of the model in \cite{Tem}.
One can think of \eqref{ap21ter} as the corresponding Malliavin version of  \eqref{ao65bisimproved}.
The next four subsections are devoted to the induction step that establishes
(\ref{ap61}) \& \eqref{ap21ter}. It will be based on the four
consecutive steps of an {\sc algebraic argument}, {\sc reconstruction}, {\sc integration},
and a {\sc three-point argument}; in this section, we shall focus on the algebraic aspects
of all steps, while the analytic ingredients will be detailed in Section~\ref{s:proofs}. 
As it turns out, these arguments do not use the fact that
$\delta\Pi$ and $\delta\Pi^{-}$ are Malliavin derivatives in direction of $\delta\xi$,
but just rely on the relations (\ref{ck69}) \& (\ref{ck70bis}), and the control (\ref{ap71}).
Hence we actually provide an a priori estimate for the inverse of the linearization
of $\dot\phi\mapsto L\dot\phi-(3\lambda\phi^2+h^{({\rho})})\dot\phi$, on our term-by-term
level, and in an $L^2$-based norm. 
Hence our approach blends solution theory and stochastic estimates.

\medskip

The induction step laid out in the next four subsections 
is restricted to multi-indices $\beta$ with
$|\beta|<$ $\lceil s+2\rceil$, which is needed for integration. 
Hence in our induction w.~r.~t.~$|\cdot|_{\prec}$, we have to ensure that we
only use the induction hypothesis under this additional constraint.
The range $|\beta|<$ $\lceil s+2\rceil$ does cover the desired range of $|\beta|<2$ 
iff $1 < s+2=\alpha+\frac{D}{2}$, which in view of
(\ref{ap43}) and $\alpha<0$ is the case iff\footnote{if $1<\alpha+\frac{D}{2}$ then $D>2$ by $\alpha<0$; conversely, \eqref{ap43} implies $\alpha>-\frac{D}{6}$, hence $\alpha+\frac{D}{2}>\frac{D}{3}\geq1$ by $D\geq3$} $D\geq3$, which we have assumed.


\subsection{The algebraic argument for 
\texorpdfstring{${\rm d}\Gamma_{x+y}^*-{\rm d}\Gamma_x^*\Gamma_{x\,x+y}^*$}{increments of dGamma}}\label{ss:algIII}

Estimate \eqref{ap21ter} will be a consequence of what in regularity structures
is called reconstruction. 
It states that a distribution like $\delta\Pi_\beta$
can be ``reconstructed'' from the family of distributions like 
$\{({\rm d}\Gamma_x^*\Pi_x^{-})_\beta\}_x$ 
that act as germs near every space-time point $x$. 
For this to be canonically feasible, 
the distributions $({\rm d}\Gamma_x^*\Pi_x^{-})_\beta$ have to satisfy a continuity condition
w.~r.~t.~the base-point $x$, see (\ref{ap53}) below.
In the next Subsection~\ref{ss:reconstr}, (\ref{ap53}) is derived from graded continuity of
$\{({\rm d}\Gamma_x^*\Gamma_x^{-*})_\beta^\gamma\}_\gamma$ in $x$.
In line with the continuity condition on modelled distributions in regularity structures,
this graded continuity is formulated in terms of smallness of the increment 
${\rm d}\Gamma^*_{x+y}-{\rm d}\Gamma_x^*\Gamma_{x\;x+y}^*$ in terms of the shift $y$, where by (\ref{ap45}) 
$\Gamma_{x\;x+y}^*=\Gamma_x^{-*}\Gamma_{x+y}^*$.
Also this graded continuity is formulated in an $L^2(B_R)$ sense: 
\begin{align}\label{ap25}
&\mbox{for populated $\gamma$}:\quad\big(\int_{B_R}dx\mathbb{E}^\frac{2}{q}|({\rm d}\Gamma^*_{x+y}
-{\rm d}\Gamma_x^*\Gamma_{x\;x+y}^*)_\beta^\gamma|^{q}\big)^\frac{1}{2}\nonumber\\
&\lesssim\left\{\begin{array}{clc}
|y|^{\alpha-|\gamma|_p+\frac{D}{2}}\hspace{-2ex}\mbox{}&(|y|+R)^{|\beta|-|\gamma|+|\gamma|_p-\alpha}
&\mbox{if}\;\alpha-|\gamma|_p+\frac{D}{2}>0\\
&(|y|+R)^{|\beta|-|\gamma|+\frac{D}{2}}&\mbox{else}
\end{array}\right\},
\end{align}
where we have set for abbreviation
\begin{align}\label{ap48}
|\gamma|_p:=\sum_{\bf n}|{\bf n}|\gamma({\bf n}),
\end{align}
and with the understanding that the
\begin{align}\label{ap46}
\lefteqn{\mbox{l.~h.~s.~of (\ref{ap25}) vanishes}}\nonumber\\
&\mbox{unless}\quad
\sum_{\bf n}\gamma({\bf n})\not=0\quad\mbox{and}\quad|\beta|-|\gamma|+|\gamma|_p-\alpha\ge 0.
\end{align}
Such a ``pointed'' Besov-type norm \eqref{ap25} was introduced in \cite[Definition~3.9]{HS23} in a very similar context.
This continuity condition (\ref{ap25}) in $y$ could be strengthened to yield
a positive H\"older exponent in all cases; however using (\ref{ap46}), estimate (\ref{ap25}) will be sufficient for reconstruction in the full range (\ref{ap43}).
Note that in fact both alternatives in \eqref{ap25} hold, no matter what the sign of $\alpha-|\gamma|_p+\frac{D}{2}$ is.

\medskip

Following (\ref{ap69}) in Subsection \ref{ss:triang2}, we distinguish the cases
of $\gamma$ purely polynomial and $\gamma$ not pp, see Subsection \ref{ss:pop} for the language. As in Subsection \ref{ss:triang2},
we first treat $\gamma$ not pp by a purely algebraic argument in this subsection.
We then treat $\gamma$ pp by an analytic argument in Subsection \ref{ss:3pIII}.

\medskip

\begin{proof}[Proof of \eqref{ap25} \textnormal{\&} \eqref{ap46} for $\gamma \neq$ pp (algebraic argument)]
\begin{sloppypar} 
The argument for (\ref{ap25}) for $\gamma\not=$ pp relies on the fact that
(\ref{temp05}) \& (\ref{temp06}) transmit to $\Gamma_{x\;x+y}^*$, 
and that (\ref{temp10}) \& (\ref{temp09}) transmit to the endomorphism 
${S:={\rm d}\Gamma_{x+y}^*-{\rm d}\Gamma_x^*\Gamma_{x\;x+y}^*}$. 
As a consequence, as for (\ref{ap19}) in Subsection \ref{ss:order1}, we have
\begin{align}
&S^\gamma_\beta=\mbox{linear combination of}\quad
S^{\delta_{{\bf n}_1}}_{\beta_1}
(\Gamma_{0\;x+y}^*)^{\delta_{{\bf n}_2}}_{\beta_2}\cdots
(\Gamma_{0\;x+y}^*)^{\delta_{{\bf n}_l}}_{\beta_l}\label{ap50bis}\\
&\mbox{where}\;k\delta_3+\beta_1+\cdots+\beta_l=\beta\;\;\mbox{if}\;\;
k\delta_3+\delta_{{\bf n}_1}+\cdots+\delta_{{\bf n}_l}=\gamma.\label{ap49}
\end{align}
\end{sloppypar}
We note that by (\ref{ao10}) and (\ref{temp08})
\begin{align}\label{ap74}
\begin{array}{cc}
(|k\delta_3|-|0|)+|\beta_1|+\sum_{j=2}^l(|\beta_j|-|0|)&=|\beta|,\\
(2k-l+1)+[\beta_1]+\sum_{j=2}^l([\beta_j]+1)&=[\beta].
\end{array}
\end{align}
Since $\gamma$ is populated but not purely polynomial we have $k\not=0$. Hence we learn from the first row 
in (\ref{ap74}) and (\ref{temp08}) that $|\beta_1|<|\beta|$. 
Since $\gamma$ is populated, we have in particular $0\le[\gamma]+1=2k-l+1$, cf \eqref{ao10}.
Since by (\ref{cw47}) for $j=2,\dots,l$, the r.~h.~s.~of (\ref{ap50bis}) vanishes unless 
$\beta_j$ is populated, and thus $[\beta_j]+1\ge 0$.
Thus we learn from the second row in (\ref{ap74}) that $[\beta_1]\le[\beta]$.
Hence $|\beta_1|_{\prec}<|\beta|_{\prec}$ by definition (\ref{ap64}).
Based on \eqref{ap17} rather than \eqref{cw47}, we also learn $|\beta_j|_{\prec}<|\beta|_{\prec}$ for $j=2,\dots,l$, so that we may use our induction hypothesis. 

\medskip

From (\ref{ao65ter}) (with $p$ replaced by a suitable exponent $p_j > p$) 
and from (\ref{ap25}) in its form (\ref{ap57}) for $\gamma$ pp, both in their induction hypothesis version
(with $q$ replaced by some exponent $q'>q$),
we obtain by H\"older's inequality in probability (provided
$\frac{1}{q'}+\sum_{j=2}^l\frac{1}{p_j}=\frac{1}{q}$) applied to (\ref{ap50bis}) 
an estimate by
\begin{align}\label{ap50}
&\min\{|y|^{\alpha-|{\bf n}_1|+\frac{D}{2}}(|y|+R)^{|\beta_1|-\alpha},
(|y|+R)^{|\beta_1|-|{\bf n}_1|+\frac{D}{2}}\}\nonumber\\
&\times (|y|+R)^{|\beta_2|-|{\bf n}_2|}\cdots(|y|+R)^{|\beta_l|-|{\bf n}_l|}.
\end{align}
According to (\ref{cw44}) we have $|\beta_j|\ge|{\bf n}_j|\ge 0$ for $j=2,\dots,l$, 
so that (\ref{ap50}) is $\le$ 
\begin{align}\label{ap76}
\min\{|y|^{\alpha-\sum_{j=1}^l|{\bf n}_j|+\frac{D}{2}}(|y|+R)^{\sum_{j=1}^l|\beta_j|-\alpha},
(|y|+R)^{\sum_{j=1}^l|\beta_j|-\sum_{j=1}^l|{\bf n}_j|+\frac{D}{2}}\} , \nonumber\\
\quad\mbox{and furthermore}\quad{\textstyle\sum_{j=1}^l}|\beta_j|-\alpha\ge 0.
\end{align}
By definitions (\ref{ao07}) and (\ref{ap48}) of $|\cdot|$ and $|\cdot|_p$, we read off (\ref{ap49}) that
\begin{align*}
\begin{array}{cl}
|\beta|-\alpha&=\sum_{j=1}^l|\beta_j|-l\alpha+2k(1+\alpha),\\
|\gamma|-\alpha&=\sum_{j=1}^l|{\bf n}_j|-l\alpha+2k(1+\alpha),
\quad\sum_{j=1}^l|{\bf n}_j|=|\gamma|_p.
\end{array}
\end{align*}
These identities show that the exponents in (\ref{ap76}) coincide with
those on the r.~h.~s.~of (\ref{ap25}). This also shows that the second
item in (\ref{ap46}) is a consequence of the second item in (\ref{ap76}). 
Finally, if $l = \sum_{\bf n} \gamma ({\bf n}) = 0$ then $\gamma = k \delta_3$ whence $S_{\beta}^{\gamma}=0$ by \eqref{temp06} and \eqref{temp09}, proving the first item in (\ref{ap46}).
\end{proof}


\subsection{Reconstruction for \texorpdfstring{$\delta\Pi^{-}-{\rm d}\Gamma_x^*\Pi_x^{-}$}{increments of deltaPi-}}
\label{ss:reconstr}

We return to our task of deriving \eqref{ap21ter},
which we shall derive from the continuity condition
\begin{align}\label{ap53}
\big(\int_{B_R}dx&\mathbb{E}^\frac{2}{q}|({\rm d}\Gamma_{x+y}^*\Pi_{x+y}^{-}
-{\rm d}\Gamma_x^*\Pi_x^{-})_{\beta r}(x+y)|^{q}\big)^\frac{1}{2}\nonumber\\
&\lesssim
r^{2\alpha}(|y|+r)^{\alpha+\frac{D}{2}}(|y|+r+R)^{|\beta|-2-3\alpha},
\end{align}
with the implicit understanding that the last exponent is non-negative unless the l.~h.~s.~vanishes.
There are now three factors on the r.~h.~s.~of \eqref{ap53} due to the presence of the three length scales $r$, $|y|$, and $R$.
To obtain \eqref{ap53} we combine \eqref{ap52} below with \eqref{ap25} and \eqref{ao65bisimpov}. 
One can check that the exponents in \eqref{ap53} are attained for $\beta=\delta_3$ and $\gamma=\delta_3+\delta_{\mathbf{0}}$ in \eqref{ap52} and the estimate is thus sharp.
We provide the reconstruction argument proper that establishes \eqref{ap21ter} 
based on the assumption \eqref{ap43} and on \eqref{ap53} in Subsection~\ref{ss:recIIIrev}, see also \cite{BL23b} for a similar reconstruction theorem in Besov spaces.

\begin{proof}[Proof of \eqref{ap53} (continuity in the base-point)]
Recalling (\ref{ap34bis}) we have 
${\rm d}\Gamma_{x+y}^*\Pi_{x+y}^{-}-{\rm d}\Gamma_x^*\Pi_x^{-}$
$=({\rm d}\Gamma_{x+y}^*-{\rm d}\Gamma_x^*\Gamma_{x\;x+y}^*)\Pi_{x+y}^{-}$.
We consider the $\beta$-component, mollify on scale $r$, and evaluate in $x+y$:
\begin{align}\label{ap52}
\lefteqn{({\rm d}\Gamma_{x+y}^*\Pi_{x+y}^{-}-{\rm d}\Gamma_x^*\Pi_x^{-})_{\beta r}(x+y)}\nonumber\\
&=\sum_{\gamma}
({\rm d}\Gamma_{x+y}^*-{\rm d}\Gamma_x^*\Gamma_{x\;x+y}^*)_\beta^\gamma\Pi_{x+y \gamma r}^{-}(x+y).
\end{align}
Because of the strict triangularities (\ref{ap15}) \& (\ref{ap16}), the induction hypothesis (\ref{ao65bisimpov}) is sufficient.
According to (\ref{temp17}),
the sum is effectively restricted to populated $\gamma$'s that are non-purely polynomial
so that we may appeal to the part of (\ref{ap25}) we established in Subsection \ref{ss:algIII}.
Hence the first r.~h.~s.~factor in (\ref{ap52}) is estimated by $|y|^{\alpha-|\gamma|_p+\frac{D}{2}}
(|y|$ $+R)^{|\beta|+|\gamma|_p-|\gamma|-\alpha}$ and by
$(|y|+R)^{|\beta|-|\gamma|+\frac{D}{2}}$. As in Subsection \ref{ss:algIII}, we combine the
estimate of both factors by H\"older's inequality to obtain that the l.~h.~s.~of (\ref{ap53})
is estimated by terms of the form
\begin{align}\label{ap80}
r^{|\gamma|-2}\min\{|y|^{\alpha-|\gamma|_p+\frac{D}{2}}(|y|+R)^{|\beta|+|\gamma|_p-|\gamma|-\alpha},
(|y|+R)^{|\beta|-|\gamma|+\frac{D}{2}}\}.
\end{align}
By the first item in (\ref{ap46}) we have
for the populated and non-pp $\gamma$'s effectively appearing in (\ref{ap52}) 
that $\gamma(3)\not=0$ and thus
\begin{align}\label{ap79}
|\gamma|-|\gamma|_p\stackrel{\eqref{ao07}}{=}
(1+\alpha)2\gamma(3)+(-\alpha)(\sum_{\bf n}\gamma({\bf n})-1)
\stackrel{\eqref{ck28}}{\ge}2(1+\alpha).
\end{align}
Hence we have in particular $|\gamma|-2\ge2\alpha$ so that
we obtain in case of $\alpha-|\gamma|_p+\frac{D}{2}$ $\ge 0$ that
the term in (\ref{ap80}) is $\le$
\begin{align*}
\lefteqn{r^{|\gamma|-2}|y|^{\alpha-|\gamma|_p+\frac{D}{2}}
(|y|+R)^{|\beta|+|\gamma|_p-|\gamma|-\alpha}}\nonumber\\
&\stackrel{\hphantom{\eqref{ap46},\eqref{ap79}}}{\le} r^{2\alpha}(|y|+r)^{-\alpha-|\gamma|_p+\frac{D}{2}+|\gamma|-2}
(|y|+R)^{|\beta|+|\gamma|_p-|\gamma|-\alpha}\nonumber\\
&\stackrel{\eqref{ap46},\eqref{ap79}}{\le}r^{2\alpha}(|y|+r)^{\alpha+\frac{D}{2}}
(|y|+r+R)^{|\beta|-2-3\alpha},
\end{align*}
as desired. In the remaining case of $\alpha-|\gamma|_p+\frac{D}{2}$ $\le 0$
we have by (\ref{ap79}) that $|\gamma|-2$ $\ge3\alpha+\frac{D}{2}$, so that by
the alternative estimate in (\ref{ap80}) 
and $|\beta|-|\gamma|+\frac{D}{2}\geq0$ (which is a consequence of \eqref{cw44} and \eqref{ap70}), 
we obtain that the term is $\le$
\begin{align*}
r^{|\gamma|-2}(|y|+R)^{|\beta|-|\gamma|+\frac{D}{2}}
&\le r^{3\alpha+\frac{D}{2}}(|y|+r+R)^{|\beta|-2-3\alpha},\\
&\le r^{2\alpha}(|y|+r)^{\alpha+\frac{D}{2}}(|y|+r+R)^{|\beta|-2-3\alpha},
\end{align*}
again as desired.
\end{proof}


\subsection{Integration for \texorpdfstring{$\delta\Pi-{\rm d}\Gamma_x^*\Pi_x$}{increments of deltaPi}}\label{ss:intIII}

It is not possible to directly pass from the estimate (\ref{ap21bis}) to 
the estimate (\ref{ap61}) via the PDE (\ref{ap27}). 
The reason is more algebraic than analytic, as we shall explain now: (\ref{ap27})
is not sufficient to even characterize $(\delta\Pi-{\rm d}\Gamma_x^*\Pi_x)_\beta$
in terms of ${(\delta\Pi^{-}-{\rm d}\Gamma_x^*\Pi_x^{-})_\beta}$ -- even when taking
the vanishing (\ref{ck96bis}) in $x$ into account. 
The reason is that there is a mismatch between the
vanishing order $2+s$ of $\delta\Pi_\beta$ $-\sum_{\gamma}({\rm d}\Gamma_x^*)_\beta^\gamma
\Pi_{x \gamma}$, and the growth of the same quantity at infinity: 
In view of (\ref{ap70}), the sum extends over all $\gamma$ with $|\gamma|<|\beta|+\frac{D}{2}$,
and thus includes terms that grow almost at order $2+\frac{D}{2}$, which is definitely
larger than $2+s$. Hence we need to truncate the sum at order $2+s$. 
For the same reason, we need to restrict to $|\beta|<\lceil 2+s\rceil$.
This means that instead of relying on (\ref{ap27}) we work with
\begin{align*}
\lefteqn{L(\delta\Pi_{\beta}-\sum_{|\gamma|<2+s}
({\rm d}\Gamma_x^*)_\beta^\gamma\Pi_{x \gamma})}\nonumber\\
&=\delta\Pi_{\beta}^{-}-\sum_{|\gamma|<2+s}({\rm d}\Gamma_x^*)_\beta^\gamma\Pi_{x \gamma}^{-}
\quad\mbox{ mod polynomials of degree}<s.
\end{align*}
This in turn requires to pre-process the input (\ref{ap21bis})
and to post-process the output (\ref{ap61}): We need an independent argument that ensures
\begin{align}
\big(\int_{B_R}dx\mathbb{E}^\frac{2}{q}|\sum_{|\gamma|\ge 2+s}
({\rm d}\Gamma_x^*)_\beta^\gamma\Pi_{x \gamma r}^{-}(x)|^q\big)^\frac{1}{2}&\lesssim r^{\alpha-2+\frac{D}{2}}
(r+R)^{|\beta|-\alpha},\label{ap71bis}\\
\big(\int_{B_R}dx\mathbb{E}^\frac{2}{q}|\sum_{|\gamma|\ge 2+s}
({\rm d}\Gamma_x^*)_\beta^\gamma\Pi_{x \gamma r}(x)|^q\big)^\frac{1}{2}&\lesssim r^{\alpha+\frac{D}{2}}
(r+R)^{|\beta|-\alpha}.\label{ap72bis}
\end{align}
In this subsection, we claim that (\ref{ap71bis}) \& (\ref{ap72bis}) follow once we establish 
\begin{align}\label{ap73bis}
\mbox{for populated $\gamma\not=\mbox{pp}$}:\quad
\big(\int_{B_R}dx\mathbb{E}^\frac{2}{q}|({\rm d}\Gamma_x^*)_\beta^\gamma|^{q}\big)^\frac{1}{2}
\lesssim R^{\frac{D}{2}+|\beta|-|\gamma|};
\end{align}
by (\ref{ap70}) the exponent is effectively positive.
Establishing 
(\ref{ap73bis}) requires a second round of algebraic argument,
integration, and three-point argument, see Subsection \ref{ss:second_round}.
On the other hand, the integration argument proper will be carried out in Subsection~\ref{ss:intIIIrevisited}.
\begin{proof}[Proof that \eqref{ap73bis} implies \eqref{ap71bis} \textnormal{\&} \eqref{ap72bis}]
We start by noting that the restriction to populated $\gamma$ that are not
purely polynomial is sufficient for (\ref{ap71bis}) \& (\ref{ap72bis}).
In case of (\ref{ap71bis}), this follows from (\ref{temp17});
in case of (\ref{ap72bis}),  this is a consequence of (\ref{cw59}).
Next, we note that by the strict triangularity (\ref{ap16}), the sum in $\gamma$
is restricted to $|\gamma|_\prec<|\beta|_\prec$ so that we may
appeal to the induction hypothesis in form of (\ref{ao65impov}) \& (\ref{ao65bisimpov}).
Finally, by (\ref{ap70}), the (finite) sum is restricted to $|\gamma|$ $<|\beta|+\frac{D}{2}$,
next to $|\gamma|$ $\ge 2+s$ $=\alpha+\frac{D}{2}$.
Hence by the triangle inequality and H\"older's inequality, (\ref{ap71bis}) is 
as desired estimated by
\begin{align*}
\max_{\alpha+\frac{D}{2}\le|\gamma|<|\beta|+\frac{D}{2}}
R^{\frac{D}{2}+|\beta|-|\gamma|}r^{|\gamma|-2}
\le r^{\alpha-2+\frac{D}{2}}(r+R)^{|\beta|-\alpha}.
\end{align*}
The argument for (\ref{ap72bis}) is very similar. 
\end{proof}


\subsection{The three-point argument for \texorpdfstring{${\rm d}\Gamma_{x+y}^*-{\rm d}\Gamma_x^*\Gamma_{x\,x+y}^*$}{increments of dGamma}}\label{ss:3pIII}

In order to close the induction step, we need to establish (\ref{ap25}) 
for purely polynomial $\gamma$, that is
\begin{align}\label{ap57}
\big(\int_{B_R}dx\mathbb{E}^\frac{2}{q}&|({\rm d}\Gamma^*_{x+y}
-{\rm d}\Gamma_x^*\Gamma_{x\;x+y}^*)_\beta^{\delta_{\bf n}}|^{q}\big)^\frac{1}{2}\nonumber\\
&\lesssim\left\{\begin{array}{lc}|y|^{\alpha-|{\bf n}|+\frac{D}{2}}(|y|+R)^{|\beta|-\alpha}
&\mbox{for}\;\alpha-|{\bf n}|+\frac{D}{2}>0\\
(|y|+R)^{|\beta|-|{\bf n}|+\frac{D}{2}}&\mbox{else}
\end{array}\right\}.
\end{align}
This will follow from estimating an appropriate norm of the polynomial the coefficients of which are given by $(\mathrm{d}\Gamma^*_{x+y}-\mathrm{d}\Gamma^*_x\Gamma^*_{x\,x+y})_\beta^{\delta_\mathbf{n}}$.
Indeed, we use (\ref{ap34bis}) to write
$({\rm d}\Gamma_{x+y}^*-{\rm d}\Gamma_x^*\Gamma_{x\;x+y}^*)\Pi_{x+y}$
$=-(\delta\Pi-{\rm d}\Gamma_{x+y}^*\Pi_{x+y})$
$+(\delta\Pi-{\rm d}\Gamma_{x}^*\Pi_{x})$; we consider the $\beta$ component,
spell out the matrix vector product and split into purely polynomial $\gamma$,
on which we use (\ref{ck38}), and the remainder to deduce the identity
\begin{align}\label{eq:3pt}
&\sum_{\bf n}({\rm d}\Gamma_{x+y}^*-{\rm d}\Gamma_x^*\Gamma_{x\;x+y}^*)_\beta^{\delta_{\bf n}}(\cdot-x-y)^{\bf n} \nonumber\\
&=(\delta\Pi-{\rm d}\Gamma_{x}^*\Pi_{x})_{\beta}(\cdot)
-(\delta\Pi-{\rm d}\Gamma_{x+y}^*\Pi_{x+y})_{\beta}(\cdot) \nonumber\\
&-\sum_{\gamma\not=\textnormal{pp}}
({\rm d}\Gamma_{x+y}^*-{\rm d}\Gamma_x^*\Gamma_{x\;x+y}^*)_\beta^{\gamma}\Pi_{x+y\gamma}(\cdot) ,
\end{align}
which involves the three points $x$, $x+y$, and an active variable $(\cdot)$, hence the name three-point-argument.

\begin{proof}[Proof of \eqref{ap57} (three-point argument)]
Let us start with the second part of the estimate, which by (\ref{ap40}) amounts
to the case of $|{\bf n}|\ge s+2$, so that by (\ref{cw59}), the l.~h.~s.~
reduces to $-\sum_{\gamma}({\rm d}\Gamma_x^*)_\beta^\gamma$ 
$(\Gamma_{x\;x+y}^*)_\gamma^{\delta_{\bf n}}$. Since by (\ref{ap16}), 
the sum is restricted to $|\gamma|_\prec<|\beta|_\prec$, we may appeal to the induction
hypothesis (\ref{ao65ter}) in conjunction with (\ref{ap47}) for the second factor.
If $\gamma$ were purely polynomial, i.~e.~of the form $\gamma=\delta_{\bf m}$,
then (\ref{cw44}) would in conjunction with (\ref{temp04}) imply 
$|{\bf m}|\ge|{\bf n}|$, so that also $|{\bf m}|\ge s+2$, and
that by (\ref{cw59}), the first factor would vanish. Hence effectively $\gamma\not=$ pp
so that for the first factor, we may appeal to (\ref{ap73bis}). Moreover, by
(\ref{cw44}) and (\ref{ap70}), the finite sum restricts to  
$|\gamma|\ge|{\bf n}|$ and $|\gamma|<|\beta|+\frac{D}{2}$.
In conclusion, we obtain by
H\"older's inequality an estimate of the l.~h.~s.~of (\ref{ap57}) by
\begin{align*}
\max_{|{\bf n}|\le|\gamma|<|\beta|+\frac{D}{2}} 
R^{|\beta|-|\gamma|+\frac{D}{2}}|y|^{|\gamma|-|{\bf n}|}\le(|y|+R)^{|\beta|-|{\bf n}|+\frac{D}{2}},
\end{align*}
as desired.

\medskip

We now turn to the first estimate in (\ref{ap57}), which we derive from (\ref{ap61}): 
We apply $(\cdot)_r(x+y)$ to \eqref{eq:3pt}, which yields the representation
\begin{align}\label{ap60}
\lefteqn{\big(\sum_{\bf n}({\rm d}\Gamma_{x+y}^*
-{\rm d}\Gamma_x^*\Gamma_{x\;x+y}^*)_\beta^{\delta_{\bf n}}(\cdot )^{\bf n}\big)_r(0)}\nonumber\\
&=(\delta\Pi-{\rm d}\Gamma_{x}^*\Pi_{x})_{\beta r}(x+y)
-(\delta\Pi-{\rm d}\Gamma_{x+y}^*\Pi_{x+y})_{\beta r}(x+y)
\nonumber\\
&-\sum_{\gamma\not=\textnormal{pp}}
({\rm d}\Gamma_{x+y}^*-{\rm d}\Gamma_x^*\Gamma_{x\;x+y}^*)_\beta^{\gamma}\Pi_{x+y\gamma r}(x+y).
\end{align}
We then take the $(\int_{B_R}dx\mathbb{E}^\frac{2}{q}|\cdot|^{q})^\frac{1}{2}$ norm.
In order to subsume the first r.~h.~s.~term under \eqref{ap61} we introduce 
$\psi^{(r)}:=\psi(\cdot+R^{-1}y)$, 
with $R^{-1}$ being the inverse of the transformation \eqref{ck29},
so that $\psi_r(\cdot+y)$ $=\psi^{(r)}_r$, and note that as long as $|y|\le r$,
the semi-norms (\ref{ck52a}) of $\psi^{(r)}$ are controlled by those of $\psi$.
Hence we obtain from (\ref{ap61}) 
\begin{align*}
&\big(\int_{B_R}dx\mathbb{E}^\frac{2}{q}|(\delta\Pi
-{\rm d}\Gamma_{x}^*\Pi_{x})_{\beta r}(x+y)|^{q}\big)^\frac{1}{2}\nonumber\\
&\lesssim r^{\alpha+\frac{D}{2}}(r+R)^{|\beta|-\alpha}\quad\mbox{provided}\;|y|\le r.
\end{align*}
For the second r.~h.~s.~term in (\ref{ap60}) 
we use the generic estimate
\begin{align}\label{ap59}
\int_{B_R}dx f^2(x+y)\le\int_{B_{R+|y|}}dx'f^2(x')
\end{align}
together with \eqref{ap61} to
obtain like for the first term
\begin{align*}
&\big(\int_{B_R}dx\mathbb{E}^\frac{2}{q}|(\delta\Pi
-{\rm d}\Gamma_{x+y}^*\Pi_{x+y})_{\beta r}(x+y)|^{q}\big)^\frac{1}{2}\nonumber\\
&\lesssim r^{\alpha+\frac{D}{2}}(r+R)^{|\beta|-\alpha}\quad\mbox{provided}\;|y|\le r.
\end{align*}
For the last term in (\ref{ap60}) we note that at this stage of
the induction step we do have access to (\ref{ap25}) for $\gamma$ not purely polynomial;
by (\ref{ap15}) \& (\ref{ap16}), only $|\gamma|_\prec<|\beta|_\prec$ are involved so that
we may appeal to (\ref{ao65impov}) on the level of the induction hypothesis.
By the triangle inequality, \eqref{temp30} and H\"older's inequality we obtain,
still for $|y|\le r$,
\begin{align*}
\lefteqn{\big(\int_{B_R}dx\mathbb{E}^\frac{2}{q}|\big(\sum_{\bf n}({\rm d}\Gamma_{x+y}^*
-{\rm d}\Gamma_x^*\Gamma_{x\;x+y}^*)_\beta^{\delta_{\bf n}}(\cdot )^{\bf n}\big)_r(0)|^{q}
\big)^\frac{1}{2}}\nonumber\\
&\lesssim r^{\alpha+\frac{D}{2}}(r+R)^{|\beta|-\alpha}
+\sum_{\substack{\gamma\colon |\gamma|_\prec < |\beta|_\prec \\ |\gamma|\ge|\gamma|_p}}
|y|^{\alpha-|\gamma|_p+\frac{D}{2}}(r+R)^{|\beta|-|\gamma|+|\gamma|_p-\alpha}
r^{|\gamma|},
\end{align*}
where the restriction to $|\gamma|\ge|\gamma|_p$ follows by definitions \eqref{ao07} and  (\ref{ap48}) from 
the first item in (\ref{ap46}).
We use this for $r=|y|$, in which it simplifies to
\begin{align*}
\big(\int_{B_R}dx\mathbb{E}^\frac{2}{q}&\big|\int d\hat x\psi(\hat x)\sum_{\bf n}|y|^{|{\bf n}|}
({\rm d}\Gamma_{x+y}^*
-{\rm d}\Gamma_x^*\Gamma_{x\;x+y}^*)_\beta^{\delta_{\bf n}}\hat x^{\bf n}\big|^{q}
\big)^\frac{1}{2}\nonumber\\
&\lesssim|y|^{\alpha+\frac{D}{2}}(|y|+R)^{|\beta|-\alpha}.
\end{align*}

\medskip 

We may conclude by a duality argument: Let $(F_x)_x$ be an arbitrary random field, then by duality, denoting by $q^*$ the H\"older conjugate of $q$,
	\begin{align*}
		& \int d \hat{x} \psi (\hat{x}) \sum\limits_{\bf n} \hat{x}^{\bf n} |y|^{|{\bf n}|} \int_{B_R} dx \mathbb{E} F_x ({\rm d}\Gamma_{x+y}^* -{\rm d}\Gamma_x^*\Gamma_{x\;x+y}^*)_\beta^{\delta_{\bf n}} \\
			& \lesssim |y|^{\alpha+\frac{D}{2}}(|y|+R)^{|\beta|-\alpha}
\big(\int_{B_R}dx\mathbb{E}^\frac{2}{q^*}|F_x|^{q^*}\big)^\frac{1}{2} .
	\end{align*}
The estimate depends on the arbitrary Schwartz function $\psi$ only from its Schwartz semi-norms, thus the right-hand-side is an upper bound on some norm of the (deterministic) polynomial
	\begin{align*}
		\hat{x} \mapsto \sum\limits_{\bf n} \hat{x}^{\bf n} |y|^{|{\bf n}|} \int_{B_R} dx \mathbb{E} F_x ({\rm d}\Gamma_{x+y}^* -{\rm d}\Gamma_x^*\Gamma_{x\;x+y}^*)_\beta^{\delta_{\bf n}} .
	\end{align*}
By equivalence of norms, also the coefficients of this polynomial are bounded by the same right-hand-side, i.~e.\ for all ${\bf n}$,
	\begin{align*}
		& |y|^{|{\bf n}|} \int_{B_R} dx \mathbb{E} F_x ({\rm d}\Gamma_{x+y}^* -{\rm d}\Gamma_x^*\Gamma_{x\;x+y}^*)_\beta^{\delta_{\bf n}} \\
		& \lesssim |y|^{\alpha+\frac{D}{2}}(|y|+R)^{|\beta|-\alpha}
\big(\int_{B_R}dx\mathbb{E}^\frac{2}{q^*}|F_x|^{q^*}\big)^\frac{1}{2} .
	\end{align*}
Since the random field $(F_x)_x$ was arbitrary we have obtained the first line in \eqref{ap57}.
\end{proof}


\subsection{Logical order of the proof}\label{ss:logic}

The round of the four arguments from Subsections \ref{ss:algIII}, \ref{ss:reconstr},
\ref{ss:intIII}, and \ref{ss:3pIII} is logically not complete: On the one hand,
the argument for the estimate (\ref{ap73bis}) on $({\rm d}\Gamma^*)_\beta^\gamma$
is still missing. On the other hand,
we still need to establish the estimates on $(\Pi_\beta,\Pi^{-}_\beta,(\Gamma^*)_\beta^\gamma)$ 
itself.
This requires two more rounds of the same sequence of four arguments,
in a specific logical order depicted in Table~\ref{t:tasks}. 

\medskip

The second round provides the estimates on $(\delta\Pi_\beta,\delta\Pi^{-}_\beta,$ 
$({\rm d}\Gamma^*)_\beta^\gamma)$; it is carried out in Subsection \ref{ss:second_round}. 
Like the first round,
it starts with an {\sc algebraic argument} to establish the estimate (\ref{ap73bis})
on $({\rm d}\Gamma^*)_\beta^\gamma$ for $\gamma$ not purely polynomial. 
Based on this, it passes from (\ref{ap21ter}) established in the first round
to an estimate of $\delta\Pi^{-}_\beta$ itself.
It then appeals to {\sc integration} to come to an estimate of $\delta\Pi_\beta$,
which in turn allows for a {\sc three-point argument} to upgrade the estimate on 
$({\rm d}\Gamma^*)_\beta^\gamma$ to one for all $\gamma$, including the purely polynomial ones.

\medskip

The third round finally yields the estimates on $(\Pi_\beta,\Pi^{-}_\beta,$ 
$({\Gamma^*})_\beta^\gamma)$ of Theo\-rem \ref{th:main};
it is carried out in Subsection \ref{ss:third_round}. Like the previous rounds,
it starts with an {\sc algebraic argument} to estimate $(\Gamma^*)^\gamma_\beta$
for ${\gamma\not=\text{pp}}$. It then proceeds to an estimate of $\Pi^{-}_\beta$,
distinguishing the cases\footnote{since the induction proceeds via $\prec$ and not
$|\cdot|$, the tow cases are intertwined} of $|\beta|<2$ and $|\beta|>2$. 
In case of $|\beta|<2$, we appeal to the spectral gap inequality and use the estimate
of $\delta\Pi^{-}_\beta$ established in the second round, 
and an estimate of $\mathbb{E}\Pi^{-}_\beta$.
In case of $|\beta|>2$, we appeal to another {\sc reconstruction} argument. 
We then use {\sc integration} to estimate $\Pi_\beta$, and finally  
a {\sc three-point argument} to estimate $({\rm d}\Gamma^*)_\beta^\gamma$ for ${\gamma=\text{pp}}$.

\newcommand\tikzmark[2]{%
\tikz[remember picture,baseline] \node[inner sep=2pt,outer sep=0] (#1){#2};%
}
\newcommand\link[2]{%
\begin{tikzpicture}[remember picture, overlay, >=stealth, shift={(0,0)}]
  \draw[->] (#1) to (#2);
\end{tikzpicture}%
}

{
\renewcommand{\arraystretch}{2}
\begin{table}[ht]
\caption{The different tasks} \label{t:tasks}
\centering
\resizebox{14cm}{!}{
\newcolumntype{M}[1]{>{\centering\arraybackslash}m{#1}}
\begin{tabular}{| M{3.5cm} || M{2cm} | M{2cm} | M{5cm} |@{}m{0pt}@{}} 
 \hline
 	 & Model itself & Malliavin derivative & Increments of Malliavin derivative&\\
 	 \hline\hline
	 \textsc{Algebraic argument} & \tikzmark{aa1}{$(\Gamma_x^*)_{\beta}^{\gamma\neq\text{pp}}$} & \tikzmark{aa2}{$({\rm d}\Gamma_x^*)_{\beta}^{\gamma\neq\text{pp}}$} & \tikzmark{aa3}{$({\rm d}\Gamma_{x+y}^*-{\rm d}\Gamma_x^*\Gamma_{x\;x+y}^*)_{\beta}^{\gamma\neq\text{pp}}$} & \\
	\hline
	\textsc{Reconstruction argument} & \tikzmark{r1}{$\Pi_{x \beta}^{-}$} & \tikzmark{r2}{$\delta\Pi_{x \beta}^{-}$} & \tikzmark{r3}{$(\delta\Pi_{x \beta}^{-}-{\rm d}\Gamma_x^*\Pi_x^-)_{\beta}$} &\\
	\hline
	\textsc{Integration argument} & \tikzmark{i1}{$\Pi_{x \beta}$} & \tikzmark{i2}{$\delta\Pi_{x \beta}$} & \tikzmark{i3}{$(\delta\Pi_{x \beta}-{\rm d}\Gamma_x^*\Pi_x)_{\beta}$} &\\
	\hline
	\textsc{3-Point argument} & \tikzmark{3p1}{$(\Gamma_x^*)_{\beta}^{\gamma=\text{pp}}$} & \tikzmark{3p2}{$({\rm d}\Gamma_x^*)_{\beta}^{\gamma=\text{pp}}$} & \tikzmark{3p3}{$({\rm d}\Gamma_{x+y}^*-{\rm d}\Gamma_x^*\Gamma_{x\;x+y}^*)_{\beta}^{\gamma=\text{pp}}$} &\\
	\hline
\end{tabular}
\link{aa1}{r1}
\link{r1}{i1}
\link{i1}{3p1}
\link{aa2}{r2}
\link{r2}{i2}
\link{i2}{3p2}
\link{aa3}{r3}
\link{r3}{i3}
\link{i3}{3p3}
\link{r2}{r1}
\link{r3}{r2}
}
\end{table}}


\subsection{A second round of algebraic argument, reconstruction, integration, and three-point argument \texorpdfstring{to estimate $(\delta\Pi,\delta\Pi^{-},$ ${\rm d}\Gamma^*)$}{}}\label{ss:second_round}

More precisely, the tasks of this subsection are:
\begin{itemize}
\item Based on the induction hypothesis in form of \eqref{ap78}, 
we establish (\ref{ap73bis}) by an {\sc algebraic argument}.
\item We post-process (\ref{ap21bis}) to 
\begin{align}\label{ap75}
\mathbb{E}^\frac{1}{q}|\delta\Pi^{-}_{\beta r}(x)|^{q}
\lesssim r^{\alpha-2-\frac{D}{2}}(r+|x|)^{|\beta|-\alpha+\frac{D}{2}}.
\end{align}
This task also contains the base case.
\item By {\sc integration}, we pass from (\ref{ap75}) to 
\begin{align}\label{ap77}
\mathbb{E}^\frac{1}{q}|\delta\Pi_{\beta r}(0)|^{q}
\lesssim r^{|\beta|}.
\end{align}
\item By a {\sc three-point argument}, we pass from (\ref{ap77}) to
the version of (\ref{ap73bis}) for purely polynomial $\gamma$ 
\begin{align}\label{ap78}
\big(\int_{B_R}dx\mathbb{E}^\frac{2}{q}
|({\rm d}\Gamma_x^*)_\beta^{\delta_{\bf n}}|^{q}\big)^\frac{1}{2}
\lesssim R^{\frac{D}{2}+|\beta|-|\delta_{\bf n}|}.
\end{align}
\end{itemize}

\medskip

\begin{proof}[Proof of \eqref{ap73bis} \textnormal{\&} \eqref{ap75} \textnormal{\&} \eqref{ap77} \textnormal{\&} \eqref{ap78}]
We just point out the differences with the previous subsections:
For (\ref{ap73bis}), we start from the re\-presentation (\ref{ap19}) and argue as in Subsection~\ref{ss:algIII}
that we may appeal to the induction hypothesis (\ref{ao65ter}) and (\ref{ap78}).
We thus obtain that the l.~h.~s.~of (\ref{ap73bis}) is $\lesssim$
$R^{\frac{D}{2}+|\beta_1|-|\delta_{{\bf n}_1}|}$
$R^{|\beta_2|-|\delta_{{\bf n}_2}|}\cdots R^{|\beta_l|-|\delta_{{\bf n}_l}|}$.
Using the additivity (\ref{temp08}), we learn from (\ref{cw71bis}) that the
sum of exponents is $=\frac{D}{2}+|\beta|-|\gamma|$.

\medskip

Turning to \eqref{ap75}, 
we first apply \eqref{ap59} to $f=\mathbb{E}^\frac{1}{q}|\delta\Pi^-_{\beta r}|^{q}$
to deduce 
\begin{align*}
\big(\int_{B_R}dx\mathbb{E}^\frac{2}{q}
|\delta\Pi^{-}_{\beta r}(x+y)|^{q}\big)^\frac{1}{2}
\lesssim \big(\int_{B_{R+|y|}}dx\mathbb{E}^\frac{2}{q}
|\delta\Pi^{-}_{\beta r}(x)|^{q}\big)^\frac{1}{2} .
\end{align*}
We rewrite
$\delta\Pi^-_{\beta r}(x)
= (\delta\Pi^- - {\rm d}\Gamma^*_{x}\Pi^-_{x})_{\beta \, r} (x)
+ ({\rm d}\Gamma^*_{x}\Pi^-_{x})_{\beta \, r}(x)$.
On the first r.~h.~s.~term we use \eqref{ap21bis}, 
while on the second r.~h.~s.~term we apply the same argument that led from \eqref{ao65bisimpov} \& (\ref{ap73bis}) to (\ref{ap71bis}) to obtain 
\begin{align}\label{ap75averaged}
\big(\int_{B_R}dx\mathbb{E}^\frac{2}{q}
|\delta\Pi^{-}_{\beta r}(x+y)|^{q}\big)^\frac{1}{2}
\lesssim r^{\alpha-2}(r+R+|y|)^{|\beta|-\alpha+\frac{D}{2}}.
\end{align}
This averaged estimate can be post processed into the pointwise (\ref{ap75}), 
the analytic details are provided in Subsection~\ref{ss:smearing}.

\medskip

For the integration of the PDE (\ref{ck70bis}) in order to pass from (\ref{ap75}) to (\ref{ap77}), 
we do not run into the problem of Subsection~\ref{ss:intIII}: The order of growth  
and the order of vanishing both agree with the non-integer $|\beta|$.
A detailed integration argument is provided in Subsection~\ref{ss:intIrevisited}.

\medskip

We now turn to the induction step for (\ref{ap78}) and start from the representation
(obtained analogously as in (\ref{ap60}))
\begin{align*}
\lefteqn{\big(\sum_{\bf n}({\rm d}\Gamma_x^*)_\beta^{\delta_{\bf n}}(\cdot)^{\bf n}\big)_r(0)}\\
&=\delta\Pi_{\beta r}(x)
- (\delta\Pi-{\rm d}\Gamma_x^*\Pi_x)_{\beta r}(x)
-\sum_{\gamma\not=\textnormal{pp}}({\rm d}\Gamma_x^*)_\beta^\gamma\Pi_{x\gamma r}(x) .
\end{align*}
Arguing as in the proof of \eqref{ck62} above, the estimate (\ref{ap77}) remains true with $0$ replaced by $x$ with $|x|\leq r$,
which we use for the first r.~h.~s.~term.
For the second r.~h.~s.~term we appeal to \eqref{ap61}.
According to (\ref{ap16}), the sum in the third r.~h.~s.~term restricts to $|\gamma|_\prec < |\beta|_\prec$, so that we
may appeal to the induction hypothesis (\ref{ao65impov}). 
For the first factor we use \eqref{ap73bis}, so that 
by (\ref{ap70}) this yields for $r\leq R$ the estimate
\begin{align*}
\lefteqn{\big(\int_{B_R}dx\mathbb{E}^\frac{2}{q}
|(\sum_{\bf n}({\rm d}\Gamma_x^*)_\beta^{\delta_{\bf n}}(\cdot)^{\bf n})_r(0)|^{q}
\big)^\frac{1}{2}}\nonumber\\
&\lesssim R^{\frac{D}{2}} r^{\alpha-\frac{D}{2}} (r\hspace{-.3ex}+\hspace{-.5ex}R)^{|\beta|-\alpha+\frac{D}{2}}
+ r^{\alpha+\frac{D}{2}}(r\hspace{-.3ex}+\hspace{-.5ex}R)^{|\beta|-\alpha}
+\max_{|\gamma|<|\beta|+\frac{D}{2}}R^{\frac{D}{2}+|\beta|-|\gamma|}r^{|\gamma|}.
\end{align*}
We use this for $r=R$ and change variables according to $y=R\hat y$ to the effect of
\begin{align*}
\big(\int_{B_R}dx\mathbb{E}^\frac{2}{q}
|\int d\hat y\psi(\hat y)\sum_{\bf n}R^{|{\bf n}|}({\rm d}\Gamma_x^*)_\beta^{\delta_{\bf n}}
\hat y^{\bf n}|^{q}
\big)^\frac{1}{2}\lesssim R^{|\beta|+\frac{D}{2}}.
\end{align*}
As in Subsection \ref{ss:3pIII}, and recalling (\ref{temp04}), this yields (\ref{ap78}).
\end{proof}

\subsection{A third round of algebraic argument, reconstruction, integration, and three-point argument \texorpdfstring{to estimate $(\Pi,\Pi^{-},\Gamma^*)$}{}}\label{ss:third_round}

More precisely, the tasks of this subsection are: 
\begin{itemize}
\item Based on the induction hypothesis, we establish \eqref{ao65ter} 
for $\gamma$ not purely polynomial by an \textsc{algebraic argument}.
\item For the estimate \eqref{ao65bisimpov} of $\Pi^-_\beta$ we distinguish two cases:\\
-- For $|\beta|>2$ we appeal to a simple \textsc{reconstruction} argument.\\
-- For $|\beta|<2$ the estimate is a consequence of the control of the expectation \eqref{cw76} and of the Malliavin derivative \eqref{ap75}, followed by an application of the spectral gap inequality \eqref{ap20ter} as outlined in Subsection~\ref{ss:spectral_gap}. 
\item By \textsc{integration}, we pass from \eqref{ao65bisimpov} to \eqref{ao65impov}.
\item By a \textsc{three-point argument}, we pass from \eqref{ao65impov} to the version of \eqref{ao65ter} for purely polynomial $\gamma$.
\end{itemize}
\begin{proof}[Proof of \eqref{ao65impov} \textnormal{\&} \eqref{ao65bisimpov} \textnormal{\&} \eqref{ao65ter}]
To obtain \eqref{ao65ter} for $\gamma$ not purely polynomial,
we proceed analogously to the argument that led to \eqref{ap73bis}, 
just replacing \eqref{ap19} by the simpler \eqref{ap09}. 

\medskip

For \eqref{ao65bisimpov} it remains to provide an argument for $|\beta|>2$. 
As in Subsection~\ref{ss:reconstr} this is the consequence of 
continuity in the base-point in form of 
\begin{align*}
\mathbb{E}^\frac{1}{p}|(\Pi_{x+y}^{-}-\Pi_{x}^{-})_{\beta r}(x+y)|^p
\lesssim r^{\alpha-2}(r+|y|)^{|\beta|-\alpha},
\end{align*}
and the qualitative 
\begin{align*}
\lim_{r\downarrow0} \Pi^-_{x\beta \, r}(x) = 0,
\end{align*}
the analytic details will be provided in Subsection~\ref{ss:recIIIrev} (in the more involved setting of Subsection~\ref{ss:reconstr}).
The latter display is a consequence of \eqref{t10} in combination with $|\beta|>2$.
For the continuity in the base-point we appeal to \eqref{temp01} in form of 
$\Pi^-_{x+y}-\Pi^-_x = ({\rm id}-\Gamma^*_{x\, x+y})\Pi^-_{x+y}$; 
by the population \eqref{temp17} of $\Pi^-$ it is enough 
to appeal to the already established \eqref{ao65ter} for $\gamma$ not purely polynomial, 
and by the triangularity \eqref{cw44} of $\Gamma^*$, we conclude with the induction hypothesis of \eqref{ao65bisimpov}.

\medskip

The (analytic) details on the integration argument leading from \eqref{ao65bisimpov} to \eqref{ao65impov} are provided in Subsection~\ref{ss:intIrevisited}.

\medskip

The three-point argument yielding \eqref{ao65ter} for purely polynomial $\gamma$ 
proceeds analogous to the one leading to \eqref{ap78}, starting from the identity
\begin{align*}
\sum_{{\bf n}}(\Gamma^*_x)_\beta^{e_{\bf n}}(\cdot-x)^{\bf n}=\Pi_{\beta}
-\sum_{\gamma\neq\textnormal{pp}}(\Gamma^*_x)_\beta^\gamma\Pi_{x\gamma} 
\end{align*}
which is a consequence of \eqref{temp02} and \eqref{ck38}
\end{proof}


\section{Proof details}\label{s:proofs}

\subsection{Semi-group convolution} \label{ss:semigroup}

Both for reconstruction and integration it is convenient to work with
a specific convolution kernel $\Psi$ in (\ref{ck30}), namely the kernel
of the semi-group generated by the positive operator $L^*L$ of Fourier symbol
$|q|^4$, cf.~(\ref{ap37}). Since ${\mathcal F}\Psi(q):=\exp(-|q|^4)$
is a Schwartz function, $\Psi$ is a Schwartz function. By definition (\ref{ck30})
of the rescaling we have ${\mathcal F}\Psi_r(q)$ $=\exp(-t|q|^4)$ provided $t=r^4$,
which motivates the short-hand notation
\begin{align}\label{ap92}
\Psi_t:=\Psi_{r=\sqrt[4]{t}}\quad\mbox{such that}\quad\partial_t\Psi_t+L^*L\Psi_t=0.
\end{align}
Since $\Psi_t*\Psi_T$ $=\Psi_{t+T}$ (as can be easily inferred on the Fourier level)
we have for any Schwartz distribution $F$
\begin{align}\label{ck30ter}
(F_{t})_T=F_{t+T}\quad\mbox{where}\quad F_t(x):=(F,\Psi_t(x-\cdot));
\end{align}
the latter definition is analogous\footnote{note however that \eqref{ck30ter} is slightly different from \eqref{ck30} because $\Psi_t$ is not the rescaling of $\Psi$ at scale $t$, but at scale $\sqrt[4]{t}$, recall \eqref{ap92}; still, below we will use the following convention when using the subscripts $t$ and $r$: $(\cdot)_t$ will refer to \eqref{ck30ter} while $(\cdot)_r$ will refer to \eqref{ck30}} to (\ref{ck30}). 


\subsection{Details on reconstruction for \texorpdfstring{$\delta\Pi^{-}-{\rm d}\Gamma_x^*\Pi_x^{-}$}{increments of deltaPi-}} \label{ss:recIIIrev}

We will thus assume \eqref{ap53} with $\psi_r$ replaced by $\Psi_t$, 
and establish \eqref{ap21ter} with $\psi_r$ replaced by $\Psi_t$. 
In Subsection~\ref{ss:change} we will argue that this 
implies \eqref{ap21ter} for an arbitrary kernel $\psi$. 

\begin{proof}[Proof that\footnotemark $\eqref{ap53}_t$ implies $\eqref{ap21ter}_t$]
\footnotetext{from now on, by the notation $\eqref{ap53}_t$ we mean the estimate $\eqref{ap53}$ with $\psi_r$ replaced by $\Psi_t$}
\begin{sloppypar} 
Since\footnote{for $\tau\leq t\leq T$, in what follows $(\cdot)_\tau$, $(\cdot)_t$ and $(\cdot)_T$ always refer to the semi-group convolution as defined in \eqref{ck30ter}} $\Psi$ and thus $\Psi_\tau$ is normalized, i.~e.~$\int dx\Psi_\tau$ $=1$,
we obtain from (\ref{ck97}) the qualitative information that
$\delta\Pi^{-}_\beta(x)$ $=\lim_{\tau\downarrow 0}({\rm d}\Gamma_x^*\Pi_x^{-})_{\beta \tau}(x)$.
Hence introducing the notation $(EF)(x)$ $:=F_x(x)$ for the diagonal evaluation of
our family ${\{F_x}$ ${:=({\rm d}\Gamma_x^*\Pi_x^{-})_\beta\}_x}$
of germs, we see that $\eqref{ap21ter}_T$ follows once we establish
\begin{align}\label{ap55}
\big(\int_{B_R}dx&\mathbb{E}^\frac{2}{q}|(EF_{\tau}
-F_{x\tau})_{T-\tau}(x)|^{q}\big)^\frac{1}{2}\nonumber\\
&\lesssim (\sqrt[4]{T})^{3\alpha+\frac{D}{2}}(\sqrt[4]{T}+R)^{|\beta|-2-3\alpha}
\quad\mbox{for}\;\tau\le T.
\end{align}
\end{sloppypar}
This is an estimate of the commutator between the evaluation operator
$E$ and the mollification operator $(\cdot)_\tau$, 
with the understanding that the mollification acts only on the active variable 
but not on the base-point when applied to $F$. 
It is here where we leverage the semi-group property \eqref{ck30ter}. 
Restricting $\tau$ to be a dyadic fraction of $T$ allows us to write the l.~h.~s.~of \eqref{ap55} as a telescoping sum over dyadic length scales:
\begin{align*}
(E F_{\tau} - F_{x \tau})_{T-\tau}(x)
= \sum_{\substack{\tau\leq t < T, \; \\ t \,\textnormal{dyadic fraction of $T$}}}
\big( (E F_{t})_{t} - E F_{2t} \big)_{T-2t} (x).
\end{align*}
Hence the claim follows once we establish
\begin{align}\label{ap54}
\big(\int_{B_R}dx&\mathbb{E}^\frac{2}{q}|((EF_{t})_{t}
-E F_{2t})_{T-2t}(x)|^{q}\big)^\frac{1}{2}\nonumber\\
&\lesssim (\sqrt[4]{t})^{3\alpha+\frac{D}{2}}(\sqrt[4]{T}+R)^{|\beta|-2-3\alpha}
\quad\mbox{for}\;t\le T/2. 
\end{align}
Indeed, since $3\alpha+D/2>0$ by assumption (\ref{ap43}), the r.~h.~s.~of (\ref{ap54}) gives rise to a convergent
geometric series as $\tau\downarrow 0$ that sums up to the r.~h.~s.~of (\ref{ap55}).

\medskip

We are thus left with establishing \eqref{ap54}, 
which is an easy consequence of (\ref{ap53}):
We write $(EF_{t}-F_{x' t})_{t}(x')$ 
$=\int dy'\Psi_{t}(y') (F_{x'-y'} -F_{x'})_{t}(x'-y')$, 
and thus $((EF_{t})_{t}-EF_{2t})_{T-2t}(x)$ 
$=\int dy\Psi_{T-2t}(y)$ $\int dy'$ $\Psi_{t}(y')$ $(F_{x-y-y'}$ $-F_{x-y})_{t}(x-y-y')$.
By the triangle inequality and using \eqref{ap59}
we learn 
from (\ref{ap53}) with $(\sqrt[4]{t},-y',R+|y|)$ playing the role of $(r,y,R)$ that the l.~h.~s.~of (\ref{ap54})
is estimated by 
\begin{align*}
(\sqrt[4]{t})^{2\alpha}\hspace{-1ex}
\int \hspace{-1ex}dy|\Psi_{T\hspace{-.3ex}-\hspace{-.3ex}2t}(y)| 
\hspace{-.5ex}\int \hspace{-1ex}dy'|\Psi_{t}(y')| (|y'| \hspace{-.5ex}+\hspace{-.5ex}\sqrt[4]{t})^{\alpha+\frac{D}{2}} (|y'|\hspace{-.5ex}+\hspace{-.5ex}\sqrt[4]{t}\hspace{-.5ex}+\hspace{-.5ex}R\hspace{-.5ex}+\hspace{-.5ex}|y|)^{|\beta|-2-3\alpha}\hspace{-.5ex}.
\end{align*}
Recalling \eqref{ck52a} and that the exponents $\alpha+\frac{D}{2}$ and $|\beta|-2-3\alpha$ are (effectively) non-negative, integrating against $\Psi_{T - 2t}(y)$, resp.~$\Psi_t(y')$ amounts to replacing $| y |$ by $\sqrt[4]{T - 2t}$ resp.~$|y'|$ by $\sqrt[4]{t}$ in the integral above, so that it can be absorbed in the r.~h.~s.~of (\ref{ap54}).
\end{proof}


\subsection{Details on the expectation \texorpdfstring{$\mathbb{E} \Pi_{\beta t}^-(0)$ for $|\beta |<2$}{E Pi beta t - for beta < 2}} \label{ss:expectation}

We claim that \eqref{ck52ter} implies 
\begin{align}\label{cw76}
|\mathbb{E}\Pi_{\beta t}^{-}(0)|\lesssim (\sqrt[4]{t})^{|\beta|-2}.
\end{align}

\begin{proof}
By \eqref{ck52ter},
it suffices to establish
\begin{align}\label{cw26}
t|\frac{d}{dt}\mathbb{E}\Pi_{\beta t}^{-}(0)|\lesssim (\sqrt[4]{t})^{|\beta|-2}.
\end{align}
By re-expansion \eqref{temp01} we have for any $\tau$,
\begin{align*}
\Pi_{\beta\tau}^{-}(x)=\sum_{\gamma}(\Gamma_{x}^*)_\beta^\gamma\Pi^{-}_{x\gamma\tau}(x) ,
\end{align*}
so that using the semi-group property (\ref{ck30ter})
in form of $(\cdot)_t$ $=(\cdot)_{t-\tau}(\cdot)_\tau$ we obtain
\begin{align*}
\Pi_{\beta t}^{-}(0)=\sum_{\gamma}\int dx\Psi_{t-\tau}(x)
(\Gamma_{x}^*)_\beta^\gamma\Pi^{-}_{x\gamma\tau}(x).
\end{align*}
Since, by stationarity (\ref{temp30}), $\mathbb{E}\Pi^{-}_{x\gamma\tau}(x)$ does not depend on $x$,
and since $\int dx\Psi_{t-\tau}(x)$ does not depend on $t$,
this yields the representation
\begin{align*}
\frac{d}{dt}\mathbb{E}\Pi_{\beta t}^{-}(0)=\sum_{\gamma}\int dx\partial_t\Psi_{t-\tau}(x)
\mathbb{E}(\Gamma_{x}^*-{\rm id})_\beta^\gamma\Pi^{-}_{x\gamma\tau}(x).
\end{align*}
Now appealing to the strict triangularity \eqref{ap15} of $\Gamma^* - {\rm id}$, in this sum effectively $\gamma \prec \beta$.
Thus, by the recursive estimates \eqref{ao65bisimpov} \textnormal{\&} \eqref{ao65ter} on $\Pi^-$ and $\Gamma^*$, we obtain 
\begin{align*}
t|\frac{d}{dt}\mathbb{E}\Pi_{\beta t}^{-}(0)|
\lesssim\sqrt[4]{\tau}^{\alpha-2}(\sqrt[4]{\tau}+\sqrt[4]{t-\tau})^{|\beta|-\alpha},
\end{align*}
which yields the desired (\ref{cw26}) when choosing $\tau=\frac{t}{2}$.
\end{proof}

\subsection{Change of kernel}\label{ss:change}

\begin{sloppypar} 
Subsections~\ref{ss:recIIIrev} and~\ref{ss:expectation} output estimates with respect to the semi-group kernel $\Psi$ introduced in Subsection~\ref{ss:semigroup}.
We need to upgrade them into estimates with respect to general Schwartz kernels $\psi$. 
This will be achieved via the
following representation formula valid for any Schwartz distribution $F$:
\begin{align}\label{ap91}
F_r(x)&=\sum_{j=0}^k\frac{1}{j!}\int dy((L^*L)^j\psi)_r(-y)F_{t=r^4}(x+y)\nonumber\\
&+\frac{1}{k!}\int_0^{r^4}\frac{dt}{t}
\Big(\frac{t}{r^4}\Big)^{k+1}\int dy((L^*L)^{k+1}\psi)_r(-y)F_t(x+y),
\end{align}
\end{sloppypar}
where the role of the arbitrary integer $k\ge0$ 
is to make the $t$-integral concentrate near $t=r^4$. 
We learn from (\ref{ap91}) that indeed $F_r(x)$ can be written as a linear combination
of $F_t(x+y)$ with essentially $t\sim r^4$ and $|y|\lesssim r$
(as we shall see e.~g.\ in the proof of \eqref{ap21ter} later in this subsection).
\begin{proof}[Proof of\eqref{ap91}]
The argument for
(\ref{ap91}) is straight-forward: By definitions (\ref{ck30}) and (\ref{ck30ter}) it
reduces to
\begin{align*}
\psi_r&=\sum_{j=0}^k\frac{1}{j!}((L^*L)^j\psi)_{r\,t=r^4}
+\frac{1}{k!}\int_0^{r^4}\frac{dt}{t}
\Big(\frac{t}{r^4}\Big)^{k+1}((L^*L)^{k+1}\psi)_{r\,t},
\end{align*}
where $(\cdot)_{r t}$ $:=((\cdot)_r)_t$ stands short for first applying the rescaling and then
the semi-group convolution, which commutes to $((\cdot)_{\hat t})_r$ where ${t=r^4\hat t}$.
Hence by a change of variables of the $t$-integral, and removing the $r$-rescaling,
the above identity follows from
\begin{align*}
\psi&=\sum_{j=0}^k\frac{1}{j!}(L^*L)^j\psi_{\hat t=1}
+\frac{1}{k!}\int_0^{1}\frac{d\hat t}{\hat t} \, 
\hat t^{k+1}(L^*L)^{k+1}\psi_{\hat t}.
\end{align*}
Because of the second item in (\ref{ap92}) in form of $(L^*L)^j\psi_{\hat t}$
$=(-\frac{d}{d\hat t})^j\psi_{\hat t}$ this follows from integration by parts. 
\end{proof}

\medskip

As a consequence of \eqref{ap91}, let us argue that the output of Subsection~\ref{ss:recIIIrev}, namely
	\begin{align}\label{ap22}
&\nonumber \big(\int_{B_R}dx\mathbb{E}^\frac{2}{q}
|(\delta\Pi^{-}-{\rm d}\Gamma_x^*\Pi_x^{-})_{\beta t}(x + y)|^{q}\big)^\frac{1}{2} \\
& \lesssim (\sqrt[4]{t})^{2 \alpha} ( \sqrt[4]{t} + | y | )^{\alpha + \frac{D}{2}}(\sqrt[4]{t}+|y|+R)^{|\beta|-2-3\alpha},
	\end{align}
implies \eqref{ap21ter}.

\begin{proof}
Let $k \geq 0$, to be adjusted later.
Denote $F_x := (\delta\Pi^{-}-{\rm d}\Gamma_x^*\Pi_x^{-})_{\beta}$, as well as $\hat\psi:=\sum_{j = 0}^k \frac{1}{j!}(L^*L)^j\psi$ and $\check{\psi}:=\frac{1}{(k+1)!}(L^*L)^{k+1}\psi$,
then by \eqref{ap91} and the triangle inequality, 
	\begin{align*}
		& \big(\int_{B_R}dx\mathbb{E}^\frac{2}{q}
|F_{x r} (x)|^{q}\big)^\frac{1}{2} \\
		& \lesssim \int dy\,\big|\hat\psi_r(-y)\big|\big(\int_{B_R}dx\mathbb{E}^\frac{2}{q} |F_{t = r^4} (x+y)|^{q}\big)^\frac{1}{2} \\
			& \quad + \int_0^{r^4}\frac{dt}{t}
\Big(\frac{t}{r^4}\Big)^{k+1}\int dy\,\big|\check{\psi}_r(-y)\big|\big(\int_{B_R}dx\mathbb{E}^\frac{2}{q} |F_{t} (x+y)|^{q}\big)^\frac{1}{2} .
	\end{align*}
Now we plug in \eqref{ap22}, recalling that the exponents $\alpha+\frac{D}{2}$ and $|\beta|-2-3\alpha$ therein are (effectively) non-negative, so that integrating against $\hat{\psi}_r$, $\check{\psi}_r$ amounts to replacing $| y |$ by $r$ in the right-hand-side of \eqref{ap22}
	\begin{align*}
		\big(\int_{B_R}dx\mathbb{E}^\frac{2}{q}
|F_{x r} (x)|^{q}\big)^\frac{1}{2}
		& \lesssim r^{3 \alpha + \frac{D}{2}} ( r + R )^{|\beta|-2-3\alpha} \\
		& \quad+\int_0^{r^4}\frac{dt}{t}
\Big(\frac{t}{r^4}\Big)^{k+1} (\sqrt[4]{t})^{2 \alpha} r^{\alpha + \frac{D}{2}}(r+R)^{|\beta|-2-3\alpha} .
	\end{align*}
Now it suffices to fix $k$ large enough so that the latter integral converges at $t = 0$, namely $k > -1 - \frac{\alpha}{2}$, yielding \eqref{ap21ter} after integration\footnote{note that by the assumption \eqref{ck28} on $\alpha$, choosing $k = 0$ was sufficient for this argument (but is not in general)}.
\end{proof}

With the same argument, we also may pass from the estimates \eqref{ao65bisimproved} 
established against the semi-group kernel $\Psi$, to the same estimates uniformly over bounded Schwartz kernels.


\subsection{Details on reconstruction for \texorpdfstring{$\delta\Pi^{-}$}{delta Pi-}}\label{ss:smearing}

In this subsection we post-process the averaged estimate \eqref{ap75averaged} 
into the pointwise \eqref{ap75}. 
This relies on the following annealed version of Sobolev's inequality, valid for any $y \in \mathbb{R}^{1 + d}$, $k > D/2$, $R > 0$, and any
(smooth)
random field $u$
	\begin{align}\label{ap95}
		\mathbb{E}^{\frac{1}{q}} | u ( y ) |^q
		& \lesssim \sum_{| \mathbf{n} | \leq k} R^{| \mathbf{n} | - \frac{D}{2}} \Big( \int_{B_R} d x \, \mathbb{E}^{\frac{2}{q}} | \partial^{\mathbf{n}} u ( x + y ) |^{q} \Big)^{\frac{1}{2}} ,
	\end{align}
where the implicit multiplicative constant depends only on $k, D$, and
which we establish now by a duality argument.
\begin{proof}[Proof of \eqref{ap95}]
Up to replacing $u$ by $u ( y + R\cdot )$, it suffices to prove the inequality when $y = 0$, $R = 1$. 
Let $q^*$ be the H\"older dual exponent to $q$, and let $F$
be an arbitrary random variable
with $\mathbb{E}^{1/q^{*}}|F|^{q^*} \leq 1$.
We apply the (standard, anisotropic) Sobolev inequality to the function $\bar{u}: x \mapsto \mathbb{E} [ u ( x ) F ]$,
to the effect of 
	\begin{align*}
		| \bar{u} ( 0 ) |^2
		& \lesssim \sum_{| \mathbf{n} | \leq k} \int_{B_1} \big| \partial^{\mathbf{n}} \bar{u} \big|^2  ,
	\end{align*}
leading by H\"older in probability to
	\begin{align*}
		| \mathbb{E} [ u ( 0 ) F ] |^2
		& \lesssim \sum_{|\mathbf{n}| \leq k} \int_{B_1} \mathbb{E}^{\frac{2}{q}} | \partial^{\mathbf{n}} u |^{q},
	\end{align*}
which yields the desired \eqref{ap95} by duality since the random variable $F$ was arbitrary.
\end{proof}

\begin{proof}[Proof that $\eqref{ap75averaged}$ implies $\eqref{ap75}$]
We apply \eqref{ap95} to $u := \delta \Pi_{\beta r}^-$, $R = r$, and $k$ being the smallest integer $>D/2$.
Let $\mathbf{n} \in \mathbb{N}_0^{1 + d}$ with $| \mathbf{n}| \leq k$.
Then
$\partial^{\mathbf{n}} \delta \Pi^-_{\beta r} = r^{- | \mathbf{n} |} \delta \Pi^-_{\beta} * \tilde{\psi}_r$ for the new Schwartz function $\tilde{\psi} = \partial^{\mathbf{n}} \psi$ whose Schwartz semi-norms \eqref{ck52a} are bounded by those of $\psi$.
Thus, by the assumption \eqref{ap75averaged}, the $\mathbf{n}$-th summand in \eqref{ap95} is bounded by 
$r^{| \mathbf{n}| - D/2 - |\mathbf{n} | + \alpha - 2} ( r + | y | )^{| \beta | - \alpha + D/2}$, 
which is the r.~h.~s.~ of \eqref{ap75} (with $y$ in place of $x$), as desired.
\end{proof}


\subsection{Abstract integration} \label{ss:abstract_int}

In preparation for estimating $\Pi$, $\delta\Pi$, and $\delta\Pi-{\rm d}\Gamma^*_x\Pi_x$ given 
estimates of $\Pi^-$, $\delta\Pi^-$, and $\delta\Pi^--{\rm d}\Gamma^*_x\Pi^-$, 
we give an abstract integration result. 

\medskip

We claim that there is a family $(\mu_{[\psi, r, \kappa]})_{\psi \in \mathcal{S}, r > 0, \kappa \in \mathbb{R} \setminus \mathbb{N}_0}$
of measures on\footnote{here we denote by $\mathcal{S}$ the space of Schwartz functions defined by the family of semi-norms \eqref{ck52a}} $\mathcal{S} \times (0, \infty)$, such that

    \begin{itemize}
    \item (Representation)
    For each bounded\footnote{w.~r.~t.~the family of semi-norms \eqref{ck52a}} set $\mathcal{B} \subset \mathcal{S}$ there is another bounded set $\tilde{\mathcal{B}} \subset \mathcal{S}$ with the following property.
	Let $u, f$ be any deterministic (Schwartz)
	distributions such that $L u = f$ modulo a polynomial of degree $\leq\kappa-2$ and
	\begin{align}\label{mb06}
        \sup_{r > 0} r^{2 - \kappa} | f * \psi_r ( 0 ) | < \infty ,
        \qquad \sup_{r > 0} r^{- \kappa} | u * \psi_r ( 0 ) | < \infty ,
    \end{align}
uniformly over $\psi$ in bounded sets in Schwartz space.
Then 
            \begin{align} \label{mb07}
                u * \psi_r ( 0 ) = \int_{\tilde{\mathcal{B}}\times(0,\infty)} f * \tilde{\psi}_{\tilde{r}} ( 0 ) \, d \mu_{[\psi, r, \kappa]} ( \tilde{\psi}, \tilde{r} ) ,
            \end{align}
            for all $\psi \in \mathcal{B}$, $r>0$.

        \item (Moment bounds) One has
            \begin{align} \label{mb03}
                \int_{\mathcal{S} \times (0, \infty)} \tilde{r}^{\kappa - 2} \, d \mu_{[\psi, r, \kappa]} ( \tilde{\psi}, \tilde{r} ) \lesssim r^{\kappa} ,
            \end{align}
            uniformly over $r > 0$ and $\psi$ in bounded sets in Schwartz space.
            In fact, $\mu$ depends on $\kappa$ only through its integer part $\lfloor \kappa \rfloor$.
Furthermore, \eqref{mb03} remains true with $\kappa$ replaced by any $\tilde{\kappa}$ provided $\lfloor \kappa \rfloor = \lfloor \tilde{\kappa} \rfloor$.
    \end{itemize}

\begin{proof}[Proof of \eqref{mb07} \textnormal{\&} \eqref{mb03}]
The proof relies on the representation formula
    \begin{align} \label{mb01}
        u * \psi_r ( 0 ) = \int_{0}^{\infty} d t \, \big( (\mathrm{id} - {\rm T}_0^{\kappa}) (L^* f * \Psi_t) \big) * \psi_r ( 0 ) ,
    \end{align}
where\footnote{we recall that $\Psi_t$ denotes the semigroup  generated by the symmetric operator $LL^*$, see Subsection~\ref{ss:semigroup}} ${\rm T}_x^{\kappa}f$ denotes the Taylor polynomial of $f$ in the
base-point $x$ of (parabolic) order $\leq\kappa$.
We justify this representation at the end of this proof, and start by estimating the right-hand side of \eqref{mb01}.

\medskip

As is typical for Schauder-type arguments,
the proof distinguishes between the ``near-field'' range $t\leq r^4$ 
and the ``far-field'' range $t\geq r^4$.
For the former we treat the contributions from $\rm id$ and ${\rm T}_0^{\kappa}$ separately, 
while for the latter we appeal to the Taylor remainder in integral form which we briefly discuss now. 
Assume first $\kappa>0$.
Fix $x \in \mathbb{R}^{1 + d}$, recall the notation $S x := ( s^2 x_0, s x_1, \cdots, s x_d )$, and consider the auxiliary function
$[0, 1]$ ${\ni s\mapsto g ( s ):=(\mathrm{id}-{\rm T}_0^{\kappa})(L^* f * \Psi_t) ( S x )}$, the derivatives of which vanish at zero up to order $\kappa$ so that by Taylor's representation,
$g ( 1 )$ ${= \int_0^1 d s \frac{(1 -s)^{k - 1}}{(k - 1)!} \frac{d^k g}{ds^k}}$, where $k$ is the smallest integer $> \kappa$.
We note that for some (generic) coefficients $c_{{\bf n}}$,
	\begin{align*}
		\frac{d^k}{ds^k} = \sum\limits_{| {\bf n} | \geq k , \sum_i n_i \leq k} c_{{\bf n}} \, s^{| {\bf n} | - k} x^{\bf n} \partial^{\bf n} ,
	\end{align*}
whence the representation
	\begin{align*}
		& (\mathrm{id}-{\rm T}_0^{\kappa})(L^* f * \Psi_t) ( x ) \\
		& = \sum\limits_{| {\bf n} | \geq k , \sum_i n_i \leq k} c_{{\bf n}} \, x^{\bf n} \int_{0}^1 d s \, (1 - s)^{k - 1} s^{| {\bf n} | - k}  \partial^{\bf n} L^* f * \Psi_t ( S x) .
	\end{align*}
In the case $\kappa<0$ the Taylor remainder $({\rm id}-{\rm T}_0^\kappa)$ simply reduces to ${\rm id}$. 
We thus rewrite the right-hand side of \eqref{mb01} as 
    \begin{align*} 
        & = \int_{0}^{r^4} d t \, L^* f * \Psi_t * \psi_r ( 0 ) \\
        & \quad - \sum_{| {\bf n} | \leq \kappa} \frac{1}{{\bf n} !} \int_{0}^{r^4} d t \, \partial^{\bf n} L^* f * \Psi_t ( 0 ) \int d x \, x^{\bf n} \psi_r ( x )  \\
        & \quad +
        \sum_{|{\bf n}|>\kappa, \sum_i n_i \leq\kappa+1} r^{| {\bf n} |} \int_{r^4}^\infty dt \int_0^1 ds \, c_{\bf n}(s) ( \cdot^{\bf n} \psi)_{s r} * \partial^{\bf n} L^* f * \Psi_t ( 0 ) \\
        & \quad + \mathbf{1}_{(-\infty,0)}(\kappa) \int_{r^4}^\infty dt\, L^* f * \Psi_t * \psi_r (0) \\
        & =: A_1 - \sum_{| {\bf n} | \leq \kappa} \tfrac{1}{{\bf n}!} A_{2, {\bf n}} + \sum_{|{\bf n}|>\kappa, \sum_i n_i\leq\kappa+1} A_{3, {\bf n}} + A_4 ,
    \end{align*}
    with the understanding that the empty sum equals $0$ and that consequently the second and third contributions are not present for $\kappa<0$.
We deal with each term separately.
We start with $A_1$, which we rewrite as
    \begin{align*}
        A_1 = f * \Big( \int_{0}^{r^4} d t \, L^* \Psi_t * \psi_r \Big) ( 0 ) = r^{2} f * \tilde{\psi}^{[\psi, r]}_r ( 0 ) ,
    \end{align*}
for the Schwartz function
    \begin{align*}
        \tilde{\psi}^{[\psi, r]} 
        = \int_0^{r^4} dt \, r^{D-2} (L^* \Psi_t * \psi_r) (r \cdot) .
    \end{align*}
One may check that the Schwartz semi-norms of $\tilde{\psi}^{[\psi, r]}$ are uniformly bounded by that of $\psi$, so that the claimed representation
    \begin{align*}
        A_1 
        = \int f * \tilde{\psi}_{\tilde{r}} ( 0 ) \, d \mu_{[\psi, r, \kappa]}^1 ( \tilde{\psi}, \tilde{r} ) 
    \end{align*}
holds with the measure
    \begin{align*}
        \mu_{[\psi, r, \kappa]}^{1} ( \tilde{\psi}, \tilde{r} ) = \tilde{r}^2 \delta_{r} ( \tilde{r} ) \, \delta_{\tilde{\psi}^{[\psi, r]}} ( \tilde{\psi} ) .
    \end{align*}
Turning to the moment bound for $\mu_{[\psi, r, \kappa]}$, one readily obtains
    \begin{align*}
        \int_{\mathcal{S} \times (0, \infty)} \tilde{r}^{\kappa - 2} \, d \mu_{[\psi, r, \kappa]}^1 ( \tilde{\psi}, \tilde{r} )
        & = r^{\kappa} ,
    \end{align*}
as desired.

\medskip

We turn to $A_{2, {\bf n}}$, which we rewrite as 
\begin{align*}
    A_{2,{\bf n}} 
    &= \int_0^{r^4} dt \, f * \partial^{\bf n} L^* \Psi_t (0)
    r^{|{\bf n}|} \int dx\, x^{\bf n}\psi(x) \\
    &= \int_0^{r^4} dt \, f * \tilde{\psi}_{\sqrt[4]{t}}^{[\psi,r]}(0) \, 
    r^{|{\bf n}|} (\sqrt[4]{t})^{-|{\bf n}|-2} \\
    &= \int_0^{r} dt \, f * \tilde{\psi}_{t}^{[\psi,r]}(0) \, 
    r^{|{\bf n}|} 4t^{-|{\bf n}|+1} 
\end{align*}
for the Schwartz function 
\begin{align*}
    \tilde{\psi}^{[\psi,{\bf n}]} = \partial^{\bf n} L^* \Psi \int dx \, x^{\bf n} \psi(x) .
\end{align*}
The Schwartz semi-norms of $\tilde{\psi}^{[\psi,{\bf n}]}$ are uniformly bounded by those of $\psi$, and we obtain the representation
\begin{align*}
    A_{2,{\bf n}} 
    = \int f * \tilde{\psi}_{\tilde{r}}(0) \, d\mu_{[\psi,r,\kappa]}^{2,{\bf n}}(\tilde{\psi}, \tilde{r}) 
\end{align*}
with the measure 
\begin{align*}
    \mu_{[\psi,r,\kappa]}^{2,{\bf n}}(\tilde{\psi}, \tilde{r})
    = \mathbf{1}_{(0,r)}(\tilde{r}) \, 
    \delta_{\tilde{\psi}^{[\psi,{\bf n}]}}(\tilde{\psi}) \,
    4 r^{|{\bf n}|} \tilde{r}^{-|{\bf n}|+1} .
\end{align*}
For the moment bound we observe
\begin{align*}
    \int_{\mathcal{S}\times(0,\infty)} \tilde{r}^{\kappa-2} \, d\mu_{[\psi,r,\kappa]}^{2,{\bf n}}(\tilde{\psi} ,\tilde{r}) 
    = \int_0^r d\tilde{r} \, 4 r^{|{\bf n}|} \tilde{r}^{\kappa-|{\bf n}|-1} ,
\end{align*}
which is integrable and bounded by $r^\kappa$ 
due to $|{\bf n}|<\kappa$, 
which in turn is a consequence of the restriction $|{\bf n}|\leq\kappa$ 
and the assumption $\kappa\not\in\mathbb{N}_0$. 

\medskip

We turn to $A_{3, {\bf n}}$, which we rewrite as
\begin{align*}
    A_{3,{\bf n}} 
    &= \int_{r^4}^\infty dt \, f * \Big(r^{|{\bf n}|} \int_0^1 ds \, c_{\bf n}(s) ( \cdot^{\bf n} \psi)_{s r} * \partial^{\bf n} L^* \Psi_t \Big) ( 0 ) \\
    &= \int_{r^4}^\infty dt \, f * \tilde{\psi}_{\sqrt[4]{t}}^{[\psi,r,t,{\bf n}]} (0) \, r^{|{\bf n}|} (\sqrt[4]{t})^{-|{\bf n}|-2} \\
    &= \int_{r}^\infty dt \, f * \tilde{\psi}_{t}^{[\psi,r,t,{\bf n}]} (0) \, r^{|{\bf n}|} 4 (\sqrt[4]{t})^{-|{\bf n}|+1}
\end{align*}
for the Schwartz function
\begin{align*}
    \tilde{\psi}_{\sqrt[4]{t}}^{[\psi,r,t,{\bf n}]} 
    = \int_0^1 ds \, c_{\bf n}(s) (\sqrt[4]{t})^{|{\bf n}|+2} 
    (\cdot^{\bf n}\psi)_{s r} * \partial^{\bf n} L^* \Psi_t ,
\end{align*}
i.~e.
\begin{align*}
    \tilde{\psi}^{[\psi,r,t,{\bf n}]} 
    = \int_0^1 ds \, c_{\bf n}(s) \int dx \, x^{\bf n} \psi(x) 
    (\partial^{\bf n} L^* \Psi)(\cdot - \tfrac{s r}{\sqrt[4]{t}}x) .
\end{align*}
One may check that the Schwartz semi-norms of $\tilde{\psi}^{[\psi,r,t,{\bf n}]}$ are bounded by those of $\psi$ (uniformly when $r^4\leq t$).
We thus obtain the representation
\begin{align*}
    A_{3,{\bf n}}
    = \int f * \tilde{\psi}_{\tilde{r}}(0) \, d\mu_{[\psi,r,\kappa]}^{3,{\bf n}}(\tilde{\psi},\tilde{r})
\end{align*}
with the measure
\begin{align*}
    \mu_{[\psi,r,\kappa]}^{3,{\bf n}} (\tilde{\psi},\tilde{r})
    = \mathbf{1}_{(r,\infty)}(\tilde{r}) \, \delta_{\tilde{\psi}^{[\psi,r,\tilde{r},{\bf n}]}}(\tilde{\psi}) \, 4 r^{|{\bf n}|} \tilde{r}^{-|{\bf n}|+1} .
\end{align*}
The moment bound follows from 
\begin{align}\label{mb05}
    \int_{\mathcal{S}\times(0,\infty)} \tilde{r}^{\kappa-2} \, d\mu_{[\psi,r,\kappa]}^{3,{\bf n}}(\tilde{\psi} ,\tilde{r}) 
    = \int_r^\infty d\tilde{r} \, 4 r^{|{\bf n}|} \tilde{r}^{\kappa-|{\bf n}|-1} ,
\end{align}
which by the restriction of $\kappa<|{\bf n}|$ is integrable and as desired bounded by $r^\kappa$.

\medskip

We turn to $A_4$, which we rewrite as 
\begin{align*}
A_4 &= \int_{r^4}^\infty dt\, f * \big( L^* \Psi_t * \psi_r ) (0) 
= \int_{r^4}^\infty dt\, f * \tilde{\psi}_{\sqrt[4]{t}}^{[\psi,r,t]} (0) (\sqrt[4]{t})^{-2} \\
&= \int_{r}^\infty dt\, f * \tilde{\psi}_{t}^{[\psi,r,t]} (0) 4t
\end{align*}
for the Schwartz function
\begin{align*}
\tilde{\psi}^{[\psi,r,t]} = \int dx\, (L^*\Psi)(\cdot-\tfrac{r}{\sqrt[4]{t}}x) \psi(x) .
\end{align*}
The Schwartz semi-norms of $\tilde{\psi}^{[\psi,r,t]}$ are bounded by those of $\psi$ (uniformly for $r^4\leq t$) and we obtain the representation 
\begin{align*}
A_4 = \int f * \tilde{\psi}_{\tilde{r}} (0) \, d\mu^4_{[\psi,r,\kappa]}(\tilde{\psi},\tilde{r}) 
\end{align*}
with the measure 
\begin{align*}
\mu^4_{[\psi,r,\kappa]}(\tilde{\psi},\tilde{r})
= \mathbf{1}_{(-\infty,0)}(\kappa) \mathbf{1}_{(r,\infty)}(\tilde{r}) \delta_{\tilde{\psi}^{[\psi,r,\tilde{r}]}}(\tilde{\psi}) 4 \tilde{r} .
\end{align*}
The moment bound follows from 
\begin{equation}\label{asdf}
\int_{\mathcal{S}\times(0,\infty)} \tilde{r}^{\kappa-2} \, d\mu^4_{[\psi,r,\kappa]}(\tilde{\psi},\tilde{r})
= \mathbf{1}_{(-\infty,0)}(\kappa) \int_r^\infty d\tilde{r} \, 4\tilde{r}^{\kappa-1},
\end{equation}
which is bounded by $r^\kappa$ as desired.

\medskip

To conclude, let us quickly justify \eqref{mb01}, the r.~h.~s.~ of which we temporarily name $\tilde{u}$.
First, the integral defining $\tilde{u}$ indeed makes sense as a distribution: 
This is because the integrand, when tested against a Schwartz function, 
is bounded as $t \to 0$, and integrable as $t \to \infty$ by virtue of the far-field estimate \eqref{mb05} and \eqref{asdf}.
Now for $0<\tau<T<\infty$, by \eqref{ap92}, 
\begin{align}
	L\int_\tau^T dt({\rm id}-{\rm T}_0^{\kappa})(L^* f * \Psi_t)
	 \nonumber = ({\rm id}-{\rm T}_0^{\kappa - 2}) f * \Psi_\tau - ({\rm id}-{\rm T}_0^{\kappa-2}) f * \Psi_T.
\end{align}
Appealing to the assumptions $\kappa\not\in\mathbb{N}_0$ and the first item of \eqref{mb06} resp.~ to the representation of the Taylor remainder above, we have in the sense of distributions ${\rm T}_0^{\kappa - 2} (f * \Psi_\tau) \to0$ as\footnote{note that in fact ${\rm T}_0^{\kappa - 2} \equiv 0$ when $\kappa<2$} $\tau \to 0$ resp.~ $({\rm id}-{\rm T}_0^{\kappa-2}) (f * \Psi_T) \to0$ as $T \to \infty$.
Thus, $L(u-\tilde{u})$ is a polynomial of degree $\leq\kappa-2$.
We have established just above in this subsection that the second item of \eqref{mb06} holds for $\tilde{u}$, so that by the Liouville argument of Subsection~\ref{ss:unique} 
we deduce $u=\tilde{u}$, as desired.
\end{proof}

\medskip

We now claim that the representation \eqref{mb07} still holds (almost surely) in the case where $u, f$ are random and \eqref{mb06} is replaced by the following annealed version: for some $p > 1$
    \begin{align}\label{mb06b}
       \sup_{r > 0} r^{2 - \kappa} \, \mathbb{E}^{\frac{1}{p}} | f * \psi_r ( 0 ) |^p < \infty ,
        \qquad \sup_{r > 0} r^{- \kappa} \, \mathbb{E}^{\frac{1}{p}} | u * \psi_r ( 0 ) |^p < \infty ,
    \end{align}
uniformly over $\psi$ in bounded sets in Schwartz space.
We argue by duality:
let $A$ be an arbitrary random variable with $\mathbb{E}^{1/p^*} | A |^{p^*} \leq 1$, where $p^* > 1$ is the dual H\"older exponent to $p$.
By H\"older's inequality in probability and \eqref{mb06b}, the assumptions \eqref{mb06} are satisfied with $u$, $f$ replaced by 
	\begin{align*}
		\tilde{f} = \mathbb{E} [A f], 
		\qquad \tilde{u} = \mathbb{E} [A u] .
	\end{align*}
Thus, the representation \eqref{mb07} holds with $\tilde{u}$, $\tilde{f}$ in place of $u, f$.
By \eqref{mb03} and Fubini, this reads
	\begin{align*}
		\mathbb{E} [A u * \psi_r ( 0 ) ] = \mathbb{E} \Big[ A \int_{\tilde{\mathcal{B}}\times(0,\infty)} f * \tilde{\psi}_{\tilde{r}} ( 0 ) \, d \mu_{[\psi, r, \kappa]} ( \tilde{\psi}, \tilde{r} ) \Big] .
	\end{align*}
But since the random variable $A$ was arbitrary, we deduce that $u$ and $f$ enjoy the representation \eqref{mb07} almost surely, as desired.

\subsection{Details on integration for \texorpdfstring{$\Pi$ and $\delta \Pi$}{Pi and delta Pi}} \label{ss:intIrevisited} 

Equipped with the result of Subsection~\ref{ss:abstract_int}, we now prove that\footnote{recall that $(\cdot)_r$ denotes convolution with $\psi_r$ for a generic kernel $\psi$}
\begin{align} \label{cw10}
\mathbb{E}^\frac{1}{p}| \Pi^-_{\beta r}(0)|^p 
\lesssim r^{|\beta|-2}
\end{align}
implies
\begin{align} \label{cw10bis}
\mathbb{E}^\frac{1}{p}| \Pi_{\beta r}(0)|^p 
\lesssim r^{|\beta|} .
\end{align}

\begin{proof}
We note that by \eqref{cw10} in combination with the purely qualitative \eqref{ck52}, the assumption \eqref{mb06b} hold with $f$ replaced by $\Pi_{\beta}^-$, $u$ replaced by $\Pi_{\beta}$ and $\kappa = | \beta| \in \mathbb{R}\setminus\mathbb{N}_0$, recall \eqref{temp04bis}.
Thus by \eqref{ao43ter} we obtain the representation
	 \begin{align*} 
                \Pi_{\beta} * \psi_r ( 0 ) = \int_{\tilde{\mathcal{B}}\times(0,\infty)} \Pi_{\beta}^{-} * \tilde{\psi}_{\tilde{r}} ( 0 ) \, d \mu_{[\psi, r, | \beta |]} ( \tilde{\psi}, \tilde{r} ) ,
     \end{align*}
so that plugging \eqref{cw10} and appealing to the moment bound \eqref{mb03} yields the desired \eqref{cw10bis}.
\end{proof}

The exact same argument allows to pass from $\mathbb{E}^\frac{1}{q}| \delta \Pi^-_{\beta r}(0)|^{q} \lesssim r^{|\beta|-2}$ to $\mathbb{E}^\frac{1}{q}| \delta \Pi_{\beta r}(0)|^{q} \lesssim r^{|\beta|}$.

\subsection{Details on integration for \texorpdfstring{$\delta\Pi-{\rm d}\Gamma_x^*\Pi_x$}{increments of deltaPi}} \label{ss:intIIIrevisited}

The purpose of this Subsection is to argue that \eqref{ap21bis},\eqref{ap73bis},\eqref{ap75}, imply \eqref{ap61}.
In fact, in view of \eqref{ap71bis} \textnormal{\&} \eqref{ap72bis} (which follow from \eqref{ap73bis} as argued in Subsection~\ref{ss:intIII}), we may add to our set of assumptions that 
	\begin{align}
		&\big(\int_{B_R}dx\mathbb{E}^\frac{2}{q}|\big(\delta\Pi_{\beta}^--\sum_{|\gamma|<2+s}
({\rm d}\Gamma_x^*)_\beta^\gamma\Pi_{x \gamma}^-\big)_{r} ( x )|^q\big)^\frac{1}{2} \nonumber \\
&\lesssim r^{\alpha-2+\frac{D}{2}}(r+R)^{|\beta|-\alpha} ,\label{ap93}
	\end{align}
and it suffices to establish
	\begin{align}
		&\big(\int_{B_R}dx\mathbb{E}^\frac{2}{q}|\big(\delta\Pi_{\beta}-\sum_{|\gamma|<2+s}
({\rm d}\Gamma_x^*)_\beta^\gamma\Pi_{x \gamma}\big)_{r} ( x )|^q\big)^\frac{1}{2} \nonumber \\
&\lesssim r^{\alpha+\frac{D}{2}}(r+R)^{|\beta|-\alpha} .\label{ap94}
	\end{align}

\begin{proof}[Proof of \eqref{ap94}]
For notational convenience, let us denote the ``rough-path increments''
	\begin{align*}
		U_x & := \delta \Pi_{\beta} - \sum_{| \gamma | < 2 + s} (\mathrm{d} \Gamma_x^*)_{\beta}^{\gamma} \Pi_{x \gamma} , \\
		F_x & := \delta \Pi_{\beta}^- - \sum_{| \gamma | < 2 + s} (\mathrm{d} \Gamma_x^*)_{\beta}^{\gamma} \Pi_{x \gamma}^- ,
	\end{align*}
so that by definition, $L U_x = F_x$.
The argument is based on the representation formula \eqref{mb07} above, which here takes the form:
provided $| \beta | < \lceil s+2 \rceil$,
for each bounded set $\mathcal{B} \subset \mathcal{S}$ there is another bounded set $\tilde{\mathcal{B}} \subset \mathcal{S}$ such that for all $\psi \in \mathcal{B}$ and $r>0$, 
	\begin{align}\label{mb02}
     	U_x * \psi_r ( x ) = \int_{\tilde{\mathcal{B}}\times(0,\infty)} F_x * \tilde{\psi}_{\tilde{r}} ( x ) \, d \mu_{[\psi, r, 2 + s]} ( \tilde{\psi}, \tilde{r} ) .
     \end{align}
As discussed at the beginning of Subsection~\ref{ss:intIII}, recall \eqref{mb06} and \eqref{mb06b}, this essentially follows from the qualitative vanishing (at $x$) and growth (at infinity), at the same order $2 + s$, of $\delta\Pi_{\beta}^{-}-\sum_{|\gamma|<2+s}({\rm d}\Gamma_x^*)_\beta^\gamma\Pi_{x \gamma}^{-}$.
We refrain from giving a detailed proof of \eqref{mb02} here, let us refer to
\cite[Proposition~4.14]{LOTT21}
where this justification is carried out in the case of a quasi-linear equation.
Here, recall that $| \beta | < \lceil s + 2 \rceil$.
We temporarily make the stronger assumption
	\begin{align} \label{mb04}
		|\beta| < s+2 .
	\end{align}
Appealing on the one hand to \eqref{ap93} when $r \leq R$, and on the other hand splitting the rough-path increment by the triangle inequality in combination with \eqref{ao65bisimpov}, \eqref{ap73bis}, \eqref{ap75}, and \eqref{mb02} when $r \geq R$, we obtain the following estimate valid for all $r, R>0$:
	\begin{align} \label{mb08}
		\Big( \int_{B_R} d x \, \mathbb{E}^{\frac{2}{q}} \big| F_x * \psi_r ( x ) \big|^{q} \Big)^{\frac{1}{2}} & \lesssim r^{s} R^{| \beta | - \alpha} .
	\end{align}
Plugging into \eqref{mb02} in combination with the moment bound \eqref{mb03} we deduce
	\begin{align*}
		\Big( \int_{B_R} d x \, \mathbb{E}^{\frac{2}{q}} \big| U_x * \psi_r ( x ) \big|^{q} \Big)^{\frac{1}{2}}
		\lesssim r^{s + 2} R^{| \beta | - \alpha} ,
	\end{align*}
which absorbs into the desired \eqref{ap94}.
We now turn to the case
	\begin{align*}
		s + 2 \leq | \beta | < \lceil 2 + s \rceil .
	\end{align*}
In that case, arguing as for \eqref{mb08}, we obtain
	\begin{align*}
		\Big( \int_{B_R} d x \, \mathbb{E}^{\frac{2}{q}} \big| F_x * \psi_r ( x ) \big|^{q} \Big)^{\frac{1}{2}} & \lesssim r^{s} R^{| \beta | - \alpha} + r^{| \beta | - 2} R^{\frac{D}{2}} ,
	\end{align*}
where the new term $r^{| \beta | - 2} R^{\frac{D}{2}}$ comes from the contribution of \eqref{ap75}.
Note that $\lfloor s \rfloor \leq s \leq | \beta | - 2 < \lceil s \rceil$, so that $\lfloor s \rfloor = \lfloor | \beta | - 2 \rfloor$.
Thus, recalling \eqref{mb03}, also the second term is subject to the moment bound and we deduce by plugging into \eqref{mb02}
	\begin{align*}
		\Big( \int_{B_R} d x \, \mathbb{E}^{\frac{2}{q}} \big| U_x * \psi_r ( x ) \big|^{q} \Big)^{\frac{1}{2}} \lesssim r^{s+2} R^{| \beta | - \alpha} + r^{| \beta |} R^{\frac{D}{2}} ,
	\end{align*}
which absorbs into the desired $r^{2+s} ( r +R )^{| \beta | - \alpha}$.
This concludes the proof of \eqref{ap94}.
\end{proof}


\subsection*{Acknowledgements}
All authors thank Pavlos Tsatsoulis, who initiated the extension to the $\phi^4$ model
and substantially contributed to the core of the first draft, before leaving academics. 
We also thank Hector Bouton, Rhys Steele, Lorenzo Zambotti and the referee for their careful reading of this manuscript. 
FO acknowledges the hospitality of the Probability Summer School in Saint-Flour and the Bernoulli Center, 
where he presented parts of this material in the summer of 2023.
MT acknowledges funding by the Deutsche Forschungsgemeinschaft (DFG, German Research Foundation) under Germany’s Excellence Strategy EXC 2044--390685587, Mathematics Münster: Dynamics--Geometry--Structure, and from the European Research Council (ERC) under the European Union’s Horizon 2020 research and innovation programme (Grant agreement No.~101045082).

\bibliographystyle{alphaurl}
\bibliography{phi4}{}

\end{document}